\documentclass[12pt]{amsart}
\usepackage{latexsym, amssymb, amsfonts}
\newenvironment{pf}{\proof[\proofname]}{\endproof}
\theoremstyle{plain}
\newtheorem{Thm}{Theorem}[subsection]
\newtheorem{Cor}[Thm]{Corollary}
\newtheorem{Main}{Main Theorem}
\newtheorem{Mainf}[Thm]{First Main Theorem}
\newtheorem{Mains}[Thm]{Second Main Theorem}
\newtheorem{Maint}[Thm]{Third Main Theorem}
\newtheorem{Mainfo}[Thm]{Fourth Main Theorem}
\newtheorem{Prop}[Thm]{Proposition}
\newtheorem{Lem}[Thm]{Lemma}
\newtheorem{Cl}[Thm]{Claim}

\theoremstyle{definition}
\newtheorem{Def}[Thm]{Definition}
\newtheorem{Rem}[Thm]{Remark}

\newtheorem{Emp}[Thm]{}

\newtheorem{Con}[Thm]{Construction}
\newtheorem{Not}[Thm]{Notation}

\numberwithin{equation}{section}

\newcommand{\pgp}{\text{$\PGL{d}(K_w)$}}
\newcommand{\gp}{\text{$\GL{d}(K_w)$}}
\newcommand{\psp}{\text{$\PSL{d}(K_w)$}}
\newcommand{\ssp}{\text{$\SL{d}(K_w)$}}
\newcommand{\B}[1]{\Bbb#1}
\newcommand{\cal}[1]{\mathcal{#1}}
\newcommand{\C}[1]{\cal#1}
\newcommand{\pgr}{\text{$\PGU{d-1,1}(\B{R})$}}
\newcommand{\gr}{\text{$\GU{d-1,1}(\B{R})$}}

\newcommand{\sr}{\text{$\SU{d-1,1}(\B{R})$}}
\newcommand{\ff}[1]{\text {$\cal F(#1)$}}
\newcommand{\isom}{\overset {\thicksim}{\to}}
\newcommand{\om}[2]{\text{$\Omega^{#1}_{#2}$}}
\newcommand{\Om}{\om{d}{K_w}}
\newcommand{\si}[2]{\text{$\Sigma^{#1}_{#2}$}}
\newcommand{\Si}{\si{d}{K_w}}
\newcommand{\SiN}{\si{d,n}{K_w}}
\newcommand{\esc}[2]{\text{$(#1,#2)$-scheme}}
\newcommand{\lra}{\longrightarrow}
\newcommand{\hra}{\hookrightarrow}
\newcommand{\wt}{\widetilde}
\newcommand{\Gm}{\Gamma}
\newcommand{\gm}{\gamma}
\newcommand{\dt}{\delta}
\newcommand{\Dt}{\Delta}
\newcommand{\bs}{\backslash}

\newcommand{\GmG}{\Gm _{G}}

\newcommand{\gmG}{\gm _{G}}
\newcommand{\DtE}{\Dt _{E}}
\newcommand{\GmE}{\Gm _{E}}
\newcommand{\dtE}{\dt _{E}}
\newcommand{\gmE}{\gm _{E}}
\newcommand{\DtS}{\Dt _{S}}
\newcommand{\GmS}{\Gm _{S}}
\newcommand{\SIE}{\text{$S\in\ff{E}$}}

\newcommand{\m}{^{\times}}
\newcommand{\mDw}{\text{$D_{w}^{\times}$}}
\newcommand{\mODw}{\text{$\C{O}_{D_w}^{\times}$}}

\newcommand{\ii}{^{int}}
\newcommand{\op}{^{opp}}
\newcommand{\al}{\alpha}
\newcommand{\af}{\B{A}_{F}^{f}}
\newcommand{\afv}{\B{A}^{f;v}_{F}}
\newcommand{\la}{\lambda}
\newcommand{\an}{^{an}}
\newcommand{\tw}{^{tw}}
\newcommand{\rl}[1]{Lemma \ref{L:#1}}
\newcommand{\rn}[1]{Notation \ref{N:#1}}
\newcommand{\rcl}[1]{Claim \ref{C:#1}}
\newcommand{\rp}[1]{Proposition \ref{P:#1}}
\newcommand{\rr}[1]{Remark \ref{R:#1}}
\newcommand{\rc}[1]{Construction \ref{C:#1}}
\newcommand{\re}[1]{\ref{E:#1}}
\newcommand{\rco}[1]{Corollary \ref{C:#1}}
\newcommand{\rcr}[1]{Corollary \ref{C:#1}}
\newcommand{\rt}[1] {Theorem \ref{T:#1}}
\newcommand{\rd}[1]{Definition \ref{D:#1}}
\newcommand{\sm}{\smallsetminus}
\newcommand{\be}{\infty}
\newcommand{\emp}{\varnothing}
\newcommand{\GU}[1]{\bold{GU_{#1}}}

\newcommand{\SU}[1]{\bold{SU}_{#1}}
\newcommand{\PGU}[1]{\bold{PGU_{#1}}}
\newcommand{\GL}[1]{\bold{GL_{#1}}}
\newcommand{\PGL}[1]{\bold{PGL_{#1}}}
\newcommand{\SL}[1]{\bold{SL_{#1}}}
\newcommand{\PSL}[1]{\bold{PSL_{#1}}}
\newcommand{\Mat}{\text{Mat}}

\newcommand{\pr}{\text{pr}}
\newcommand{\Ker}{\text{Ker}}
\newcommand{\val}{\text{val}}
\newcommand{\Det}{\text{det}}
\newcommand{\im}{\text{Im}}
\newcommand{\Spec}{\text{Spec}}
\newcommand{\Aut}{\text{Aut}}

\newcommand{\Ad}{\text{Ad}}
\newcommand{\Gal}{\text{Gal}}
\newcommand{\Tr}{\text{Tr}}
\newcommand{\Fr}{\text{Fr}}

\newcommand{\PG}{\bold{PG}}
\newcommand{\PH}{\bold{PH}}
\newcommand{\Lie}{\text{Lie}}
\newcommand{\ad}{\text{ad}}

\newcommand{\Rep}{\text{Rep}}
\newcommand{\vect}{\text{Vec}}
\newcommand{\Tor}{\text{Tor}}
\newcommand{\End}{\text{End}}

\newcommand{\Hom}{\text{Hom}}
\newcommand{\Span}{\text{Span}}
\begin{document}


\title[$P$-adic uniformization]%
{$P$-adic uniformization\\
   of Unitary Shimura Varieties}
\author[Yakov Varshavsky]{Yakov Varshavsky}
\address{Department of Mathematics, University of Toronto\\ 
100 St. George St.\\ 
Toronto, ON, M5S 3G3\\ 
CANADA}
\email{vyakov@math.toronto.edu}
\date{\today}
\maketitle


\centerline{\bf Introduction}

Let $\Gm\subset\bold{PGU_{d-1,1}}(\B{R})^0$ be a torsion-free cocompact lattice. Then 
$\Gm$ acts on the unit ball $B^{d-1}\subset \B{C}^{d-1}$ by holomorphic
automorphisms. The quotient $\Gm\bs B^{d-1}$ is a complex manifold, which has a 
unique structure of a complex projective variety $X_{\Gm}$ 
(see \cite[Ch. IX, $\S$3]{Sha}). 

 Shimura had proved that when $\Gm$ is an arithmetic congruence subgroup, 
$X_{\Gm}$ has a canonical structure of a projective variety over some number field
$K$ (see \cite{De} or \cite{Mi}). For certain arithmetic problems it is
desirable to know a description of the reduction of $X_{\Gm}$ modulo $w$, 
where $w$ is some prime of $K$.
In some cases it happens that the projective variety $X_{\Gm}$ has a
$p$-adic uniformization. By this we mean that the $K_w$-analytic space 
$(X_{\Gm}\otimes_K K_w)\an$ is isomorphic to $\Dt\bs\Omega$ for some $p$-adic
analytic symmetric space $\Omega$ and some group $\Dt$, acting on $\Omega$  discretely.
Then a formal scheme structure on $\Dt\bs\Omega$ gives us an $\C{O}_{K_w}$-integral
model for $X_{\Gm}\otimes_K K_w$.

Cherednik was the first who obtained a result in this direction. Let $F$ be a 
totally real number field, and let $B/F$ be a quaternion algebra, which is definite at 
all infinite places, except one, and ramified at a finite prime $v$ of $F$.
Then Cherednik proved in \cite{Ch2} that the Shimura curve corresponding to $B$ 
has a $p$-adic uniformization by the $p$-adic upper half-plane $\om{2}{F_v}$, 
constructed by Mumford (see \cite{Mum1}), when the subgroup defining the
curve is maximal at $v$. Cherednik's proof is based on the method
of elliptic elements, developed by Ihara in \cite{Ih}.

The next significant step was done by Drinfel'd in \cite{Dr}. First he 
constructed certain covers of $\om{2}{F_v}$ (see below). Then, when $F=\B{Q}$, 
he proved the existence of a $p$-adic uniformization by some of his covers
for all Shimura curves, described in the previous paragraph, without the 
assumption of maximality at $v$. 
The basic idea of Drinfel'd's proof was to invent some moduli problem, whose 
solution is the Shimura curve as well as a certain $p$-adically uniformized curve,
showing, therefore, that they are isomorphic.

Developing Drinfel'd's method, Rapoport and Zink (see \cite {RZ,Ra}) obtained
some higher-dimensional generalizations of the above results.

In this paper we generalize Cherednik's method and prove that certain unitary 
Shimura varieties and automorphic vector bundles over them have a $p$-adic 
uniformization. Our results include all previously known results as particular
cases.

We now describe our work in more detail. Let $F$ be a totally real number field
of degree $g$ over $\B{Q}$, and let $K$ be a totally imaginary quadratic 
extension of $F$. Let $D$ and $D\ii$ be central simple algebras 
of dimension $d^2$ over $K$ with involutions of the second kind $\al$ and 
$\al\ii$ respectively over $F$. Let $\bold{ G}:=\GU{}(D,\al)$ and 
$\bold{ G\ii}:=\GU{}(D\ii,\al\ii)$ be the corresponding algebraic groups 
of unitary similitudes (see \rd{inv} and \rn{gu} for the notation).

Let $v$ be a non-archimedean prime of $F$ that splits in $K$, let $w$ and
$\bar{w}$ be the primes of $K$ that lie over $v$, and let $\be_1$ be an archimedean prime of $F$. Suppose that 
$D\ii\otimes_{K}{K_w}$ has Brauer invariant $1/d$, that  
$D\otimes_{K}K_w\cong \Mat_d(K_w)$, and that the pairs $(D,\al)\otimes_{F}F_u$ and 
$(D\ii,\al\ii)\otimes_{F}F_u$ are isomorphic for all primes $u$ of $F$, except $v$ and $\be_1$.
Assume also that $\al$ is positive definite at all archimedean places
$F_{\be_i}\cong\B{R}$ of $F$, that is that $\bold{ G}(F_{\be_i})\cong\GU{d}(\B{R})$ 
for all $i=1,...,g$, and that the signature of $\al\ii$ at $\be_1$ is $(d-1,1)$, so that 
$\bold{ G\ii}(F_{\be_1})\cong \GU{d-1,1}(\B{R})$.

Let $\af$ and $\afv$ be the ring of finite adeles of $F$
and the ring of finite adeles of $F$ without the $v$-th component respectively.
Set $E':=F_v\m\times \bold{G}(\afv)$, and fix a central simple algebra $\wt{D}_w$ over $K_w$ 
of dimension $d^2$ with Brauer invariant $1/d$. Then  $\bold{G\ii}(\af)\cong\wt{D}_w\m\times E'$
and $\bold{G}(\af)\cong\GL{d}(K_w)\times E'$. In particular, the
group $\GL{d}(K_w)$ acts naturally on $\bold{G}(\af)$ by left multiplication.

Let $\Om$ be the Drinfel'd's $(d-1)$-dimensional upper half-space over $K_w$ constructed 
 in \cite{Dr1}, and let  $\{\SiN\}_{n\in\B{N}\cup\{0\}}$ be the projective system of 
\'etale coverings of $\Om$ constructed in \cite{Dr}. This system is equipped with an equivariant  
action of the group $\GL{d}(K_w)\times\wt{D}\m_w$ such that if $T_n$ denotes 
the $n$-th congruence subgroup of $\C{O}_{\wt{D}_w}\m$, then
we have $T_n\bs\Sigma_{K_w}^{d,m}\cong\SiN$ for all $m\geq n$ (see \re{anspace} and \re{Dr} for our notation 
and conventions, which differ from those of Drinfel'd).

Denote by $\bold{ G\ii}(F)_+$ the set of all $d\in (D\ii)\m$ such that 
$d\cdot\al\ii(d)$ is a totally positive element of $F$. Choose an embedding 
$K\hra\B{C}$, extending $\be_1:F\hra\B{R}$. It defines us an embedding 
$\bold{ G\ii}(F)_{+}\hra \GU{d-1,1}(\B{R})^0=\Aut(B^{d-1})$. Choose finally an 
embedding of $K_w$ into $\B{C}$, extending that of $K$.

For each compact and open subgroup $S$ of $E'$ and each non-negative integer $n$
let $X_{S,n}$ be the weakly-canonical model over $K_w$ of the Shimura 
variety corresponding to the complex analytic space 
$(T_n\times S)\bs[B^{d-1}\times \bold{ G\ii}(\af)]/\bold{ G\ii}(F)_+$ and to the 
morphism $h:\bold{S}\to\bold{ G\ii}\otimes_{\B{Q}}\B{R}$, described in \re{Shim}
(see \rd{canmod} and \rr{canmod} for the definitions). The experts
might notice that our $h$ is not the one usually used in moduli problems of abelian varieties.

Let $V_{S,n}$ be the canonical model of the automorphic 
vector bundle on $X_{S,n}$ (see \cite[III]{Mi} or the last paragraph of the proof of 
\rp{eq} for the definitions), 
corresponding to the complex analytic space $(T_n\times S)\bs
[\beta_{\B{R}}^*(W\tw\otimes_{K_w}\B{C})\an\times \bold{ G\ii}(\af)]
/\bold{ G\ii}(F)_+$ (see \re{int} for the necessary notation).

Let $P_{S,n}$ be the canonical model of the standard principal
bundle over $X_{S,n}$ (see \cite[III]{Mi} or \rcr{canm} for the definitions),
corresponding to the complex analytic space $(T_n\times S)\bs
[B^{d-1}\times(\bold{PG\ii}\otimes_{\B{Q}}\B{C})\an\times\bold{G\ii}(\af)]
/\bold{G\ii}(F)_+$ (see \re{int} for the necessary notation).

\begin{Main} \label{M:0}
For each compact and open subgroup $S$ of $E'$ and each $n\in\B{N}\cup\{0\}$ we have
isomorphisms of $K_w$-analytic spaces:

a) $(X_{S,n})\an\isom\GL{d}(K_w)\bs[\SiN\times(S\bs\bold{G}(\af)/\bold{G}(F))]$;

b) $(V_{S,n})\an\isom\GL{d}(K_w)\bs[\beta_{w,n}^*(W\an)\times(S\bs\bold{G}(\af)
/\bold{G}(F))]$
(see \re{int} for the necessary notation), where the group  $\GL{d}(K_w)$ acts on
$\beta_{w,n}^*(W\an)$ as the direct factor of $(\bold{G}\otimes_{\B{Q}}K_w)(K_w)$, corresponding
to the natural embedding $K\hra K_w$;

c) $(P_{S,n})\an\isom [\GL{d}(K_w)\bs[\SiN\times((\bold{PG}\otimes_{\B{Q}}K_w)
\an\times(S\bs\bold{G}(\af)))/\bold{G}(F)]]\tw$ (see \re{int} for the definition
of the twisting $( )\tw$), where the group $\GL{d}(K_w)$ acts trivially on
$(\bold{PG}\otimes_{\B{Q}}K_w)\an$.

These isomorphisms commute with the natural projections for 
$S_1\subset S_2,\; n_1\geq n_2$ and with the action of $\bold{G\ii}(\af)\cong
\wt{D}_w\m\times F_v\m\times \bold{G}(\afv)$.
\end{Main}

The idea of the proof is the following. Consider the $p$-adic analytic varieties
$\wt{Y}_{S,n}$ of the right hand side of a) of the Main Theorem. They form
a projective system, and each of them has a natural structure $Y_{S,n}$ of
a projective variety over $K_w$. Kurihara proved in \cite{Ku}
that for every torsion-free cocompact lattice $\Gm\subset\PGL{d}(K_w)$ the Chern numbers 
of $\Gm\bs\Om$ are proportional to those
of the $(d-1)$-dimensional projective space and that the canonical class of 
$\Gm\bs\Om$ is ample. The result of Yau (see \cite{Ya}) then implies that 
$B^{d-1}$ is the universal covering of each connected component of the complex
analytic space $(Y_{S,n}\otimes_{K_w}\B{C})\an$ for all sufficiently small 
$S\in\ff{E}$ and all embeddings $K_w\hra\B{C}$.

It is technically better to work with the inverse limit of the $Y_{S,n}$'s equipped with the action
of the group $\bold{G\ii}(\af)\cong\wt{D}\m_w\times E'$ on it rather then to work with each 
$Y_{S,n}$ separately. Generalizing the ideas of Cherednik \cite{Ch2} we prove 
that there exists a subgroup 
$\Dt\subset \GU{d-1,1}(\B{R})^0\times \bold{ G\ii}(\af)$ such that $(Y_{S,n}
\otimes_{K_w}\B{C})\an\cong(T_n\times S)\bs(B^{d-1}\times\bold{G\ii}(\af))/\Dt$ 
for all compact open subgroups $S\subset E'$ and all 
$n\in\B{N}\cup\{0\}$.
 
Using Margulis' theorem on arithmeticity we show that the groups $\Dt$
and $\bold{ G\ii}(F)_+$ are almost isomorphic modulo centers. More precisely, we 
show that
$(Y_{S,n}\otimes_{K_w}\B{C})\an$ is isomorphic to a finite covering of 
$(X_{S,n}\otimes_{K_w}\B{C})\an$. Using Kottwitz' results \cite{Ko} on local 
Tamagawa measures we find that the volumes of $(Y_{S,n}\otimes_{K_w}\B{C})\an$ and
$(X_{S,n}\otimes_{K_w}\B{C})\an$ are equal. It follows that the varieties 
$Y_{S,n}\otimes_{K_w}\B{C}$ and $X_{S,n}\otimes_{K_w}\B{C}$ are isomorphic over
$\B{C}$. Comparing the action of the Galois group on the set of special points on
both sides we conclude that  $Y_{S,n}$ and $X_{S,n}$ are actually isomorphic over $K_w$.

Notice that if one considers only Shimura varieties corresponding to subgroups which are 
maximal at $w$, then the use of Drinfel'd's covers in the proof of the 
$p$-adic uniformization is very minor. (We use them only for showing that the 
$p$-adically uniformized Shimura varieties have Brauer invariant $1/d$ at $w$, 
that probably can be done directly.) In this case the proof would be technically much easier
but contain all the essential ideas.

The proof of the $p$-adic uniformization of standard principal bundles is similar.
In addition to the above considerations it uses the connection
on principal bundles. Using the ideas from \cite[III]{Mi} we show that  
the $p$-adic uniformization of standard principal bundles implies the $p$-adic 
uniformization of automorphic vector bundles. In fact Tannakian arguments show
(see \cite{DM}) that these statements are equivalent.

This paper is organized as follows. In the first section we introduce certain
constructions of projective systems of projective algebraic varieties, give
their basic properties and do other technical preliminaries.

In the second section we give two basic examples of such systems. Then we 
formulate and prove the complex version of our Main Theorem for Shimura varieties.

The third and the forth sections are devoted to the proof of the theorem on the 
$p$-adic uniformization of Shimura varieties and of automorphic vector bundles 
respectively.


Our proof appears to be very general. That is starting from any reasonable $p$-adic symmetric space,
whose quotient by an arithmetic cocompact subgroup is algebraizable, 
we find Shimura varieties uniformized by it. 
For example, in another work (\cite{Va}) we extend our results to Shimura varieties uniformized
by the product of Drinfel'd's upper half-spaces.
Hence it would be interesting to have more examples of such $p$-adic symmetric spaces.
 
Our result on the $p$-adic uniformization of automorphic vector bundles is not
complete, because we prove the $p$-adic uniformization only under the 
assumption that the center acts trivially. In fact our proof of the
complex version of the theorem works also in the general case, but to get an 
isomorphism over $K_w$ one should understand better the action of the Galois
group on the set of special points.

After this work was completed, it was pointed out to the author that Rapoport and Zink 
have recently obtained similar results concerning the uniformization of Shimura
varieties by completely different methods (see \cite{RZ2}).

\centerline{\bf Notation and conventions}
	
1) For a group $G$ let $Z(G)$ be the center of $G$, let $PG:=G/Z(G)$
be the adjoint group of $G$, and let $G^{der}$ be the derived group of $G$.

2) For a Lie group or an algebraic group $G$ let $G^0$ be its connected component of the identity.

3) For a totally disconnected topological group $E$ let $\ff{E}$ be the set of all compact and 
open subgroups of $E$, and let $E^{disc}$ be the group $E$ with the discrete
topology.

4) For a subgroup $\Gm$ of a group $G$  let $Comm_{G}(\Gm)$ be the commensurator of $\Gm$ in $G$.

5) For a subgroup $\Gm$ of a topological group $G$ let $\overline{\Gm}$ be the
closure of $\Gm$ in $G$.

6) For a set $X$ and a group $G$ acting on $X$ let $X^G$ be the set of
all elements of $X$ fixed by all $g\in G$.

7) For a set $X$, a subset $Y$ of $X$ and a group $G$ acting on $X$ let $Stab_{G}(Y)$
be the set of all elements of $G$ mapping $Y$ into itself.
 
8) For an analytic space or a scheme $X$ let $T(X)$ be the tangent bundle on $X$.

9) For a vector bundle $V$ on $X$ and a point $x\in X$ let $V_x$ be the fiber of $V$ over $x$.

10) For an algebra $D$ let $D\op$ be the opposite algebra of $D$.

11) For a finite dimensional central simple algebra $D$ over a field let $SD\m$
be the subgroup of $D\m$ consisting of elements with reduced norm $1$.
 
12) For a number field $F$ and a finite set $N$ of finite primes of $F$ let $\af$ be the ring of finite adeles of $F$, 
and let $\B{A}_F^{f;N}$ be the ring of finite adeles of $F$ without the components from $N$.

13) For a field extension $K/F$ let $\bold{R}_{K/F}$ be the functor of the restriction of scalars from $K$ to $F$.


14) For a natural number $n$ let $I_n$ be the $n\times n$ identity matrix, and let $B^n\subset\B{C}^n$ be 
the $n$-dimensional complex unit ball.
 
15) For a scheme $X$ over a field $K$ and a field extension $L$ of $K$ write $X_L$ or $X\otimes_K L$ instead of $X\times_{\Spec K}\Spec L$.
 
16) For an analytic space $X$ over a complete non-archimedean field $K$ and a for a  complete non-archimedean
field extension $L$ of $K$ let $X\widehat{\otimes}_K L$ be a field extension from $K$ to $L$. (A completion sign will be omitted in the case of a finite extension).

17) By a $p$-adic field we mean a finite field extension of $\B{Q}_p$ for some prime number $p$. Let $\B{C}_p$ be the completion of the algebraic closure of
$\B{Q}_p$.

18) By a $p$-adic analytic space we mean an analytic space over a $p$-adic field in the sense of Berkovich \cite{Be1}.
 
19) For an affinoid algebra $A$ let $\C{M}(A)$ be the affinoid space associated to it.

\centerline{\bf Acknowledgements}
 
First of all the author wants to thank Professor R. Livn\'e for formulating the problem, for 
suggesting the method of the proof and for his attention and help during all 
stages of the work. He also corrected an earlier version of this 
paper.

I am also grateful to Professor V. Berkovich for his help on $p$-adic analytic
spaces, to Professor J. Rogawski for the reference to \cite{Cl}, to Professor Th. Zink for 
his corrections and interest, and to the referee for his suggestions.  

The work forms part of the author's Ph.D. Thesis in the Hebrew University of 
Jerusalem, directed by Professor R. Livn\'e.

The revision of the paper was done while the author enjoyed the hospitality of the Institute for the Advanced Study at Princeton 
and was supported by the NSF grant DMS 9304580.

\tableofcontents

\section{Basic definitions and constructions} \label{S:bd}

\subsection{General preparations} \label{SS:denpre}
\begin{Def} \label{D:lpg}
 A {\em locally profinite group} is a locally compact
totally disconnected topological group. In such a group $E$,
the set $\ff{E}$ forms a fundamental system of neighbourhoods
of the identity element, and $\underset {S\in \ff{E}} {\bigcap} S=\{1\}$.
\end{Def}

\begin{Lem} \label{L:lpg}
Let $E$ be a locally profinite group, and let $X$ be a separated topological 
space with a continuous action $E \times X \to X$ of $E$. For each 
$S \in \ff{E}$, set
$X_S:=S\backslash X$. Then $\{X_S\}_S$ is a projective system, and 
$X \cong\underset{\underset {S}{\longleftarrow}}{\lim}\,X_S$.
\end{Lem}
\begin{pf}
\cite[Ch. II, Lem. 10.1]{Mi}
\end{pf}

This lemma motivates the following definition.

\begin{Def} \label{L:Esc}
Let $X$ be a separated scheme over a field $L$,
let $E$ be a locally profinite group, acting $L$-rationally on $X$.
We call $X$ an {\em \esc{E}{L}} (or simply an {\em $E$-scheme} if $L$ 
is clear or is not important) if for each $S\in \ff{E}$ there exists
a quotient $X_S:=S\backslash X$, which is a projective scheme over $L$,
and $X\cong\underset{\underset {S}{\longleftarrow}}{\lim}\,X_S$.
\end{Def}

The following remarks show that $E$-schemes are closely related to
projective systems of projective schemes, indexed by $\ff{E}$.

\begin{Rem}  \label{R:rem1}
If $X$ is an $E$-scheme or merely a topological space with a
continuous action of $E$, then for each $g\in E$ and each
$S,T \in \ff{E}$ with $S\supset gTg^{-1}$ we have a morphism
$\rho_{S,T}(g):X_T \to X_S$, induced by the action of $g$ on $X$ and
satisfying the following conditions:\par
a) $\rho_{S,S}(g)=$Id if $g \in S$;\par
b) $\rho_{S,T}(g)\circ \rho_{T,R}(h)=\rho_{S,R}(gh)$;\par
c) if $T$ is normal in $S$, then $\rho_{T,T}$ defines the action
 of the finite group $S/T$ on $X_T$, and $X_S$ is isomorphic to the 
quotient of $X_T$ by the the action of $S/T$.
\end{Rem}

\begin{Rem} \label{R:rem2} 
Conversely, suppose that for each $S \in\ff{E}$  there is given a
scheme $X_S$, and for each $g\in E$
and each $S,T \in \ff{E}$ with $S\supset gTg^{-1}$, there is given a morphism
$\rho_{S,T}(g):X_T \to X_S$, satisfying the conditions a)--c) of
\ref{R:rem1}. Then for each $T \subset S$ there is a map
$\rho_{S,T}(1):X_T \to X_S$, which is finite, by condition c).
In this way we  get a projective system of schemes and we can form an inverse 
limit scheme $X :=\underset{\underset {S}{\longleftarrow}}{\lim}\,X_S$.
Then there is a unique action of $E$ on $X$ such that for each $g\in E$ and 
each $S \in \ff{E}$ the action of $g$ on $X$
 induces an isomorphism $\rho_{gSg^{-1},S}:X_S \isom X_{gSg^{-1}}$. 
It follows from c) that
$X_S\isom S\backslash X$ for each $S \in \ff{E}$.
\end{Rem}

\begin{Def} \label{D:eqis}
Let $\wt{E}$ be a topological group, which is isomorphic to
$E$ under an isomorphism $\Phi :E \isom \wt{E}$. We say that an
 \esc{E}{L} $X$ is {\em $\Phi$-equivariantly isomorphic} to
an \esc{\wt{E}}{L} $\wt{X}$ if there exists an isomorphism
$\phi :X \isom \wt{X}$ of schemes over $L$ such that for each
$g \in E$ we have $\phi \circ g=\Phi(g) \circ \phi$.
If in addition $E=\wt{E}$ and $\phi$ is the identity, then we say
that $\phi$ is an {\em isomorphism of \esc{E}{L}s}.
\end{Def}


\begin{Def} \label{D:dec}
Let $L_2/L_1$ be a field extension. We say that  an \esc{E}{L_1} 
$X$ is  {\em an $L_2/L_1$-descent} of an \esc{E}{L_2} $Y$  if the \esc{E}{L_2}s
$X_{L_2}$ and $Y$ are isomorphic. 
\end{Def}

Suppose from now on that $E$ is a noncompact locally profinite group.

\begin{Not} \label{N:index}
For a topological group $G$ and a subgroup $\Gm \subset G\times E$ let
$\pr_{G}$ and $\pr_{E}$ be the projection maps from $\Gm $ to $G$ and $E$ 
respectively. 
 Set $\Gm _{G}:=\pr_{G}(\Gm ),\,\Gm_{E}:=\pr_{E}(\Gm )$ and  $\Gm _{S}:=\pr_{G}(\Gm\cap(G\times S))$
for each $S \in \ff{E}$. For each $\gm \in \Gm $ set $\gm_{G}:=\pr_{G}(\gm )$ and $\gm_{E}:=\pr_{E}(\gm )$.
\end{Not}

\begin{Lem} \label{L:basic}
Let  $\Gm \subset G\times E$ be a cocompact lattice. Suppose that $\pr_G$ is
injective. Then for each $S \in \ff{E}$ we have the following: \par
a) $|S\bs E/\Gm_E|<\infty$; \par
b) $[\Gm_{G}:\Gm_{S}]=\infty$; \par
c) $\Gm_{S}$ is a cocompact lattice of $G$; \par
d) $\Gm_{G} \subset Comm_{G}(\Gm_{S})$.
\end{Lem}

\begin{pf} 
a) Since the the double quotient
$(G\times S)\bs(G\times E)/\Gm\cong S\bs E/\Gm_{E}$ is compact and discrete, it is finite.
 
b) The group $E$ is noncompact, therefore $|S\bs E|=\infty$. Hence 
$[\Gm_{E}:S \cap \Gm_{E}]=|S\bs S\Gm_{E}|=\be$. But $\GmG =\pr_{G}(\Gm)=
\pr_{G}(\pr_{E}^{-1}(\GmE ))$, and likewise $\GmS =\pr_{G}(\pr_{E}^{-1}(\GmE \cap S))$.
Since  $\pr_{G}$ is injective, we are done.

c) The group $\Gm$ is a cocompact lattice in $G\times E$, hence $\Gm\cap (G\times S)$
is a cocompact lattice in $G \times S$, and the statement follows by projecting to $G$
(see \cite[Prop. 1.10]{Shi}).\par

d) Let $\gm\in\Gm$, and set $S'=\gmE S \gmE ^{-1} \in \ff{E}$. Then 
$\gm (\Gm\cap(G\times S))\gm^{-1}=\Gm\cap(G\times S')$. 
But $S \cap S' \in \ff{E}$ is a subgroup of finite index in both $S$ and $S'$,
hence $\gmG \GmS \gmG ^{-1} \cap \GmS=\Gm_{S \cap S'}$ is a subgroup of finite 
index in both $\GmS$ and $\gmG \GmS \gmG ^{-1}$.
\end{pf}

Suppose that $d\geq 2$, and take $G$ equal to $\pgp $ for some $p$-adic field $K_w$ 
or to $\pgr^0$. We shall call these the $p$-adic and
the real (or the complex) cases respectively. 
 
\begin{Prop} \label{P:basic}
Under the assumptions of Lemma \ref{L:basic} we have: \par
a) $\overline{\GmG} \supset G^{\text{der}} $; \par
b) $\pr_{E}$ is injective; \par
c) for each $S \in \ff{E}$, the group $\GmS $ is an arithmetic subgroup of $G$ in the
  sense of Margulis (see \cite[pp.292]{Ma});\par
d) if $S \in \ff{E}$ is sufficiently small, then the subgroup $\Gm_{aSa^{-1}}$
is torsion-free for each $a\in E$.
\end{Prop}

\begin{pf}

a) For each $S \in \ff{E},\;\GmS $ is cocompact in $G$ and $[\GmG :\GmS ]=
\infty$. It follows that $\overline{\GmG}$ is a closed non-discrete cocompact subgroup
of $G$. Therefore its inverse image
$\pi ^{-1}(\overline{\GmG})$ in $\sr $ (resp. $\ssp $) is also closed, non-discrete and
cocompact, hence by \cite[Ch. II, Thm. 5.1]{Ma} it is all of $\sr $ (resp. $\ssp $).
This completes the proof.\par

b) Set $\Gm _0 :=\pr_{G}(\Ker\,\pr_{E})$. This is a discrete (hence a closed) subgroup  
of $G$, which is normal in $\GmG $. Therefore it is normal in 
$\overline{\GmG} \supset G^{\text{der}}$. It follows that each $\gm\in\Gm _0$ must commute with some open 
neighborhood of the identity in $G^{\text{der}}$, hence $\Gm_0$ is trivial. \par

c) is a direct corollary of \cite[Ch. IX, Thm. 1.14]{Ma} by b)--d) of
Lemma \ref{L:basic}. \par

d) (compare the proof of \cite[Lem. 1.3]{Ch1})
Choose an $\SIE$, then $\GmS\subset G$ is a cocompact lattice.

\begin{Lem} \label{L:fin}
 The torsion elements of $\GmS$ comprise a finite number of 
conjugacy classes in $\GmS$.
\end{Lem} 

We first complete the proof of the proposition assuming the lemma. Let $a_1,...,a
_n\in  E$ be representatives of double classes $\GmE\bs E/S$ (use \rl{basic}, a)).
For each $i=1,...,n$ let $M_i\subset \Gm_{a_i S a_i^{-1}}$ be a finite set of 
representatives of conjugacy classes of torsion non-trivial elements 
of $\Gm_{a_i S a_i^{-1}}$. Then the image of all non-trivial
torsion elements of $\Gm_{a_i S a_i^{-1}}$ under the natural injection 
$j_i:\Gm_{a_i S a_i^{-1}}\isom\Gm\cap (G\times a_i S a_i^{-1})\hra a_i S a_i^{-1}
\isom S$ is contained
in the set $X_i=\{s\cdot j_i(M_i)\cdot s^{-1}|s\in S\}$, which is compact and does
 not contain $1$. Hence there 
exists $T\in\ff{E}$ not intersecting any of the $X_i$'s. By taking a smaller 
subgroup we may suppose that $T$ is a normal subgroup of $S$. Since  all the $j_i$'s
are injective, the subgroup $\Gm_{a_i T a_i^{-1}}=j_i^{-1}(T)$ is 
torsion-free for each $i=1,...,n$. For each $a\in E$ there exist $i\in\{1,...,n\},
\,s\in S$ and $\gm\in\Gm$ such that $a=\gmE a_i s$. Hence the
subgroup
\begin{align}
\Gm_{aTa^{-1}} & \cong\Gm\cap(G\times aTa^{-1})=\Gm\cap(\gmG G\gmG^{-1}\times 
\gmE a_i T a_i^{-1}\gmE^{-1})\notag\\
  & \cong\gm(\Gm\cap(G\times a_iTa_i^{-1}))\gm^{-1}\cong\Gm_{a_i T a_i^{-1}}\notag 
\end{align}
is torsion-free.
\end{pf}

\begin{pf} 
(of the lemma)
The group $G$ acts continuously and isometrically
on some complete negative curved metric space $Y$. Indeed, in
the real case $Y=B^{d-1}$ with the hyperbolic metric. In the $p$-adic case
$Y$ is a geometric realization (see \cite[Ch. I, appendix]{Br}) of the 
Bruhat-Tits building $\Dt$ of $\SL{d}(K_w)$. This is a locally finite 
simplicial complex of dimension $d-1$ which can be described as follows. 
Its vertices are the equivalence classes of free $\C{O}_{K_w}$-submodules of rank 
$d$ of the vector space $K_w^n$, where $M$ and $N$ are equivalent when there exists $a\in K_w\m$
such that $M=aN$. The distinct vertices $\Dt_1,\Dt_2,...,\Dt_k$ form a simplex 
when there exist for them representative lattices $M_1,M_2,...,M_k$ such that 
$M_1\supset M_2\supset  ...\supset M_k\supset \pi M_1$. For more information 
see \cite[$\S$1]{Mus} or \cite[Ch. V, $\S$8]{Br}. 

The geometric realization $Y$ of $\Dt$ has a canonical metric, that makes $Y$ 
a complete metric space with negative curvature (see \cite[Ch. VI, $\S$3]{Br}).
Moreover, the natural action of $\pgp$ on the set of vertices of $\Dt$ can be
(uniquely) extended to the simplicial, continuous and isometric action on $Y$.

Now the Bruhat-Tits fixed point theorem (see \cite[Ch. VI, $\S$4, Thm. 1]{Br}) 
implies that each compact subgroup of $G$ has a fixed point on $Y$. In 
particular, each torsion element of $G$ has a fixed point on $Y$. 
Notice that in the $p$-adic case it then
stabilizes the minimal simplex, containing the fixed point.

Conversely, the stabilizer in $G$ of each point of $Y$ is compact. In the 
real case this is true, since the group $\pgr^0$ acts transitively on $B^{d-1}$ and 
the group $K=Stab_{B^{d-1}}(0)\cong \bold{U_{d-1}}(\B{R})$ is compact. 
In the $p$-adic
case the group $\pgp$ acts transitively on the set of vertices,
and the stabilizer of the equivalence class of $\C{O}_{K_w}^d\subset K_w^d$ is 
$\PGL{d}(\C{O}_{K_w})$, hence it is compact. 
Since the stabilizer in $G$ of each point $y\in Y$ must stabilize the minimal
simplex $\sigma$ containing $y$, it must permute the finitely many of 
vertices of $\sigma$, so that it is also compact.
It follows that the stabilizer of each point of $Y$ in $\GmS$ is compact and discrete, 
hence it is finite.

To finish the proof of the lemma in the real case we note that for each 
$x\in B^{d-1}$ there exists an open neighbourhood $U_x$ of $x$ such that 
$\Gm_x:=\{g\in\GmS|
g(U_x)\cap U_x\neq \o\}=\{g\in\GmS|g(x)=x\}$ is finite (see \cite[Prop. 1.6 and
1.7]{Shi}).
The space $\GmS\bs B^{d-1}$ is compact, hence there exist a finite number of 
points $x_1,x_2,...,x_m$ of $B^{d-1}$ such that 
$\GmS(\bigcup_{i=1}^{m}U_{x_i})=B^{d-1}$. If $\gamma$ is a torsion element of 
$\GmS$, then it fixes some point of  $B^{d-1}$. By conjugation we may assume 
that it fixes a point in some $U_{x_i}$, therefore $\gamma$ is conjugate to
an element of the finite set $\bigcup_{i=1}^{m}\Gm_{x_i}$.

In the $p$-adic case we first assert that $\Dt$ has only a finite number 
of equivalence classes of simplexes under the action of $\GmS$. Since $\Dt$ is
locally finite, it is enough to prove this assertion for vertices. The 
group  $G$ acts transitively on the set of vertices, and $G=\GmS\cdot K$ for 
some compact set $K\subset G$. Hence if $v$ is a vertex of $\Dt$, then 
$K\cdot v$ is a compact and discrete (because the set of all vertices of $\Dt$ 
 is a discrete set in $Y$) subset of $Y$, and our assertion follows. Now the 
same considerations as in the real case complete the proof.
\end{pf}

\subsection{GAGA results} \label{SS:GAGA}

In what follows we will need some GAGA results. Let $L$ be equal to $K_w$ in the $p$-adic case and to $\B{C}$ in the complex case. 
We will call both the complex and the $p$-adic ($L$-)analytic spaces simply ($L$-)analytic spaces.
Recall that for each scheme $X$ of locally finite type over $L$ and each
coherent sheaf $F$ on $X$ a certain $L$-analytic space $X\an$ and a
coherent analytic sheaf $F\an$ on $X\an$ can be associated (see 
\cite[Thm. 3.4.1]{Be1} in the $p$-adic case and \cite[Exp. XII]{SGA1} in the 
complex one).

\begin{Thm} \label{T:GAGA}
Let $X$ be a projective $L$-scheme. The functor $F\mapsto F\an$ from the
category of coherent sheaves on $X$ to the category of coherent analytic 
sheaves
on $X\an$  is an equivalence of categories.
\end{Thm}

\begin{pf}
 In the complex case the theorem is proved in \cite[$\S$12, Thm. 2 and 3]{Se4}, in
the $p$-adic one the proof is the same. One first shows by a direct computation
that the $p$-adic analytic and the algebraic cohomology groups of $\B{P}^n$
coincide. Next one concludes from Kiehl's theorem (see \cite[Prop. 3.3.5]{Be1})
that the cohomology group of an analytic coherent sheaf on $\B{P}^n$
is a finite dimensional vector space. Now the arguments of Serre's proof 
in the complex case hold in the $p$-adic case as well. 
See \cite[3.4]{Be1} for the relevant definitions and basic properties.
\end{pf}
 
\begin{Cor} \label{C:GAGA1}
a) If $X$ is an algebraic variety over $L$ and $X'$ is a compact $L$-analytic
subvariety of $X$, then $X'$ is a proper $L$-algebraic subvariety of $X$.\par
b) The functor which associates to a proper $L$-scheme $X$
the analytic space $X^{an}$ is fully faithful.
\end{Cor}

\begin{pf} 
Serre's arguments (see \cite[\S19, Prop. 14 and 15]{Se4}) hold in both the complex and
the $p$-adic cases.
\end{pf}

\begin{Cor} \label{C:GAGA2}
 Let $X$ be a projective $L$-scheme. The functor $X'\mapsto (X')^{an}$ 
induces an equivalence between:

a) the category of vector bundles of finite rank on $X$ and 
the category of analytic vector bundles of finite rank on $X\an$;

b) the category of finite schemes over $X$ and the category
of finite $L$-analytic spaces over $X^{an}$, if $L$ is a $p$-adic field.
\end{Cor}

\begin{pf}
a) To prove the statement we first notice that the category of vector 
bundles of finite rank is equivalent to the category of locally free sheaves
of finite rank. In the algebraic case this is proved in \cite[II, Ex. 5.18]{Ha}.
In the analytic case the proof is similar. Now the corollary would follow from
the theorem, if we show that locally free analytic sheaves of finite rank
correspond to locally free algebraic ones. The analytic structure sheaf is 
faithfully flat over the algebraic one (see \cite[\S 2, Prop. 3]{Se4} and
\cite[Thm. 3.4.1]{Be1}). Therefore the statement follows from the fact that
an algebraic flat coherent sheaf is locally free (see \cite[Thm. 2.9]{Mi2}).

%
%

b) We first show that the correspondence $(\varphi:Y\to X)\mapsto
\varphi_{*}(\C{O}_{Y})$ (resp. $(\varphi:\wt{Y}\to X\an)\mapsto
\varphi_{*}(\C{O}_{\wt{Y}})$) gives an equivalence between the category
of finite schemes (resp. analytic spaces) over $X$ (resp. $X\an$)
and the category of coherent $\C{O}_{X}-$ (resp. $\C{O}_{X\an}-$)algebras. 
In the algebraic case this is proved in \cite[II, Ex. 5.17]{Ha}. In the analytic
case the proof is exactly the same, because a finite algebra over an
affinoid algebra has a canonical structure of an affinoid algebra
(see \cite[Prop. 2.1.12]{Be1}).
\end{pf}

\begin{Rem}
 If $X'$ is finite over $X$, then it is projective over $X$, therefore if, 
in addition, $X$
is projective over $K_w$, then $X'$ is also projective over $K_w$. 
\end{Rem}

\begin{Cor} \label{C:GAGA3}
Let $X$ and $Y$ be projective $L$-schemes, and let $W$ and $V$ be algebraic vector 
bundles 
of finite ranks on $X$ and $Y$ respectively. Then for each analytic map of 
vector bundles $\wt{f}:W\an\to V\an$ covering some map $f:X\to Y$ there exists a unique algebraic 
morphism $g:W\to V$ such that $g\an =\wt{f}$.
\end{Cor}	

\begin{pf}
By definition, $\wt{f}$ factors uniquely as 
$W\an\overset{\wt{f}'}{\lra}V\an\times
_{Y\an}X\an\cong (V\times_{Y} X)\an\overset{\text{proj}}{\lra}V\an$. \rco{GAGA2} implies that there 
exists a unique $g':W\to V\times_Y X$ such that $(g')\an=\wt{g}'$. 
Set $g:=\text{proj}\circ g'$. 

For the uniqueness observe that if $h:W\to V$ 
satisfies $h\an=\wt{g}$, then it covers $f$. Hence $h$ factors as
$W\overset{h'}{\lra}V\times_Y X\overset{\text{proj}}{\lra}V$. 
Since $\wt{f}'$ and $g'$ are unique, we have $h'=g'$ and $h=g$.
\end{pf}

\begin{Rem}
Using the results and ideas of \cite[Exp. XII]{SGA1} one can replace in the above 
results the assumption of projectivity by properness.
\end{Rem}
 
We now introduce two constructions of  $E$-schemes which are basic for this
work.

\subsection{First construction} \label{SS:fircon}
   
\begin{Emp} \label{E:anspace}
Let $\om{d}{K_w}$ be an open $K_w$-analytic subset of $(\B{P}_{K_w}^{d-1})\an$, obtained by removing from $(\B{P}_{K_w}^{d-1})\an$ the
union of all the  $K_w$-rational hyperplanes (see \cite{Be1} and \cite{Be3} for the definition and basic 
properties of analytic spaces). It is called the $(d-1)$-dimensional  Drinfel'd upper half-space over $K_w$
(see also \cite[\S 6]{Dr1}). Then
 $\om{d}{K_w}$ is the generic fiber of a certain formal scheme $\Hat{\Omega}^{d}_{K_w}$ over $\cal O_{K_w}$, 
constructed in \cite{Mus,Ku}, generalizing \cite{Mum1}. 

The group $\pgp$ acts naturally on $\Om$. (It will be convenient for us to consider $\B{P}^{d-1}$ as the set of
lines in $\B{A}^d$ and not as the set of hyperplanes, as Drinfel'd does. Therefore our action differs by transpose 
inverse  from that of Drinfel'd.) Moreover, this action naturally extends to the $\cal O_{K_w}$-linear action 
of $\pgp$ on $\Hat{\Omega}^{d}_{K_w}$. Furthermore, $\pgp$ is the group of all formal scheme automorphisms of 
$\Hat{\Omega}^{d}_{K_w}$ over $\cal O_{K_w}$ (see \cite[Prop. 4.2]{Mus}) and of all analytic 
automorphisms of $\Om$ over $K_w$ (see \cite{Be2}). Though the action of  $\pgp$ on $\Om$ is far 
from being transitive, we have the following

\begin{Lem} \label{L:inv}
There is no non-trivial closed analytic subspace of 
$\Om\Hat{\otimes}_{K_w}\B{C}_p$, invariant under the subgroup
$$
U=\left\{U_{\bar{x}}:=\left( \text{
\begin{tabular}{c|c}
$I_{d-1}$ & $\bar{x}$\\ \hline
$0$       & $1$
\end{tabular}}
\right) \Bigg|\: \bar{x}\in K_w^{d-1}\right\}\subset\PGL{d}(K_w).
$$
\end{Lem}

\begin{pf}
Suppose that our lemma is false. Let $Y$ be a non-trivial $U$-invariant
closed analytic subset of $\Om\Hat{\otimes}_{K_w}\B{C}_p$. Then 
dim$\,Y<$dim$\,\Om\Hat{\otimes}_{K_w}\B{C}_p=d-1$. Choose a regular point 
$y\in Y(\B{C}_p)$ (the set of regular points is open and non-empty). Then 
dim$\,T_y(Y)=$dim$\,Y<d-1$. Next we identify $\Om\Hat{\otimes}_{K_w}\B{C}_p$ with an open
analytic subset of $(\B{A}^{d-1}_{\B{C}_p})\an$ by the map $(z_1:...:z_d)\mapsto (\frac{z_1}{z_d},...,
\frac{z_{d-1}}{z_d})$. Then $U_{\bar{x}}(z)=z+\bar{x}$ for every $z\in\B{A}^{d-1}_{\B{C}_p}$ and every 
$\bar{x}\in\B{A}^{d-1}(K_w)$. In particular, $y+\bar{x}\in Y$ for every 
$\bar{x}\in\B{A}^{d-1}(K_w)$, contradicting the assumption that dim$\,T_y(Y)<d-1$.
\end{pf}

Recall also that the group $\pgr^0$ acts transitively on $B^{d-1}$ and that it is the group of all 
analytic (holomorphic) automorphisms of  $B^{d-1}$
(see \cite[Thm. 2.1.3 and 2.2.2]{Ru}).
\end{Emp}

In what follows we will need the notion of a pro-analytic space.

\begin{Def} \label{L:proan}
A {\em pro-analytic space} is a projective system $\{X_{\al}\}_{\al\in I}$ of 
analytic spaces such that for some $\al_0\in I$ all transition maps 
$\varphi_{\beta\al}:X_{\beta}\to X_{\al};\;\beta\geq\al\geq\al_0$ are \'etale and surjective.
\end{Def}

\begin{Def} \label{L:point}
By a {\em point} of $X:=\{X_{\al}\}_{\al\in I}$ we mean a system 
$\{x_{\al}\}_{\al\in I}$, where $x_{\al}$ is a point of $X_{\al}$ for all 
$\al\in I$, and $\varphi_{\beta\al}(x_{\beta})=x_{\al}$ for all $\beta\geq\al$ in 
$I$. For a point $x=\{x_{\al}\}_{\al\in I}$ of $X=\{X_{\al}\}_{\al\in I}$
let $T_x(X):=\{v=\{v_{\al}\}_{\al\in I}|\,v_{\al}\in T_{x_{\al}}(X_{\al}),\,
d\varphi_{\al\beta}(v_{\beta})=v_{\al}  
\text{ for all }\beta\geq\al \text{ in }I\}$ be the {\em tangent space} of $x$ in 
$X$.
\end{Def}

\begin{Def} \label{L:morph}
Let $X=\{X_{\al}\}_{\al\in I}$ and $Y=\{Y_{\beta}\}_{\beta\in J}$ be two
pro-analytic spaces. To give a {\em pro-analytic morphism} $f:X\to Y$ is to give
an order-preserving map $\sigma:I\to J$, whose image is cofinal in $J$, and a 
projective system of 
analytic morphisms $f_{\al}:X_{\al}\to Y_{\sigma(\al)}$. A morphism $f$ is called
{\em \'etale} if there exists $\al_0\in I$ such that for each $\al\geq\al_0$
the morphism $f_{\al}$ is \'etale.
\end{Def}

\begin{Con} \label{C:con1}
Suppose that $\Gm \subset G\times E$ satisfies the conditions of Lemma
\ref{L:basic}. We are going to associate to $\Gm$ a certain $(E,L)$-scheme. 

Let  $X^0$ be $B^{d-1}$ in the real case and $\Om$ in the 
$p$-adic one. Consider the $L$-analytic space $\wt{X}:=(X^0\times E^{disc})/\Gm$,
where $\Gm $ acts on $X^0 \times E^{disc}$ by the natural right action:
$(x,g)\gm:=(\gmG ^{-1}x,g\gmE )$. Then $E$ acts analytically on $\wt{X}$ 
by left multiplication.

\begin{Prop} \label{P:sav1}
For each $S\in\ff{E}$ the quotient $S\bs \wt{X}=S\bs(X^0\times E)/\Gm$ exists
and has a natural structure of a projective scheme $X_S$ over $L$.
\end{Prop}

\begin{pf}
 First take $S\in \ff{E}$ satisfying part d) of Proposition
\ref{P:basic}. Then $S\bs \wt{X}$ has $|S\bs E/\GmE|<\infty$ connected 
components, each of them is isomorphic to $\Gm_{aSa^{-1}}\bs X^0$ for 
some $a\in E$. By c), d) of \rp{basic}, each $\Gm_{aSa^{-1}}$ is a 
torsion-free arithmetic cocompact lattice of $G$. 

By \cite[Prop. 1.6 and 1.7]{Shi}, \cite[Ch. IX, 3.2]{Sha} in the real case and by
\cite{Mus} or \cite{Ku} in 
the $p$-adic one, each quotient $\Gm_{aSa^{-1}}\bs X^0$ exists and has a unique structure of a projective 
algebraic variety over $L$. Therefore there exists 
a projective scheme $X_S$ over $L$ such that $X_S^{an}\cong S\bs \wt{X}$.

 Take now an arbitrary $S\in\ff{E}$. It has a normal subgroup $T\in\ff{E}$
which satisfies part d) of Proposition \ref{P:basic}.
The finite group $S/T$ acts on $T\bs\wt{X}\cong X_T^{an}$ by analytic
automorphisms, and $S\bs\wt{X}\cong (S/T)\bs X_T^{an}$.
\rco{GAGA1} implies that the analytic action of $S/T$ on
 $X_T^{an}$ defines an algebraic action on $X_T$ and that the 
projective scheme $X_S:=(S/T)\bs X_T$ (the quotient exists by
\cite[\S 7]{Mum2}) satisfies $(X_S)^{an}\cong S\bs\wt{X}$. Moreover,
the same corollary implies also that the algebraic structure on $S\bs\wt{X}$
is unique.
\end{pf}

 For all $g\in E$ and all $S,T\in\ff{E}$ with $S\supset gTg^{-1}$
we obtain by \rr{rem1} analytic morphisms $\rho_{S,T}(g): X_T^{an}\to X_S^{an}$.
They give us by \rco{GAGA1} uniquely determined algebraic morphisms $\rho_{S,T}(g):X_T\to X_S$, 
which provide us by Remark \ref{R:rem2} an 
\esc{E}{L} $X:=\underset{\underset{S}{\longleftarrow}}{\lim}\, X_S$.
\end{Con}


\begin{Prop} \label{P:prop1}
a) There exists the inverse limit $X\an$ of the $X\an_S$'s in the
category of $L$-analytic spaces, which is isomorphic to $\wt{X}$.

b) $Stab_E(X^0\times\{1\})=\GmE$.

c) Let $X_0$ be the connected component of $X$ such that $X_0^{an}\supset
X^0\times\{1\}$ (note that $X^0\times\{1\}$ is a connected component of $X^{an}$,
and that the analytic topology is stronger then the Zariski topology). Then 
$Stab_E(X_0)=\overline{\GmE}$.

d) The group $E$ acts faithfully on $X$.

e) For each $x\in X$ the orbit $E\cdot x$ is (geometrically) Zariski dense. In particular, 
$E$ acts transitively on the set of geometrically connected components of $X$.

f) For each $S\in \ff{E}$ satisfying part d) of \ref{P:basic},
the map $X\to X_S$ is \'etale; 

g) For each embedding $K_w\hra \B{C}$ and each $S\in\ff{E}$ as in c), $B^{d-1}$ 
is the universal covering of each connected component of $(X_{S,\B{C}})^{an}$
 in the $p$-adic case and of $X_{S}^{an}$ in the complex one.
\end{Prop}

\begin{pf}
a) We start from the following 
\begin{Lem} \label{L:discr}
a) Let $\Pi$ be a torsion-free discrete subgroup of $G$. Then the natural projection 
$X^0\to \Pi\bs X^0$ is an analytic (topological) covering.

b) For each $x\in X^0$ the stabilizer of $x$ in $G$ is compact.
\end{Lem}

\begin{pf}
a) follows from \cite[Prop. 1.6 and 1.7]{Shi} in the real case and from 
\cite[Lem. 4 and 6]{Be2} in the $p$-adic one.

b) By \cite[$\S$6]{Dr} there exists a $\PGL{d}(K_w)$-equivariant map from
$\Om$  to the Bruhat-Tits building $\Dt$ of $\SL{d}(K_w)$, thus it suffice 
to show the required property for stabilizers of points in $\Dt$ and $B^{d-1}$. 
This was done in the proof
of \rl{fin}.
\end{pf}

The lemma implies that for each sufficiently small $\SIE$
the analytic space $X_S\an$
admits a covering by open analytic subsets $U_i$ satisfying the following condition: 
for each $i$ and each subgroup 
$S\supset T\in\ff{E}$ the inverse image 
$\rho_T^{-1}(U_i)$ of $U_i$ under the natural projection $\rho_T:X_T\an\to X_S\an$ 
splits as a disjoint union 
of analytic spaces, each of them isomorphic to $U_i$ under $\rho_T$.

Now we will define a certain $L$-analytic space $X\an$ associated to $X$.
As a set it is the inverse limit of the underlying sets of the $X_S\an$'s.
To define an analytic structure on $X\an$ 
consider subsets  $V_{\al}\subset X\an$ such that for some (hence for every) sufficiently small 
$\SIE$, the natural projection 
$\pi_S:X\an\to X_S\an$ induces a bijection of  $V_{\al}$ with an open analytic subset 
$\pi_S(V_{\al})$ of $X_S\an$, described in the previous paragraph. Provide then such a  
$V_{\al}$ with an analytic structure by requiring that $\pi_S:V_{\al}\to
\pi_S(V_{\al})$ is an analytic isomorphism. Then the analytic structure of the $V_{\al}$'s does not depend 
on the choice of the $S$'s, and there exists
a unique $L$-analytic structure on $X\an$ such that each $V_{\al}$ is an open
analytic subset of $X\an$. 

By the construction,  $X\an$ is the inverse limit of the $X_S\an$'s in the category of
$L$-analytic spaces. Hence there exists a unique
$E$-equivariant analytic map $\pi:\wt{X}\to X\an$
such that for each $\SIE$ the natural projection $\wt{X}\to X_S\an$ factors as 
$\wt{X}\overset{\pi}{\lra}X\an\overset{\pi_S}{\lra}X_S\an$, where by $\pi_S$ 
we denote the natural projection. It remains to show
that $\pi$ is an isomorphism.

For each $\SIE$ satisfying part d) of \rp{basic} the natural projection
$X^0\to\GmS\bs X^0$ is a local isomorphism. Hence the projections $\wt{X}\to X_S$
and $\pi$ are local isomorphisms as well.
 
The map $\pi_S\circ \pi$ is surjective, hence for each
$x\in X\an$ there exists a point $y\in\wt{X}$ such that 
$\pi_S(x)=\pi_S\circ \pi(y)$.
Therefore, $\pi(y)=sx$ for some $s\in S$. Since $\pi$ is $E$-equivariant,
we conclude that $\pi(s^{-1}(y))=x$. Hence $\pi$ is surjective.

Suppose that $\pi(y_1)=\pi(y_2)$ for some $y_1,\,y_2\in\wt{X}$. 
Let $(x_1,g_1)$ and $(x_2,g_2)$ be their representatives in $X^0\times E$.
Then for each $\SIE$ there exist $s\in S$ and $\gm\in\Gm$ such that $x_1=\gm_G^{-1}(x_2)$
and $g_1=sg_2\gmE$. Such $\gmG$'s belong to the set $\{g\in G | g(x_1)=x_2\}\cap\Gm_{g_2^{-1}Sg_1}$,
which is compact (by the lemma) and discrete, hence finite. Therefore we can choose sufficiently small $\SIE$ 
such that $g_1\gm_E^{-1}g_2^{-1}=s\in S$ must be equal to $1$. This means that $y_1=y_2$.
Thus $\pi$ is a surjective, one-to-one local isomorphism, hence
it is an isomorphism. 

b) is clear.\par

c) For each $\SIE$ let $Y_S$ be the connected component of $X_S$ 
such that $Y_S\an$ is the image of $X^0\times\{1\}\subset X\an$
under the natural projection $\pi_S:X\an\to X_S\an$. Then 
$X_0=\underset{\underset{S}{\longleftarrow}}{\lim}\, Y_S$. It follows that $g\in E$ satisfies
$g(X_0)=X_0 $ if and only if $ g(Y_S)=Y_S$ for each $\SIE$ if and only if $ X^0\times\{g\}\subset 
S(X^0\times 1)\Gm$ 
for each $S\in\ff{E} $ if and only if $ g\in S\GmE$ for each $\SIE$ if and only if $ g\in\underset{\SIE}
{\bigcap}S\GmE=\overline{\GmE}$. 

d) If $g\in E$ acts trivially on $X$, then it acts trivially on 
$X^{an}\cong(X^0\times E^{disc})/\Gm$. By b), $g=\gmE$ for some $\gm\in\Gm$, and $\gmG$ acts trivially on $X^0$. Since
$\pr_{G}$ is injective, $\gm =g=1$.\par

e) Let $Y$ be the Zariski closure of $E\cdot x$. Then $Y$ is $E$-invariant, and, therefore, 
$Y\an\cap(X^0\times\{1\})$ is a closed $\Gm$-invariant analytic subspace of $X^0\times\{1\}\cong X^0$.  
By \rp{basic} a), it is $G^{der}$-invariant. Since $G^{der}$ acts transitively on $X^0$ in the real case and
by \rl{inv} in the $p$-adic one, $Y\an\cap(X^0\times\{1\})$ has to be all of $X^0\times\{1\}$. It follows that $Y=X$.

f) holds, since the projection $\pi_S:X\an\to X_S\an$ is a local isomorphism (see the proof of a)).

g) The real case is clear, the $p$-adic case is deep. It uses Yau's theorem 
(see \cite[Rem. 2.2.13]{Ku}).
\end{pf}

\begin{Rem} \label{R:an}
The functorial property of projective limits implies that $X\an$ satisfies
the functorial properties of analytic spaces associated to schemes 
(see \cite[Thm. 3.4.1]{Be1} or \cite[Exp. XII, Thm. 1.1]{SGA1}). 
\end{Rem}

\begin{Lem} \label{L:conn}
Let $\Gm\subset G\times E$ and $X$ be as above, let $E'$ be a compact normal subgroup of $E$,
and let $\Gm '\subset G\times(E'\bs E)$ be the image of $\Gm$ under the natural projection. 
Then we have the following:\par
a) the map $\varphi :\Gm\to\Gm '$ is an isomorphism; \par

b) $\Gm '$ satisfies the conditions of Lemma \ref{L:basic}; \par

c) the quotient $E'\bs X$ exists and is isomorphic to the \esc{E'\bs E}{L} corresponding to $\Gm '$.
\end{Lem}

\begin{pf}
a) The composition map $\Gm\overset{\varphi}{\to}\Gm'\overset{\pr_{G}}{\to} G$
is injective, therefore $\varphi$ is an isomorphism, and $\pr_{G}:\Gm '\to G$
is injective. \par

b) $\Gm '$ is clearly cocompact. Let $U\times S\subset G\times(E'\bs E)$ be an 
open neighbourhood of the identity with a compact closure. 
Then $\varphi^{-1}(U\times S)$ is an open neighbourhood of the 
identity of $G\times E$ with a compact closure. It follows that 
$\varphi^{-1}(U\times S)\cap\Gm$ is finite, thus
$(U\times S)\cap\Gm '$ is  also finite. Hence $\Gm '$ is discrete.\par

c) Since $E'$ is compact and normal, we have $E'S=SE'\in\ff{E}$ for each $\SIE$.
Hence $E'\bs X:=\underset{\underset{S}{\longleftarrow}}{\lim}\,X_{E'S}$ is the required
quotient. Next we notice that for each $\SIE$ the subgroup $\bar{S}:=S\bs E'S$ belongs
to $\ff{E'\bs E}$ and that each $T\in\ff{E'\bs E}$ is of this form. Since
$X_{E'S}\an\cong E'S\bs[X^0\times E]/\Gm\cong\bar{S}\bs[X^0\times (E'\bs E)]/\Gm'$,
we are done.
\end{pf}

\subsection{Drinfel'd's covers} \label{SS:Drin}

\begin{Emp} \label{E:Dr}
 Now we need to recall some Drinfel'd's results \cite{Dr} concerning
covers of $\Om$. (A detailed treatment is given in \cite{BC} for $d=2$ and in \cite{RZ2}
for the general case.)

 Let $K_w$ be as before, and let $D_w$ be a central skew field over $K_w$ with
 invariant $1/d$. Let $\cal O_{D_w}\subset D_w$ be the ring of integers.
Fix a maximal commutative subfield  $K_{w}^{(d)}$ of $D_w$ unramified over $K_w$. 
Let $\pi\in K_w$ be a uniformizer, and let
$Fr_w$ be the Frobenius automorphism of $K_{w}^{(d)}$ over $K_w$. Then $D_w$ is
generated by $K_w$ and an element $\Pi$ with the following defining relations: 
$\Pi^d=\pi$, $\Pi\cdot a=Fr_w(a)\cdot\Pi$ for each $a\in K_{w}^{(d)}$. \par 

 Denote by $\Hat{\cal O}_{w}^{nr}$ the ring of integers of the completion of 
the maximal unramified extension $\Hat{K}_{w}^{nr}$ of $K_w$. Drinfel'd had 
constructed a commutative formal group $Y$ over 
$\Hat{\Omega }^{d}_{K_w}\Hat{\otimes}_{\cal O_{w}}\Hat{\cal O}_{w}^{nr}$ 
with an action of $\cal O_{D_w}$ on it. For a natural number $n$ denote 
by $\Gm_n$ the kernel of the homomorphism $Y\overset{\pi^n}{\to} Y$.
 Let $\C{X}_{n}:=\Gm_n\Hat{\otimes}_{\C{O}_w}K_w$ be the generic fiber of $\Gm_n$,
and let $\C{X}_{n-1/d}\subset\C{X}_{n}$ be the kernel of $\Pi^{nd-1}(=\pi^{n-1/d})$.\par

 Put $\si{d,n}{K_w}:=\C{X}_{n}-\C{X}_{n-1/d}$, and set $T_n:=1+\pi^n\cal O_{D_w}\in
\ff{D_w\m}$. Then 
$\si{d,n}{K_w}$ is an \'etale Galois covering of $\si{d,0}{K_w}:=\Om\Hat
{\otimes}
_{K_w} \Hat{K}_{w}^{nr}$ with Galois group $(\cal O_{D_w}/\pi^n)^{\times}\cong
\cal O_{D_w}^{\times}/T_n$.
We also denote $\cal O_{D_w}^{\times}$ by $T_0$. The action of $\pi$
induces \'etale covering maps $\pi_n:\si{d,n}{K_w}\to\si{d,n-1}{K_w}$,
 giving a  $\Hat{K}_w^{nr}$-pro-analytic space $\Si:=\{\SiN\}_n$.
The group $\cal O_{D_w}^{\times}$ acts naturally on $\Si$ and 
we have $\SiN\cong T_n\bs\Si$ for each $n\in\B{N}\cup\{0\}$. Moreover,
Drinfel'd had also constructed an action of the group $\gp\times D_{w}^
{\times}$ on $\Si$, viewed as a pro-analytic space over $K_w$, which extends the action of
 $\cal O_{D_w}^{\times}$ and satisfies the following properties (notice that our conventions \re{anspace} differ from
those of Drinfel'd):\par

a) the diagonal subgroup $\{(k,k)\in \gp\times D_{w}\m|k\in K_w\m\}$ acts  trivially;

b) $\gp $ (resp. $D_{w}^{\times}$) acts on $\si{d,0}{K_w}=\Om\Hat{\otimes}_{K_w}\Hat{K}_{w}^{nr}$
by the product of the natural action of $\pgp$ on $\Om$ (resp.  the trivial action on  $\Om$)
and the Galois action $g\mapsto Fr_{w}^{\val_w(\Det(g))}$ (resp. $g\mapsto Fr_{w}^{-\val_w(\Det(g))}$) 
on $\Hat{K}_{w}^{nr}$.
\end{Emp}

\begin{Emp} \label{E:deq1}
In the case $d=1$ Drinfel'd's coverings can be described explicitly. Let $L$ be a $p$-adic field.
Then by property a) above, the action of $L\m\times L\m$ on $\si{1}{L}$ is determined uniquely by its restriction to 
the second factor. Denote  by $\theta_{L}:L^{\times}\to \Gal(L^{ab}/L)$ the Artin homomorphism
(sending the uniformizer to the arithmetic Frobenius automorphism). 


\end{Emp}

\begin{Lem} \label{L:Galact}
One has $\si{1}{L}\cong\C{M}(\Hat{L}^{ab})$, and the action of
$(1,l)\in\{1\}\times L\m$ on $\si{1}{L}$ is given by the action of 
$\theta_L(l)^{-1}\in \Gal(L^{ab}/L)$ on $L^{ab}$. 
\end{Lem}

\begin{pf}
This follows from the fact that Drinfel'd's construction  for $d=1$
is equivalent to the construction of Lubin-Tate of the maximal
abelian extension of $L$ (see, for example, \cite[Ch. VI, $\S$3]{CF}).
\end{pf}

\begin{Emp} \label{E:pmaxtor}
Let $L$ be an extension of $K_w$ of degree $d$ and of ramification index $e$. 
For every embeddings $L\hra \Mat_d(K_w)$, $L\hra D_w$ (such exist by 
\cite[Ch. VI, \S 1, App.]{CF})
and $K_{w}^{nr}\hra L^{nr}$ and for every $n\in\B{N}\cup\{0\}$ there exists a closed $L$-rational
embedding $i_n:\si{1,en}{L}\hra \SiN$, which is $(L\m\times L\m)$-equivariant and commutes 
with the projections $\pi_n$.
Moreover, $i_0:\om{1}{L}\Hat{\otimes}_{L}\Hat{L}^{nr}\hra\Om\Hat{\otimes}_{K_w}
\Hat{K}_{w}^{nr}$ is the product of our embedding
$\Hat{K}_{w}^{nr}\hra\Hat{L}^{nr}$ and a closed embedding $i:\om{1}{L}\hra\Om$,
with image $(\Om)^{L^{\times}}$ (see \cite[Prop. 3.1]{Dr}).
Taking an inverse limit we obtain an embedding 
$\wt{\imath}:\si{1}{L}\hra\Si$.
\end{Emp}

\begin{Lem} \label{L:image}
 Let $H$ be a subgroup of $\bold{R}_{L/K_w}(\B{G}_m)(K_w)\cong L\m$, Zariski dense in $\bold{R}_{L/K_w}(\B{G}_m)$. Then $\text{Im}\,\wt{\imath}=
\{x\in\Si |(l,l)x=x \text{ for every }l\in H\}$.
\end{Lem}

\begin{pf}
Since for each $l\in{H}\subset L\m$ the action of $(l,l)$ on $\si{1}{L}$ is trivial, 
and since
$\wt{\imath}$ is $(L\m\times L\m)$-equivariant, $\text{ Im}\,\wt{\imath}$ is contained in the set of fixed points of 
$(l,l),\,l\in{H}$.

Conversely, if $x\in\Si$ is fixed by all $(l,l),\,l\in{H}$, then its
 image $\bar{x}\in\Om$ under the natural projection $\rho:\Si\to\Om$ belongs to
$(\Om)^{H}=(\Om)^{L^{\times}} =i(\om{1}{L})$. Since $\rho(\im\,\wt{\imath})=\im\,i$, 
there exists $y\in \im\,\wt{\imath}$ such that $\rho (y)=\bar{x}(=\rho (x))$. 
Recall that $\Om =\mDw\bs\Si$. Therefore $y=\dt x$ for some $\dt\in\mDw$. 
It follows that $(l,\dt l\dt^{-1})y=y$ for each $l\in H$, hence also
$(1,\dt l\dt^{-1}l^{-1})y=y$. Since the covering
$\Si\to\mODw\bs\Si$ is \'etale, the group $\mODw$ acts freely on $\Si$.
Therefore $\dt l\dt^{-1}l^{-1}=1$ for each $l\in H$. Hence
$\dt$ belongs to  the centralizer of $H$ in $D_w\m$, so that to $L\m$.
It follows that $x=\dt^{-1}y\in L\m\cdot \im\,\wt{\imath}=\im\,\wt{\imath}$.
\end{pf}

\begin{Prop} \label{P:padact}
For each $n\in\B{N}\cup\{0\}$ the group $SD_w\m\cap T_1$ acts trivially on the set $\pi_0$
of connected components of $\si{d,n}{K_w}\Hat{\otimes}_{K_w}\B{C}_p$.
\end{Prop}

\begin{pf}
 Recall (see \re{pmaxtor}) that each maximal commutative subfield $L\subset D_w$ 
gives us (after some choices) a closed $L$-rational $L\m$-equivariant embedding
$i_n:\si{1,en}{L}\to\SiN$.  Let $\C{Y}$ be a  connected component of 
$\si{1,en}{L}\Hat{\otimes}_L\B{C}_p$. Take  $\C{X}\in\pi_0$, which contains $i_n(\C{Y})$.
Then by \rl{Galact}, $\C{X}$ is defined over $\Hat{L}^{ab}$, and
\begin{equation} \label{EE:conn}
l(\C{X})=(\theta_{L}(l))^{-1}(\C{X}) \text{ for each } l\in L\m.
\end{equation}

Fix a $\C{X}_0\in \pi_0$, and let $M$ be the field of definition of $\C{X}_0$. 
Then $M\supset\Hat{K}_w^{nr}$. Since the quotient $D_w\m\bs\SiN\cong\Om$ is geometrically connected,
$D_w\m$ acts transitively on $\pi_0$. Since the action of $D_w\m$ on $\pi_0$ is $K_w$-rational,
$M$ is the field of definition of every $\C{X}\in\pi_0$. In particular, 
$M$ is the closure of a Galois extension of $K_w$, and 
$M\subset\Hat{L}^{ab}$ for every extension $L$ of $K_w$ of degree $d$. Taking $L$ be unramified 
we see that the group $\Aut^{cont}_{K_w}(M)$ of continuous automorphisms of $M$ over $K_w$
is meta-abelian (=extension of two abelian groups).
Set $H:=\{\dt\in D_w\m|$ there exists a $\sigma(\dt)\in \Aut^{cont}_{K_w}(M)$ such that
$\dt(\C{X_0})=\sigma(\dt)^{-1}(\C{X_0})\}$. Then $H$ is a group,
and $\sigma:H\to\Aut^{cont}_{K_w}(M)$ is a well-defined homomorphism.

We claim that $H=D_w\m$. 
Take a $\dt\in D_w\m$, then $K_w[\dt]$ is a 
commutative subfield of $D_w$. Let $L$ be a maximal commutative subfield of $D_w$
containing $\dt$. Then by (\ref{EE:conn}), 
$\dt(\C{X})=(\theta_L (\dt))^{-1}(\C{X})$ for some $\C{X}\in\pi_0$. Take
$\dt'\in D_w\m$ such that $\C{X}=\dt'(\C{X_0})$.
Then  $(\dt')^{-1}\dt\dt'(\C{X_0})=(\dt')^{-1}\circ(\theta_L(\dt))^{-1}\circ\dt'(\C{X_0})
=(\theta_L(\dt))^{-1}(\C{X_0})$, so that $(\dt')^{-1}\dt\dt'\in H$.
Thus each element of $D_w\m$ is conjugate to some element of $H$. In particular, $Z(D_w\m)\subset H$. 
Since $T_n$ acts trivially on 
$\SiN$, it is also contained in $H$. Hence $H\supset T_n\cdot Z(D_w\m)$ has a
finite index in $D_w\m$. Therefore our claim follows from the following 

\begin{Lem} \label{L:gps}
 Let $G$ be a group, and let $H$ be a subgroup of $G$ of finite index. Suppose
that $G=\underset{g\in G/H}{\bigcup}g H g^{-1}$. Then $G=H$.
\end{Lem}

\begin{pf}
 Set $K:=\underset{g\in G}{\bigcap}gHg^{-1}$. Then $K$ is a normal subgroup of $G$
of finite index, and $G/K=\underset{g\in G/H}{\bigcup} g(H/K)g^{-1}=
\underset{g\in G/H}{\bigcup} (g(H/K)g^{-1}-\{1\})\cup \{1\}$. Hence 
$|G/K|\leq |G/H|(|H/K|-1)+1=|G/K|-|G/H|+1$, therefore $G=H$.  
\end{pf}

Now the proposition follows from the fact that
$SD_w\m$ is the derived group of $D_w\m$ (see \cite[1.4.3]{PR}) and that $T_1\cap 
SD_w\m$ is the derived group of $SD_w\m$ (see \cite[1.4.4, Thm. 1.9]{PR}). 
\end{pf}

\subsection{Second construction} \label{SS:seccon}

\begin{Con} \label{C:con2}
 Suppose that a subgroup $\Gm\subset\gp\times E$ satisfies the following conditions:\par
a) $Z(\Gm)=Z(\gp\times E)\cap\Gm$;

b) the subgroup $\overline{Z(\Gm)}\subset Z(\gp\times E)$ is cocompact;

c) $P\Gm\subset\pgp\times PE$ satisfies the assumptions of Lemma 
\ref{L:basic}
(this imply, in particular, that the closure of $\Gm$  is cocompact in $\gp\times E$);

d) the intersection of $Z(\Gm)$ with $Z(\gp)\times\{1\}$ is trivial.

\noindent We are going to associate to $\Gm$ a certain $(D_w\m\times E, K_w)$-scheme.

Consider the quotient $\wt{X}:=(\Si\times E)/\Gm.$ The group $D_w\m\times E$ acts on $\wt{X}$ 
by the product of the natural action of $D_w\m$ on $\Si$ and the left multiplication  by $E$.

\begin{Prop} \label{P:sav2}
For each $S\in\ff{\mDw\times E}$ the quotient $S\bs\wt{X}=S\bs(\Si\times E)/\Gm$ has a natural structure
of a $K_w$-analytic space, which has a unique structure $X_S$ of a projective scheme over $K_w$.
\end{Prop}


\begin{pf}
 First take $S=T_n\times S'$ for some $n\in\B{N}\cup\{0\}$ and some sufficiently small  $S'\in\ff{F}$ 
(to be specified later). Then $S\bs\wt{X}=S'\bs(\SiN\times E)/\Gm$ 
is a disjoint union of $|S\bs E/\GmE|<\infty$ (as in Lemma \ref{L:basic})
quotients of the form
$\Gm_{aS'a^{-1}}\bs\SiN$ with $a\in E$. Thus it remains to prove the statement for quotients
$\Gm_{aS'a^{-1}}\bs\SiN$. For simplicity of notation we assume that $a=1$.
Set $$\Gm_{S',0}:=\Gm_{S'}\cap \pr_{G}(Z(\Gm))= \pr_{G}(\Gm\cap 
[Z(\gp)\times(Z(E)\cap S')]).$$ 
First we construct the quotient $\Gm_{S',0}\bs\SiN$.
Assumptions a) and b) of \ref{C:con2} imply that the closure of $\Gm_{S',0}$ is cocompact in 
$Z(\gp)\cong K_w\m$, hence $\val_w(\Det(\Gm_{S',0}))=dk\B{Z}$ for some $k\in\B{N}$.
Let $K^{(dk)}_{w}$ be the unique unramified extension  of $K_w$ of degree $dk$,
then  $\Gm_{S',0}\bs\Om\Hat{\otimes}_{K_w}\Hat{K}^{nr}_{w}\cong
\Om\otimes_{K_w} K^{(dk)}_{w}$. 

Consider the natural \'etale projection $\pi_n:\SiN\to\si{d,0}{K_w}\to\Om$. Let $\{\C{M}(A_i)\}_{i\in I}$ 
be an affinoid covering of $\Om$. Since the projection $\SiN\to\si{d,0}{K_w}$ is finite, 
each $\pi_n^{-1}(\C{M}(A_i))\subset\SiN$  
is finite over the affinoid space $\C{M}(A_i\Hat{\otimes}_{K_w}\Hat{K}_w^{nr})$. 
Hence it is isomorphic to an affinoid space $\C{M}(B_i)$ for a certain 
$\Hat{K}_w^{nr}$-affinoid algebra $B_i$, finite over $A_i\Hat{\otimes}_{K_w}\Hat{K}_w^{nr}$.
Since $\pi_n$ is $D_w\m$-invariant, we have a natural action of $\Gm_{S',0}$ on $B_i$. 
%
Set $C_i:=B_i^{\Gm_{S',0}}$. Since an affinoid algebra is noetherian, we see that
$C_i$ is finite over the $K_w$-affinoid algebra $A_i$. 
Hence $C_i$ has a canonical structure of a $K_w$-affinoid algebra (see \cite[Prop. 2.1.12]{Be1}). 
Gluing together the $\C{M}(C_i)$'s we obtain a $K_w$-analytic space $\Gm_{S',0}\bs\SiN$, finite and 
\'etale over $\Om$.

Put $\bar{S}:=S\cdot Z(E)/Z(E)\subset PE$. Then $\bar{S}\in\ff{PE}$. To construct 
$\Gm_{S'}\bs\SiN$ we observe that the action of 
$P\Gm_{\bar{S}}=\Gm_{S',0}\bs\Gm_{S'}$ on $\Gm_{S',0}\bs\SiN$ covers its action
 on $\Om$. Suppose that $S'$ is so small that $\bar{S}$ satisfies part d) of \rp{basic}. 
Recall that by \rl{discr} each $x\in\Om$ has an open analytic neighbourhood
$U_x$ such that $\gm(U_x)\cap U_x\neq\emp$ for all $\gm\in P\Gm_{\bar{S}}
\sm \{1\}$, and, as a consequence, $P\Gm_{\bar{S}}\bs\Om$ is obtained by gluing the $U_x$'s. 
Let $\bar{\pi}_n$ be the natural projection from $\Gm_{S',0}\bs\SiN$ to $\Om$. 
For each  $y\in\Gm_{S',0}\bs\SiN$ set $V_y:=\bar{\pi}_n^{-1}(U_{\bar{\pi}_n(x)})$.
 Then the quotient $K_w$-analytic space $P\Gm_{\bar{S}}\bs
(\Gm_{S',0}\bs\SiN)=\Gm_{S'}\bs\SiN$ is obtained by gluing the $V_y$'s. 

Since $\Gm_{S'}\bs\SiN$ is a finite (and \'etale) covering of 
$P\Gm_{\bar{S}}\bs\Om$, which has a structure of a projective scheme over 
$K_w$ by \cite{Mus,Ku},  $\Gm_{S'}\bs\SiN$ also has such a structure by \rco{GAGA2} and the 
remark following it.

Finally consider an arbitrary $S\in\ff{D_w\m\times E}$. It has a normal subgroup 
$\wt{S}$ of the form $\wt{S}= T_n\times S'$ with sufficiently small $S'\in\ff{E}$,
therefore to complete the proof we can use the same considerations as in the end of 
the proof of \rp{sav1}.
\end{pf}

The same argument as in Construction \ref{C:con1} gives us a
\esc{\mDw\times E}{K_w} $X=\underset{\underset{S}{\longleftarrow}}{\lim}\,X_{S}$.
\end{Con}

\begin{Prop} \label{P:prop2}
a) The kernel $E_0$ of the action of $D_w\m\times E$ on $X$ is the closure of the subgroup $Z(\Gm)\subset Z(\gp\times E)=Z(D_w\m\times E)$ after the natural identification $Z(\gp)=K_w\m=Z(\mDw)$.

b) Let $\wt{E}_0$ be the closure of $Z(\Gm)$ in $E$, and let 
$\Gm'\subset\pgp\times(\wt{E}_0\bs E)$ be the image of $\Gm$ under the
natural projection. Then $\Gm'$
satisfies the assumptions of \rl{basic}.

c)  The quotient $D_w\m\bs X $ exists and is isomorphic to the $\wt{E}_0\bs E$-scheme corresponding to $\Gm'$ by \rc{con1}.

d) The quotient $(\mDw\times Z(E))\bs X$ exists and is isomorphic to
the \esc{PE}{K_w} $X'$ corresponding to $P\Gm$ by \rc{con1}.

e) For each $x\in X$ the orbit $(D_w\m\times E)x$ is Zariski dense in $X$.

f) For each sufficiently small $\SIE$  and each $n\in\B{N}\cup\{0\}$ the map 
$X\to X_{T_n\times S}$ is \'etale, and $B^{d-1}$ is the universal covering of each 
connected component of $(X_{T_n\times S,\B{C}})^{an}$ for each embedding 
$K_w\hra \B{C}$. In particular, the projective system $X\an:=\{X_T\an\}_{T\in\ff{D_w\m\times E}}$,
 associated to $X$, is a $K_w$-pro-analytic space.
\end{Prop}

\begin{pf}
a) Notice that $g\in E_0$ if and only if  $g$ acts trivially on $X_S$ 
(or, equivalently, on  $X_{S}^{an}=S\bs(\Si\times E)/\Gm$) and normalizes $S$
for each $S\in\ff{\mDw\times E}$.
For each $\gm\in Z(\Gm)\subset Z(D_w\m\times E)$ let  $\gm_w$ be the projection of $\gm$ to the first factor.
Since $(\gmG\times\gm_w)$ acts trivially on $\Si$, we have $\gm([x,e])=[\gm_w(x),\gmE e]\sim
[(\gmG\times\gm_w)(x),e]=[x,e]$ 
for each $x\in\Si$ and $e\in E$, that is $\gm$ acts trivially on each $X_S\an$. 
Since $\gm$ is central, it certainly normalizes $S$. This shows that the closure 
of $Z(\Gm)$ is contained in $E_0$.

Conversely, suppose that some $(g_1,g_2)\in D_w\m\times E$ with $g_1\in\mDw$ and $g_2\in E$ belongs to $E_0$.
Choose  $S'\in\ff{E}$ and $n\in\B{N}\cup\{0\}$. It suffice to show that 
$(g_1,g_2)\in (T_n\times S')Z(\Gm)$. Since
$(g_1,g_2)$ acts trivially on $S'\bs(\SiN\times E)/\Gm$, we have
$[g_1(x),g_2]\sim [x,1]$ for each $x\in\SiN$. 
This means that there exists an element $\gm=\gm_x\in\Gm$
such that $g_1(x)=\gm_{G}^{-1}(x)$ and $g_2\in S'\gm_{E}$.
Let $x'$ be the projection of $x$ to $\si{d,0}{K_w}$, and let $x''$ be its
projection to $\Om$. The group $\mDw$ acts trivially on $\Om$, therefore
$\gm_{G}^{-1}(x'')=x''$. Choose $x$ so that no non-trivial element of 
$\pgp$ fixes $x''$, then $\gm_{G}\in Z(\gp)\cong K_w\m$. 

Assumption c) of \ref{C:con2} implies that $\gm\in Z(\Gm)$. Since
$g_1(x)=\gmG^{-1}(x)=\gm_w(x)$, we conclude that $g^{-1}_1\gm_w(x)=x$. 
Hence $x'=(g_1^{-1}\gm_w)(x')=Fr_{w}^{-\val(\Det(g_1^{-1}\gm_w))}(x')$,
so that $g_1^{-1}\gm_w\in\mODw$.
Since $\SiN$ is an \'etale Galois covering
of $\si{d,0}{K_w}$ with Galois group $\mODw/T_n$, the equality
$(g_1^{-1}\gm_w)(x)=x$ implies that $g_1^{-1}\gm_w\in T_n$.
It follows that 
$(g_1,g_2)\in (T_n\times S')(\gm_w,\gmE)\subset(T_n\times S')Z(\Gm)$, as claimed.

b) The natural projection $\pgp\times (\wt{E}_0\bs E)\to\pgp\times PE$ induces an 
isomorphism $\Gm'\isom P\Gm$. Hence $\Gm'$ is discrete and has
injective projection to $\pgp$. It is cocompact, because so is
$\Gm\subset\pgp\times E$.

c) Notice first that for each open subgroup $E_0\subset S\subset D_w\m\times E$ ,
compact modulo $E_0$, the quotient $S\bs X$ exists and is projective.
Assumption b) of \ref{C:con2} implies that $E_0$ is cocompact in $D_w\m\times Z(E)$.
Therefore for each $S\in\ff{D_w\m\times E}$ the quotient
$D_w\m S\bs X=(D_w\m E_0 S)\bs X=(D_w\m\times \wt{E}_0)S\bs X$ exists.
Set $\bar{S}:=(D_w\m\times \wt{E}_0)\bs(D_w\m\times \wt{E}_0)S\in\ff{\wt{E}_0\bs E}$.
Then $(D_w\m S\bs X)\an\cong(D_w\m\times \wt{E}_0)S\bs[\Si\times E]/\Gm\cong 
\bar{S}\bs[\Om\times(\wt{E}_0\bs E)]/\Gm'$, and the statement follows as in the proof of \rl{conn} c).

d) follows from c) and \rl{conn} c).

e) follows from c) and \rp{prop1} e).

f) Take $T\in\ff{PE}$ satisfying part d) of Proposition
\ref{P:basic}. Then there exists  $\SIE$ such that $Z(E)\bs S\cdot Z(E)=T$.
Since we have shown in the proof of \rp{sav2} that
$X_{T_n\times S}$ is \'etale over $T\bs X'$ for each $n\in\B{N}\cup\{0\}$, 
the statement follows immediately from \rp{prop1} f), g).
\end{pf}

\begin{Cor} \label{C:etcom}
For each $a\in E$ the composition map 
$$\rho_a:\Si\isom\Si\times\{a\}\hra(\Si\times E^{disc})/\Gm\to X^{an}$$ 
of pro-analytic spaces over $K_w$ is \'etale and one-to-one.
\end{Cor}

\begin{pf}
The \'etaleness is clear. Let $x_1$ and $x_2$ be points of $\Si$ such that $\rho_a(x_1)=\rho_a(x_2)$, and 
let $\bar{a}\in PE'$ be
projection of $a$. Since the composition map 
$$\Om\isom\Om\times\{\bar{a}\}\hra(\Om\times (PE')^{disc})/P\Gm
\isom (X')^{an}$$ is injective, we conclude 
that $x_1$ and $x_2$ have the same projection $y\in\Om$. 

Choose $\SIE$ so small that the group $P\Gm_{aSa^{-1}}$ is 
torsion-free (use \rp{basic} d)). Then no non-central element of $\Gm_{aSa^{-1}}$ fixes $y$.
For each $n\in \B{N}$ let $\pi_{n,S}$ be the projection 
$X^{an}\to (X_{T_n\times S})^{an}$. Then the image of $\pi _{n,S}\circ \rho_a$ is
isomorphic to $\Gm_{aSa^{-1}}\bs\SiN$. Hence there exists $\gm_n\in\Gm_{aSa^{-1}}$
such that the projections of $\gm_n(x_1)$ and $x_2$ to $\SiN$ coincide. Therefore
$\gm_n(y)=y$, so that $\gm_n\in Z(\GL{d}(K_w))=K_w\m$. It follows that the sequence
$\{\gm_n\}_n$ converges to some $\gm\in K_w\m$, which satisfies $\gm(x_1)=x_2$.
Then $(\gm,1)\in Z(D_w\m\times E)$ fixes $z:=\rho_a(x_1)=\rho_a(x_2)$. Since $(\gm,1)$ is central,
it then fixes the whole $(D_w\m\times E)$-orbit of $z$. Hence by \rp{prop2} c), it acts trivially on  
$X$. Therefore by \rp{prop2} a), the element $(\gm,1)$ belongs to 
$\overline{Z(\Gm)}\subset Z(\GL{d}(K_w)\times E)$. Assumption d) of \ref{C:con2} implies that $\gm=1$,
hence $x_1=x_2$.
\end{pf}

\subsection{Relation between the $p$-adic and the real constructions}

The following proposition (and its proof) is a modification of Ihara's theorem
(see \cite[Prop. 1.3]{Ch2}). It will allow us to establish the connection between the 
$p$-adic (\ref{C:con1}, \ref{C:con2}) and the real (or complex) (\ref{C:con1})
constructions.

\begin{Prop} \label{P:relation}
 Let $X$ be an \esc{E}{\B{C}}. Suppose that \par
a) $E$ acts faithfully on $X$;\par 

b) $E$ acts transitively on the set of connected components of $X$; \par

c) there exists $\SIE$ such that the projection $X\to X_S$ is \'etale, 
and $B^{d-1}$ is the universal covering of each connected component of $X_{S}\an$.
\par
Then $X$ can be obtained from the real case of Construction \ref{C:con1}.
\end{Prop}

\begin{Rem} \label{R:remrel}
 a) It follows from Proposition \ref{P:prop1} that all the above
 conditions are necessary.\par

b) Let $X$ be an \esc{E}{\B{C}}, and let $E_0$ be the kernel of the action of $E$ on $X$.
Then $X$ is an \esc{E_0\bs E}{\B{C}} with a faithful action of $E_0\bs E$.
Conversely, every \esc{E_0\bs E}{\B{C}} can be viewed as an \esc{E}{\B{C}}
with a trivial action of $E_0$.\par

c) Let $X$ be an \esc{E}{\B{C}}, and let $X_0$ be a connected component of $X$.
Put $X':=\underset{g\in E}{\bigcup}g(X_0)$. Then $X'$ is an \esc{E}{\B{C}} with 
a transitive action of $E$ on the set of its connected components,
and $X$ is a disjoint union of such \esc{E}{\B{C}}s.\par
Remarks b) and c) show that assumptions a) and b) of the proposition are
not so restrictive.
\end{Rem}

\begin{pf}
We start the proof with the following
\begin{Lem} \label{L:prlim}
Suppose that $\{X_{\al}\}_{\al\in I}$ is a projective system of complex 
manifolds such that the transition maps $X_{\beta}\to X_{\al}$, where 
$\al,\beta\in I$ with $\beta\geq\al$, are analytic coverings. Then there exists a 
projective limit $X$ of the $X_{\al}$'s in the category of complex manifolds.
\end{Lem}

\begin{pf}
Choose an $\al\in I$. Cover $X_{\al}$ by open balls $\{U_{\al\beta}\}_{\beta\in J}$,
and let $\pi:X'\to X_{\al}$ be an analytic covering. Then the
inverse image $\pi^{-1}(U_{\al\beta})$ of each $U_{\al\beta}$ is a disjoint union
of analytic spaces, each of them isomorphic to $U_{\al\beta}$ under $\pi$.
Hence the construction of the projective limit from the proof of 
\rp{prop1}, a) can be applied.
\end{pf}

Now we return to the proof of the proposition. By assumption c), $X_S\an$ is a complex manifold 
for each sufficiently small $\SIE$, and the natural covering $X_S\an\to X_T\an$ is 
\'etale (analytic) for each $T\subset S$ in $\ff{E}$. Therefore by the lemma there exists an 
analytic space $X\an:=\underset{\underset{S}{\longleftarrow}}{\lim}\,X_S\an$.

Since $X_S$ is a complex projective scheme for each $\SIE$,
the set of its connected components coincides with
the set of connected components of $X_S\an$. Hence assumption b) implies that
the group $E$ acts  transitively on the set of  connected components of $X^{an}$.

 Let $M$ be a connected component of $X^{an}$. Denote by $\GmE$
 the stabilizer of $M$ in $E$. Then $\GmE$ acts naturally on $M$, and
the transitivity statement above implies that 
$X^{an}\cong (M\times E^{disc})/\GmE$.

For each $\SIE$ the analytic space $X_{S}^{an}\cong S\bs (M\times E)/\GmE$
is compact. Therefore, as in the proof of Lemma \ref{L:basic}, 
$|S\bs E/\GmE|<\infty$ and $[\GmE:\GmE\cap S]=\infty$. Note that 
$M_S:=(\GmE\cap S)\bs M$ is a connected component of  $X_{S}^{an}$. Suppose that 
$S$ satisfies condition c), then the map $M\to M_S$ is \'etale, and 
$B^{d-1}$ is the universal covering of $M_S$. Hence it is also the universal 
covering of $M$. It follows that $\GmE\subset \Aut(M)$ can be lifted to
$\Gm _{\B{R}}\subset\Aut(B^{d-1})=\pgr^0$.

The kernel $\Dt$ of the natural homomorphism $\pi:\Gm _{\B{R}}\to\GmE$ is the
fundamental group of $M$. Let $\GmS\subset\pgr^0$ be the fundamental group
of the compact analytic space $M_S$, then  $\GmS$ is a cocompact lattice in
$\pgr^0$, satisfying $\GmS =\pi^{-1}(\GmE\cap S)$. It follows that $[\Gm _{\B{R}}:
\GmS]=[\GmE:\GmE\cap S]=\infty$. Therefore, as in the proof of Proposition
\ref{P:basic} a), we see that $\Gm _{\B{R}}$ is dense in $\pgr^0$. The group $\Dt$ 
is discrete in $\pgr^0$ and normal in  $\Gm _{\B{R}}$, thus it is trivial
 (compare the proof of Proposition \ref{P:basic} b)). In particular,
$M\cong B^{d-1}$ and $\pi$ is an isomorphism.

Put $\Gm:=\{(\gm ,\pi(\gm))|\gm\in\Gm _{\B{R}}\}\subset\pgr^0\times E$.  
Since $\GmS$ is discrete in $\pgr^0$, so is $\Gm$ in $\pgr^0\times E$. Let 
$K\subset\pgr^0$ be the stabilizer of $0\in B^{d-1}$. Then $X_{S}^{an}\cong S\bs 
(B^{d-1}\times E)/\Gm=(K\times S)\bs(\pgr^0\times E)/\Gm$.
Since $K,S$ and $X_S\an$ are compact, $\Gm$ is cocompact in 
$\pgr^0\times E$. Since $\Ker(\pr_G)$ equals the kernel of the action of $E$ on $X$,
the projection $\pr_G$ is injective. This shows that $\Gm$ satisfies all the assumptions 
of Construction \ref{C:con1}.
\end{pf}

\begin{Cor} \label{C:rel}
Choose an embedding $K_w\hra\B{C}$. Let $X$ be an \esc{E}{K_w} obtained by the
$p$-adic case of Construction \ref{C:con1} or an \esc{\wt{E}}{K_w} obtained by
Construction \ref{C:con2}. Then $X_{\B{C}}$ can be constructed by the real case
of Construction \ref{C:con1}.
\end{Cor}

\begin{pf}
 This is an immediate consequence of Propositions \ref{P:relation}, \ref{P:prop1}
and \ref{P:prop2}.
\end{pf}

\subsection{Elliptic elements} \label{SS:ellele}

\begin{Def} \label{D:ell}
Suppose that a group $G$ acts on a (pro-)analytic space (or a scheme) $X$. An element $g\in G$
is called {\em elliptic} if it has a fixed point $x$ such that the linear
transformation of the tangent space of $x$, induced by $g$, has no non-zero 
fixed vectors.
In such a situation we call $x$ an {\em elliptic point} of $g$.
\end{Def}

\begin{Lem} \label{L:lintr}
 Let $\la_1,\la_2,...,\la_d$ be the eigenvalues of some element $g\in \GL{d}(L)$
(with multiplicities).
 Let $v\in\B{P}^{d-1}(\overline{L})$ be one of the fixed points of $g$ corresponding to $\la_1$.
 Then $\frac{\la_2}{\la_1},\frac{\la_3}{\la_1},...,
\frac{\la_d}{\la_1}$ are the eigenvalues of the linear transformation of
the tangent  space of $v$, induced by $g$.
\end{Lem}

\begin{pf}
Simple verification.
\end{pf}

\begin{Prop} \label{P:ell}
 The set of elliptic elements of $\pgr^0$ with respect to its action on $B^{d-1}$ and of $\psp$ with respect to its action on $\Om$ is open and non-empty.
\end{Prop}

\begin{pf}
In the real case we observe that an element $g:=diag(\la_1,\la_2,...,\la_d)\in\pgr^0$ fixes
$(0,0,...,0)\in B^{d-1}$. Therefore by Lemma \ref{L:lintr}, $g$ is elliptic
if $\la_i\neq\la_d$ for all $i\neq d$. It follows that the set of elliptic elements
is non-empty. It is open, because if $g$ has a fixed point in $B^{d-1}$
corresponding to an eigenvalue of $g$ appearing with multiplicity $1$, then the same is
true in some open neighbourhood of $g$.

In the $p$-adic case we start with the following 
\begin{Lem} \label{L:padell}
 An element $g\in\gp$ is elliptic (acting on $\Om$) if and only if its characteristic
 polynomial is irreducible over $K_w$.
\end{Lem} 
\begin{pf}
 Suppose that the characteristic polynomial $\chi_g$ of $g$ is irreducible
over $K_w$. Then $g$ has $d$ distinct eigenvalues. Let $\la$ be some eigenvalue 
of $g$, let $v\neq 0$ be the eigenvector of $g$ corresponding to $\la$,
and let $\bar{v}\in \B{P}^{d-1}(\overline{K}_w)$ be the fixed point of $g$
corresponding to $v$. By \rl{lintr}, the linear 
transformation of the tangent space of $\bar{v}$, induced by $g$, has no fixed 
non-zero
vector. So it remains to be shown that $\bar{v}\in\Om$. If $\bar{v}\notin\Om$,
then it lies in a $K_w$-rational hyperplane. Therefore there exist elements
$a_1,...,a_d\in K_w$, not all $0$ (say $a_d\neq 0$) such that $(a_1,...,a_d)\cdot
v=0$. 
We also know that $(g-\la I)v=0$. Let $A$ be the matrix obtained from $g-\la I$ by
replacing the last row by $(a_1,...,a_d)$. Then $Av=0$, so that $\Det\,A=0$. 
The determinant of $A$ is a polynomial in $\la$ of degree $(d-1)$ with 
coefficients in $K_w$ with leading coefficient $(-1)^{(d-1)}a_d\neq 0$. This
contradicts the fact that the minimal polynomial of $\la$ over $K_w$ 
has degree $d$.

Suppose now that the characteristic polynomial $\chi_g$ of $g$ equals the product
$f_1\cdot...\cdot f_k$ of polynomials irreducible over $K_w$ ($k>1$).
Consider the matrix $f_1(g)$. If  $f_1(g)=0$, then the minimal
polynomial $m_g$ of $g$ divides $f_1$. Hence each root of $\chi_g$, being
a root of $m_g$, is a root of $f_1$. Each $f_i$ has only simple roots,
therefore $f_i|f_1$ for each $i$. Since $f_1$ is irreducible, 
all the $f_i$'s are equal up to a constant. Hence $\chi_g=c\cdot f_1^k$
for some $c\in K_w\m$. In particular, each root of $\chi_g$ is at least double.
\rl{lintr} then implies that $g$ is not elliptic.

Hence we can suppose that $f_i(g)\neq 0$ for all $i=1,2,...,k$. Let 
$\la$ be an 
eigenvalue of $g$, let $v$ be the eigenvector corresponding to $\la$, and let
$\bar{v}\in\B{P}^{d-1}(\overline{K}_w)$ be the fixed point of $g$ corresponding 
to $v$. Choose  $i\in\{1,...,k\}$ such that $\la$ is a root of $f_i$. Then
$f_i(g)v=f_i(\la)v=0$. The matrix $f_i(g)\neq 0$ has all its entries in $K_w$, 
hence $\bar{v}$ lies in a $K_w$-rational hyperplane. Therefore $g$ is not
elliptic.
\end{pf}

Now we return to the proof of the proposition. Embed an extension $L=K_w(\la)$  of 
$K_w$ of degree $d$ into $\Mat_d(K_w)$. Then $\lambda\in L\m\subset
\gp$ has an irreducible characteristic polynomial over $K_w$. Therefore the set of
elliptic elements of $\pgp$ is non-empty. It is open because by Krasner's lemma
if $g\in\gp$ has 
a characteristic polynomial irreducible over $K_w$, then each $g'\in\gp$, close
 enough to $g$, has the same property (see \cite[Ch. II, $\S$3, Prop. 4]{La}).
It follows that the set of elliptic elements of $\psp$ is also open. For showing 
that it is non-empty observe 
that if an element $g\in\pgp$ is elliptic but $g^d$ is not elliptic, then by
\rl{lintr} the characteristic polynomial of each representative of $g^d$ in
$\gp$ has at least two equal roots. Hence such a $g$ belongs to some proper Zariski 
closed subset of $\PGL{d}$. It follows that there exists an elliptic element
$g\in\pgp$ such that $g^d$ is elliptic as well. Since $g^d$ always belongs to $\psp$, we are done.
\end{pf}

\begin{Prop} \label{P:padell}
a) An element $(g,\dt)\in \gp\times D_w\m$ is elliptic with respect to its 
action on $\Si$ (viewed as a pro-analytic space over $K_w$) if and only if the 
characteristic polynomials of $g$ and $\dt$ are $K_w$-irreducible and 
coincide.\par

b) For every element $g\in\gp$, elliptic with respect to its action on $\Om$, there 
exists a $\dt\in D_w\m$ such that $(g,\dt)$ is elliptic with respect to
its action on $\Si$.
\end{Prop}

\begin{pf}
a) Let $x\in\Si$ be an elliptic point of $(g,\dt)$, and let $\bar{x}\in\Om$ be its
 image under the natural projection $\rho:\Si\to\Om$. Since $\rho$ is \'{e}tale,
it induces an isomorphism of tangent spaces (up to an extension of scalars).
Hence $g$ is elliptic with respect to its action on $\Om$. 
By \rl{padell}, $g$ generates maximal commutative subfield $L:=K_w(g)$ of $\Mat_d(K_w)$.

Choose an embedding $j:K_w(g)\hra D_w$ (such exists by \cite[Ch. VI,
 $\S$1, App.]{CF}). It defines an
$L\m$-equivariant embedding $\wt{\imath}:\si{1}{L}\hra\Si$ (see \re{Dr}).
We know that $\bar{x}\in(\Om)^{L\m}=\rho\circ \wt{\imath}(\si{1}{L})$. 
In particular, there exists $y\in\wt{\imath}(\si{1}{L})$
such that $\rho(y)=\bar{x}$.  Since $\wt{\imath}$ is
$L\m$-equivariant, the element $(g,j(g))\in\gp\times D_w\m$ fixes $y$.
Using the fact that $x\in\rho^{-1}(\bar{x})$ and that $ D_w\m\bs\Si=\Om$, we have 
 $y=d_{0}x$ for some $d_0\in D_w\m$. Hence elements
$(g,d_{0}^{-1}j(g)d_0)\in\gp\times D_w\m$ and $\wt{d}:=d_{0}^{-1}j(g)d_{0}\dt^{-1}\in D_w\m$
fix $x$. In particular, $\wt{d}\in D_w\m$ fixes some point (the projection of $x$) on $\si{d,0}{K_w}
\Hat{\otimes}_{K_w}\B{C}_p$, therefore $\wt{d}\in\C{O}_{ D_w}\m$.
Since the Galois covering $\Si\to\C{O}_{ D_w}\m\bs\Si$ is \'etale, $\C{O}_{ D_w}\m$ acts freely on $\Si$. 
It follows that $\wt{d}=1$, hence $\dt=d_{0}^{-1}j(g)d_0$.
This completes the proof of the implication ``only if'', because $g\in\Mat_d(K_w)$ and $j(g)\in D_w$
have the same characteristic polynomials.

Conversely, suppose that the characteristic polynomials of $g$ and $\dt$
are $K_w$-irreducible and coincide. Then the subfields $K_w(g)\subset\Mat_d(K_w)$
and $K_w(\dt)\subset D_w$ have degree $d$ over $K_w$ and are isomorphic under
the $K_w$-isomorphism sending $g$ to $\dt$. Using this isomorphism we obtain 
embeddings of the field $L:=K_w(g)$ into $\Mat_d(K_w)$ and into $D_w$.
These embeddings define (by \re{pmaxtor}) an $(L\m\times L\m)$-equivariant embedding
$\wt{\imath}:\si{1}{L}\hra\Si$ such that every point $x\in\wt{\imath}(\si{1}{L})$ 
is fixed by all elements of the form 
$(l,l)\in L\m\times L\m\subset\gp\times D_w\m$.
In particular, $x$ is a fixed point of $(g,\dt)$. As before, the action of
$(g,\dt)$ on the tangent space of $x$ coincides with the action  of $g$ on the 
tangent space of $\bar{x}$. Since $\bar{x}$ is an elliptic point of $g$
(by \rl{padell}), $x$   is an elliptic point of $(g,\dt)$.

b) If an element $g\in\gp$ is elliptic, then by \rl{padell} it has an irreducible 
characteristic polynomial over $K_w$. Therefore $K_w(g)\subset\Mat_d(K_w)$ is 
a field extension of $K_w$ of degree $d$. Then for every embedding $j$ of $K_w(g)$
into $D_w$ the element $(g,j(g))$ is elliptic by a). 
\end{pf}

\subsection{Euler-Poincar\'e measures and inner twists}  \label{SS:inner}

Here we give a brief exposition of Kottwitz' result \cite[\S 1]{Ko}.

\begin{Emp} \label{E:EP}
Let $L$ be a local field of characteristic $0$, and let $\bold{H}$ be a connected reductive group
over $L$. Serre \cite{Se1} proved that there exists a unique invariant measure 
(called the Euler-Poincar\'e measure)
$\mu_{\bold{H}}$ on $\bold{H}(L)$ such that $\mu_{\bold{H}}(\Gm\bs \bold{H}(L))$ is equal to
the Euler-Poincar\'e characteristic $\chi_E (\Gm)$ of $H^{*}(\Gm,\B{Q})$ for every torsion-free 
cocompact lattice $\Gm$ in $\bold{H}(L)$.
In particular, $\mu_{\bold{H}}(\bold{H}(L))=1$ if the group $\bold{H}(L)$ is
compact. The Euler-Poincar\'e measure is either always negative, always positive or 
identically zero. It is non-zero if and only if $\bold{H}$ has an anisotropic 
maximal $L$-torus. (A result of Kneser shows that in the $p$-adic case this 
happens if and only if the connected center of $\bold{H}$ is anisotropic.)
\end{Emp}

\begin{Emp} \label{E:Tam}
Let $\bold{G}$ be an inner form of $\bold{H}$.
Choose an inner twisting $\rho:\bold{H}\to\bold{G}$ over $\overline{L}$. Choose a 
non-zero invariant differential form $\omega_{\bold{G}}$ of top degree on $\bold{G}$. 
Set $\omega_{\bold{H}}:=\rho^{*}(\omega_{\bold{G}})$. Using the fact that $\bold{H}$ is reductive, 
that the twisting is inner and that $\omega_{\bold{G}}$ is invariant, we see that $\omega_{\bold{H}}$  is 
invariant, defined over $L$ and does not depend on $\rho$. Hence $\omega_{\bold{G}}$  and $\omega_{\bold{H}}$ define 
invariant measures $|\omega_{\bold{G}}|$ and $|\omega_{\bold{H}}|$ on $\bold{G}(L)$ and $\bold{H}(L)$ 
respectively (see \cite[2.2]{We2}). 
\end{Emp}

\begin{Def} \label{D:compat}
 The invariant measures $\mu$ on $\bold{H}(L)$ and $\mu'$ on $\bold{G}(L)$ are called
{\em compatible} if there exists some $c\in\B{R}$ such that 
$\mu=c|\omega_{\bold{H}}|$ and $\mu'=c|\omega_{\bold{G}}|$. 
\end{Def}

\begin{Emp} \label{E:ass}
Now suppose that $\bold{H}$ has an anisotropic maximal $L$-torus $\bold{T}$, so that the  
Euler-Poincar\'e measure $\mu_{\bold{H}}$ on $\bold{H}(L)$ is non-trivial. (Notice that
for semisimple groups of type $A_n$ this assumption is satisfied automatically). 
Denote by $|\mu_{\bold{H}}|$ the absolute value of $\mu_{\bold{H}}$. Write 
$\C{D}(\bold{T},\bold{H})$ for the finite set $\Ker[\text{H}^1(L,\bold{T})\to\text{H}^1(L,\bold{H})]$,
and write $|\C{D}(\bold{T},\bold{H})|$ for its cardinality. It is well known that $\bold{T}$ transfers to
$\bold{G}$, thus we can also consider the finite set $\C{D}(\bold{T},\bold{G})$.
\end{Emp} 

\begin{Prop} \label{P:compat}
(\cite[Thm. 1]{Ko}) The invariant measure $|\C{D}(\bold{T},\bold{H})|^{-1}|\mu_{\bold{H}}|$ on 
$\bold{H}(L)$ is compatible with the invariant measure $|\C{D}(\bold{T},\bold{G})|^{-1}|\mu_{\bold{G}}|$ on 
$\bold{G}(L)$.
\end{Prop}

\begin{Rem} \label{R:compat}
a) In the $p$-adic case, the sets $\C{D}(\bold{T},\bold{H})$ and $\C{D}(\bold{T},\bold{G})$ always have
the same cardinality.

b) In the real case, $\C{D}(\bold{T},\bold{H})=\Omega(\bold{H}(\B{C}),\bold{T}(\B{C}))/
\Omega(\bold{H}(\B{R}),\bold{T}(\B{R}))$, where $\Omega$ stands for the Weyl group.
In particular, $|\C{D}(diag,\PGU{d})|=1$ and $|\C{D}(diag,\PGU{d-1,1})|$ is $d$ (resp. $1$) if $d>2$ 
(resp. $d=2$).
\end{Rem}

\subsection{Preliminaries on torsors (=principal bundles)} \label{SS:proptor}

\begin{Def} \label{D:tor}
Let $\bold{G}$ be an affine group scheme over a field $L$ (resp. an $L$-analytic group), 
and let $X$ be an $L$-scheme (resp. an $L$-analytic space). 
A {\em $\bold{G}$-torsor } over $X$ is a scheme (resp. an analytic space) 
$T$ over $X$ with
an action $\bold{G}\times T\to T$ of  $\bold{G}$ on $T$ over $X$ such that for some surjective
\'etale covering $X'\to X$ the fiber product $T\times_X X'$ is the trivial $\bold{G}$-torsor over $X'$
(that is isomorphic to $\bold{G}\times X'$).
\end{Def}
 
\begin{Rem} \label{R:compltor}
Since each \'etale morphism of complex analytic spaces is a local isomorphism, 
our definition is this case coincides with the classical one.
\end{Rem}

\begin{Lem} \label{L:tor}
a) If $T$ is a $\bold{G}$-torsor over $X$, then the map 
$\varphi_T:\bold{G}\times T\to T\times_X T\;(\varphi_T(g,t)=(gt,t))$ is an isomorphism. 

b) Let $T$ and  $T'$ be two $\bold{G}$-torsors over $X$ and $Y$ respectively.
Then for each $\bold{G}$-equivariant map $f:T\to T'$ the natural morphism
$T\to T'\times_Y X$ is an isomorphism.
\end{Lem}

\begin{pf}
a) Since the problem is local for the \'etale topology on $X$ (see \cite[Ch. I, Rem. 2.24]{Mi2} 
in the algebraic case, \cite[Prop. 4.1.3]{Be3} in the $p$-adic analytic and \rr{compltor} in the complex one),
we may suppose that $T$ is trivial. Then our morphism $(g,(h,x))\mapsto ((gh,x),(h,x))$ is invertible.

b) For the trivial torsors the statement is clear. The general case follows as in a).
\end{pf}

\begin{Rem} \label{R:altor}
By \cite[Ch. I, Rem. 2.24 and Prop. 3.26]{Mi2} our definition in the algebraic case is equivalent 
to the standard one. In particular, a $\bold{G}$-torsor over $X$ is affine and faithfully flat over $X$.
\end{Rem}

\begin{Lem} \label{L:descent}
Let $X$ be a separated scheme over a field $L$, let $\bold{G}$ and $\bold{H}$ be 
two affine group schemes over $L$, let $T$ be a $\bold{G}$-torsor over $X$, and 
let $\pi:T\to X$ be the natural projection.\par

a) The functor $\C{F}\mapsto \pi^* \C{F}$ defines an equivalence between
the category of quasi-coherent sheaves on $X$ and the category of 
$\bold{G}$-equivariant quasi-coherent sheaves on $T$, that is, quasi-coherent sheaves
on  $T$ with a $\bold{G}$-action that lifts the action of $\bold{G}$ on $T$.\par

b) The functor $Z\mapsto Z\times_X T$ defines an equivalence between the following categories:

\quad i) the category of vector bundles of finite rank on $X$  and the category of $\bold{G}$-equivariant
vector bundles of finite rank on $T$;\par

\quad ii) the category of $\bold{H}$-torsors over $X$ 
and the category of  $\bold{G}$-equivariant $\bold{H}$-torsors over $T$;

\quad iii) (if $X$ is noetherian and regular) the category of $\B{P}^n$-bundles  
on $X$ and the category of $\bold{G}$-equivariant $\B{P}^n$-bundles on $T$.

The quasi-inverse functor is $\wt{Z}\mapsto\bold{G}\bs\wt{Z}$. 
\end{Lem}

\begin{pf}
This is a consequence of a descent theory.\par

a) Abusing notation we will write $\C{F}\times_{Y_2}Y_1$ instead of $\rho^* 
\C{F}$ for every morphism $\rho:Y_1\to Y_2$ and every sheaf of modules $\C{F}$ on
$Y_2$. Let $\wt{\C{F}}$ be a $\bold{G}$-equivariant quasi-coherent sheaf on $T$. 
Define an isomorphism $\phi:(\wt{\C{F}}\times_{T}T)\times_X T\isom T\times_X 
(\wt{\C{F}}\times_T T)$ over $T\times_X T$ by the formula 
$\phi(f,gt,t)=(gt,g^{-1}f,t)$ for all $g\in \bold{G},\; t\in T$ and $f\in\wt{\C{F}}_t$
 (use \rl{tor}). Then $\phi$ satisfies the descent
conditions of \cite[Prop. 2.22]{Mi2}. Since $T\to X$ is affine and faithfully 
flat, there is a unique quasi-coherent sheaf $\C{F}$ on $X$ 
such that $\wt{\C{F}}\cong\C{F}\times_X T$. 
Since the construction of descent is functorial
(see \cite[2.19]{Mi2}), we obtain an equivalence of categories.
Notice that $\C{F}\cong \bold{G}\bs (\C{F}\times_X T)$.\par

b) follows from a) in a standard way (use \cite[II, Ex. 5.18, 5.17 and 7.10]{Ha}).
\end{pf}

From now on we suppose that the reader is familiar with basic definitions of 
tensor categories (see \cite{DM}).

\begin{Not} \label{N:cat}
For a field $L$, an affine group scheme (resp. an analytic group) $\bold{G}$ over $L$
and a scheme (resp. an analytic space) $X$ over $L$:\par
a) let $\underline{\Rep}_L(\bold{G})$ be the category of finite-dimensional 
representations of $\bold{G}$ over $L$ ;\par

b) let $\underline{\vect}_X$ be the category of vector bundles of finite rank
on $X$;\par

c) let $\underline{\Tor}_X(\bold{G})$ be the category of $\bold{G}$-torsors over $X$.

\noindent We will sometimes identify categories with the sets of their objects.
\end{Not}

\begin{Def} \label{D:fibre}
Let $L$ be a field, and let $\bold{G}$ be an affine group scheme over $L$. A 
{\em $\bold{G}$-fibre functor with values in a separated scheme (resp. analytic 
space) $X$} over $L$ is an 
exact faithful tensor functor from $\underline{\Rep}_L(\bold{G})$ to 
$\underline{\vect}_X$.
\end{Def}

\begin{Rem} \label{R:fib}
If $X=\Spec\,R$ is affine, then $\underline{\vect}_X$ is equivalent to 
the category of finitely generated projective modules over $R$, hence our 
definition is a global version of that of \cite[3.1]{DM}.
\end{Rem}

\begin{Emp} \label{E:tor1}
Let $T$ be a $\bold{G}$-torsor over $X$, then by \rl{descent}, the correspondence 
$V\mapsto G\bs (V\times T)$ defines a $\bold{G}$-fibre functor with values in $X$.
This correspondence defines a functor $\nu$ from $\underline{\Tor}_X(\bold{G})$ to
the category of $\bold{G}$-fibre functors with values in $X$.
\end{Emp}

\begin{Thm} \label{T:torsor}
 The functor $\nu$ determines an equivalence between $\underline{\Tor}_X(\bold{G})$ and
the category of $\bold{G}$-fibre functors with values in $X$.
\end{Thm}

\begin{pf}
The local version is \cite[Thm. 2.11 and 3.2]{DM}. The gluing works because $X$ 
is separated.
\end{pf}

\begin{Emp} \label{E:tor2}
Later on we will use the following description of the quasi-inverse 
functor $\tau$ of $\nu$. Let $\eta$ be a $\bold{G}$-fibre functor with values in $X$.
For each morphism $\pi_0:T_0\to X$ we define two tensor functors 
$\eta_1:V\mapsto V\times T_0$ and $\eta_2=\pi_0^{*}\circ\eta$ from 
$\underline{\Rep}_L(\bold{G})$ to $\underline{\vect}_{T_0}$. Let $\mu(T_0,\pi_0)
:=\text{Isom}(\eta_2,\eta_1)$ be the set of isomorphisms of tensor functors.
The action of $\bold{G}$ on the first factor of $V\times T_0$ defines an action of $\bold{G}$
on $\eta_1$, and a fortiori defines an action of $\bold{G}$ on $\mu(T_0,\pi_0)$. Thus $\mu$ is a functor from the category of schemes over $X$ to the category of
sets with a $\bold{G}$-action. \rt{torsor} says that this functor is 
representable by a $\bold{G}$-torsor $\tau(\eta)$ over $X$ 
(see \cite[Thm. 2.11 and 3.2]{DM} and their proofs).
\end{Emp}

\begin{Emp} \label{E:tor3}
Let $T$ be a $\bold{G}$-torsor over $X$. For each 
$V\in\underline{\Rep}_L(\bold{G})$ the identity map of $T$, viewed as a
$T$-valued point of $T$, corresponds to a certain isomorphism 
$\varphi_V : V\times T\isom (\bold{G}\bs(V\times T))\times_X T$. Then
$\varphi_V$ is the quotient of the $\bold{G}$-equivariant isomorphism
$Id_V\times\varphi_T:V\times\bold{G}\times T\isom V\times T\times_X T$
(for the diagonal action of  $\bold{G}$ on the first two factors on both sides)
by the action of $\bold{G}$. Explicitely, $\varphi_V(v,t)=([v,t],t)$.
\end{Emp}

\begin{Prop} \label{P:torsor}
Let $L$ be equal to $K_w$ or to $\B{C}$ as in \re{anspace}. Let $X$ be a projective 
$L$-scheme, and let $\bold{G}$ be a linear algebraic group over $L$. The functor $T\mapsto T\an$ induces an 
equivalence between the category of $\bold{G}$-torsors over $X$ and the category of $\bold{G\an}$-torsors 
over $X\an$.
\end{Prop}

\begin{pf}
A quasi-inverse functor can be described as follows. Let $\wt{\pi}:\wt{T}\to X\an$
be a $\bold{G}\an$-torsor. Then the map $V\mapsto \bold{G}\an\bs(V\an\times \wt{T})$ defines a 
$\bold{G}$-fibre functor with values in $X\an$.
Since the correspondence described in \rcr{GAGA2}
commutes with tensor products, the tensor categories $\underline{\vect}_X$ and
$\underline{\vect}_{X\an}$ are equivalent. Therefore \rt{torsor} gives us an 
algebraic $\bold{G}$-torsor ${\pi}:T\to X$.

It remains to show that there exists a canonical isomorphism $\wt{T}\isom T\an$. 
By the definition of $T$ we have for each $V\in \underline{\Rep}_L(\bold{G})$ 
a canonical isomorphism $\psi_V:\bold{G}\an\bs(V\an\times\wt{T})\isom\bold{G}\an\bs(V\an\times T\an)$.
We also have (as in \re{tor3}) natural isomorphisms $\wt{T}\times V\an\isom \wt{T}\times_{X} 
(\bold{G}\an\bs(\wt{T}\times V\an))$ mapping $(t,v)$ to $(t,[t,v])$. Hence each
point $t_0$ of $\wt{T}$ defines canonical isomorphisms $V\an\cong \{t_0\}\times V\an
\isom \{t_0\}\times_{X}(\bold{G}\an\bs(\wt{T}\times V\an)):\:v\mapsto (t_0,[t_0,v])$.
Since $t_0$ defines a point of $X\an$ and therefore of $X$, it gives us 
by the universal property of $T$ (see \re{tor2}) a point 
$\psi(t_0)\in{T}\an$, satisfying $\psi_V([t_0,v])=
[\psi(t_0),v]$ for all $V\in \underline{\Rep}_L(\bold{G})$. 

Taking $V$ be a faithful representation of $\bold{G}$, we obtain that the map 
(of sets) $\psi:\wt{T}\to T\an$ is $\bold{G}\an$-equivariant, therefore it is
one-to-one and surjective.
It remains to show that the maps $\psi$ and $\psi^{-1}$ are analytic.
Let us prove it, for example, for $\psi$.
Let $\rho:X'\to X\an$ be an \'etale surjective covering such that
$\rho^*(T\an)\cong\bold{G\an}\times X'$. By \cite[Prop. 4.1.3]{Be3} in the $p$-adic case and by \rr{compltor}
in the complex one, it will suffice to show that 
$\rho^*\psi:\rho^*(\wt{T})\to\rho^*(T\an)\cong\bold{G\an}\times X'$ (or just its projection 
to the first factor $\pi':\rho^*(\wt{T})\to\bold{G\an}$) is analytic.
Consider the map 
\begin{equation} 
\wt{\psi}_V:V\an\times\rho^*(\wt{T})\overset{\text{proj}}{\lra}\bold{G}\an\bs[V\an\times\rho^*(\wt{T})]
\overset{\rho^*(\psi_V)}{\lra}\cr \notag
\bold{G}\an\bs[V\an\times\rho^*(T\an)]\cong\bold{G}\an\bs(V\an\times\bold{G\an}\times X')\cong 
V\an\times X'\overset{\text{proj}}{\lra}V\an.\notag
\end{equation}
It is analytic, and satisfies $\wt{\psi}_V(v,t)=(\pi'(t))^{-1}v$. Hence $\pi'$ is analytic as well.
\end{pf}

\begin{Cor} \label{C:torsor}
Let $X$ and $Y$ be projective $L$-schemes, let $\bold{G}$ and $\bold{H}$ be algebraic groups over 
$L$, and let $\psi:\bold{G}\to \bold{H}$ be an algebraic group homomorphism over $L$.
If $T\in\underline{\Tor}_X(\bold{G})$ and $S\in\underline{\Tor}_Y(\bold{H})$, then for every 
$\psi$-equivariant analytic map $\wt{f}:T\an\to S\an$ (that is, satisfying
$\wt{f}(gt)=\psi(g)\wt{f}(t)$ for all $g\in \bold{G}\an$ and $t\in T\an$), there is a 
unique algebraic morphism $f:T\to S$ such that $f\an\cong\wt{f}$.
\end{Cor}

\begin{pf} 
(compare the proof of \rcr{GAGA3})
Since $\wt{f}$ is $\psi$-equivariant, it covers some algebraic morphism $\bar{f}:X\to Y$
(use \rcr{GAGA1}). Therefore  $\wt{f}$ factors through $S\an\times_{Y\an}X\an\cong
(S\times_X Y)\an$. Hence we may suppose, replacing
$S$ by $S\times_X Y$, that $X=Y$ and that $\bar{f}$ is the identity.

Consider the $\bold{H}$-torsor $\bold{H}\times T$ over $T$ equipped with the following $\bold{G}$-action:
$g(h,t)=(h\psi(g)^{-1},gt)$ for all $g\in \bold{G},\;h\in \bold{H}$ and $t\in T$. By \rl{descent}, 
there exists an $\bold{H}$-torsor $\bold{H}\times_{\bold{G}} T:=\bold{G}\bs(\bold{H}\times T)$
over $X$. Let $i$ be the composition of the embedding $t\mapsto(1,t)$ of $T$ into $\bold{H}\times T$ with
the natural projection to $\bold{H}\times_{\bold{G}} T$. Then by the definition,
every $\psi$-equivariant algebraic morphism $\mu:T\to S$ factors as a composition of $i$ 
with the unique $\bold{H}$-equivariant map $\bold{H}\times_{\bold{G}} T\to T$ 
(defined by $[h,t]\mapsto h\mu(t)$). Therefore 
$(\bold{H}\times_{\bold{G}} T)\an\cong \bold{H}\an\times_{\bold{G}\an}
T\an$ is an $\bold{H}\an$-torsor over $X\an$ having the same functorial property.

Now we are ready to prove our corollary. From the $\psi$-equivariance of 
$\wt{f}$ we conclude that it factors uniquely as $\wt{f}:T\an\overset{i\an}
{\lra}(\bold{H}\times_{\bold{G}} T)\an\overset{\wt{f}'}{\lra}S\an$. By the 
proposition, $\wt{f}'$ has a unique underlying algebraic morphism 
$f':\bold{H}\times_{\bold{G}} T\to S$. Set $f:=f'\circ i$. The uniqueness can be 
derived from the above considerations as in the proof of \rcr{GAGA3}.
\end{pf}

Now we recall the notion and basic properties of connections on torsors (following
\cite[Ch. VI, $\S$1]{St}).

\begin{Def} \label{D:con}
Let $X$ be a smooth scheme or an analytic space, and let $\pi:P\to X$ be a $\bold{G}$-torsor. 
A {\em connection on $P$} is a $\bold{G}$-equivariant vector
subbundle $\C{H}$ of the tangent bundle $T(P)$ of $P$ such that ${\pi_*}_{|\C{H}_p}$ is an isomorphism
$\C{H}_p\isom T_{\pi(p)}(X)$ for each $p\in P$.
\end{Def}

\begin{Emp} \label{E:con}
Starting from the isomorphism $\varphi_P:\bold{G}\times P\isom P\times_{X}P$ we obtain an isomorphism 
of tangent spaces $(\varphi_P)_*:T_e(\bold{G})\times T_p(P)\isom T_p(P)\times_{T_{\pi(p)}(X)}T_p(P)$ 
and an identification $(X\mapsto\text{proj}_1((\varphi_P)_*(X,0))$ 
of $\C{G}:=\Lie(\bold{G})=T_e(\bold{G})$ with the tangent space to the
 fiber through $p\in P$. Therefore a connection $\C{H}$ on $P$ gives us a canonical
decomposition $T_p(P)=\C{G}\oplus\C{H}_p$ for each $p\in P$. Now considering the projection 
of $T_p(P)$ onto $\C{G}$ with kernel $\C{H}_p$ for each $p\in P$ we get a certain 
$\C{G}$-valued differential $1$-form $\Omega=\Omega(\C{H})$, called  the {\em connection
form} of $\C{H}$.
\end{Emp}

\begin{Def} \label{D:curv}
Let $\C{H}$ be a connection on a $\bold{G}$-torsor $P$, whose connection form is 
$\Omega$. Let $h$ be the natural projection of $T_p(P)$ on $\C{H}_p$ for all 
$p\in P$.
The {\em curvature} of the connection $\C{H}$  is the $2$-form $D\Omega$ defined 
by $\langle X\wedge Y|D\Omega\rangle:=\langle h(X)\wedge h(Y)|d\Omega\rangle$.
A connection with zero curvature is called {\em flat}.
\end{Def}

\begin{Rem} \label{R:trivcon}
The trivial torsor $P\cong \bold{G}\times X$ has a natural flat connection, consisting
of vectors, tangent to $X$. We will call such a connection {\em trivial}.
\end{Rem}

\begin{Lem} \label{L:trivcon}
Let $X$ be a simply connected complex manifold, let $\pi:P\to X$ be a $\bold{G}$-torsor, 
and let $\C{H}$ be a flat connection on $P$. Then there exists
a unique decomposition $P\isom\bold{G}\times X$ such that  $\C{H}$ corresponds to the
trivial connection on $\bold{G}\times X$.
\end{Lem}
\begin{pf}
By \cite[Ch. VII, Thm. 1.1 and 1.2]{St}, there exists a unique $\bold{G}$-equivariant diffeomorphism 
$\varphi:P\isom\bold{G}\times X$ over $X$ which maps $\C{H}$ to the
trivial connection. Hence $\varphi$ induces complex isomorphism between tangent spaces
$T_p(P)=\C{G}\oplus\C{H}_p$ and $T_{\varphi(p)}(\bold{G}\times X)=\C{G}\oplus T_{\pi(p)}(X)$ for each 
$p\in P$. In other words, both $\varphi$ and  $\varphi^{-1}$ are almost
complex mappings between complex manifolds. \cite[Ch. VIII, p.284]{He} then implies that $\varphi$ is
biholomorphic.
\end{pf}

\section{First Main Theorem}

\subsection{Basic examples} \label{SS:basex}

\begin{Def} \label{D:inv}
Let $K/k$ be a quadratic field extension, and let $D$ be a central simple
algebra over $K$. We say that  $\al:D\to D$ is an {\em involution of the 
second kind over $k$} if $\al(d_1+d_2)=\al(d_1)+\al(d_2),\:\al(d_1 d_2)=
\al(d_2)\al(d_1)$ for all $d_1,d_2\in D$, and the restriction of $\al$ to $K$ is
the conjugation over $k$.
\end{Def}

\begin{Not} \label{N:gu}
For $k,D$ and $\al$ as in \rd{inv}, let $\bold{G}=\GU{}(D,\al)$ be the algebraic group over $k$ 
of unitary similitudes, that is $\bold{G}(R)=\{d\in (D\otimes_k R)\m|d\al(d)\in R\m\}$
for each $k$-algebra $R$. Define  the {\em similitudes homomorphism} $\bold{G}\to\B{G}_m$
by $x\mapsto x\al(x)$. Notice also that by the Skolem-Noether theorem the group $\bold{G}$ 
satisfies $\bold{PG}(L)=\bold{G}(L)/Z(\bold{G}(L))$ for every field extension $L$ of $k$.
\end{Not}

\begin{Emp} \label{E:ex1}
{\bf First basic example.} Let $F$ be a totally real field of degree $g$ over $\B{Q}$, 
let $K$ be a totally imaginary quadratic extension on $F$. 
Let $D$ be a central simple algebra of dimension $d^2$ over $K$ with an involution of the second kind $\al$ over $F$. 
Set $\bold{G}:=\GU{}(D,\al)$, and put $D_u:=D\otimes_K K_u$ for each prime $u$ of $K$. 
Let $v$ be a (non-archimedean) prime of $F$ that splits in $K$, and 
let $w$ and $\bar{w}$ be the primes of $K$ that lie over $v$. 
Then $D\otimes_{F}F_v\cong D_w\oplus D_{\bar{w}}$,
and the projection to the first factor together with the similitude
homomorphism induce an isomorphism $\bold{G}(F_v)\isom\mDw\times F_{v}\m$.
We identify $\bold{G}(F_v)$ with $\mDw\times F_{v}\m$ by this isomorphism.

Suppose that $D_w\cong\Mat_d(K_w)$. Identifying $D_w$ with $\Mat_d(K_w)$
by some isomorphism we identify $\bold{G}(F_v)$ with $\GL{d}(K_w)\times F_{v}\m$.
Suppose that $\al$ is positive definite, that is $\bold{G}(F_{\infty_i})
\cong \GU{d}(\B{R})$ for all archimedean completions $F_{\infty_i}\cong\B{R}$
of $F$. Put $E':=F_v\m\times \bold{G}(\afv)$, then $E'$ is a noncompact locally
profinite group. Set $\Gm:=\bold{G}(F)\subset\bold{G}(\af)=\GL{d}(K_w)\times E'$, 
embedded diagonally.

\begin{Prop} \label{P:test1}
The subgroup $\Gm\subset \bold{G}(\af)=\GL{d}(K_w)\times E'$ satisfies the
assumptions of Construction \ref{C:con2}.
\end{Prop}

\begin{pf}
a) is trivial.\par
b) is true, because the closure of $Z(\Gm)\cong K\m$ is cocompact in $Z(\bold{G}(\af))\cong 
(\B{A}_K^f)\m$. \par

c) Since $PE'=\PG(\afv)$ and $\PG(F_v)\cong \PGL{d}(K_w)$,
we have to show that $P\Gm(=\PG(F))$ is a cocompact lattice in
$\PG(\af)$.

\begin{Lem} \label{L:lat}
 If $\bold{H}$ is an $F$-anisotropic group, then $\bold{H}(F)$ is a cocompact lattice in 
$\bold{H}(\B{A}_{F})$. 
\end{Lem}
\begin{pf}
see \cite[Thm. 5.5]{PR}
\end{pf}

Since $\PG$ is anisotropic over each $F_{\infty_i}$, it is 
anisotropic over $F$. Hence by the lemma, $\PG(F)$ is a cocompact lattice in
$ \PG(\B{A}_{F})$. The compactness of the $\PG(F_{\be_i})$'s implies also that 
the projection of $\PG(F)$ to $\PG(\af)$ is a cocompact lattice as well 
(see \cite[Prop. 1.10]{Shi}).
Observe also that the projection $\PG(F)\to \PG(F_v)\cong \PGL{d}(K_w)$ 
is injective.

d) Since $Z(\Gm)\cong K\m$ and $Z(\bold{G}(\af))\cong(\B{A}_K^{f})\m$,
we have to show that the intersection of $\overline{K}\m\subset 
(\B{A}_K^f)\m$ 
with $K_w\m\times\{1\}$ is trivial. This can be shown either
by the direct computation or using the relation between global and local Artin 
maps (see \cite[Ch. VII, Prop. 6.2]{CF}).
\end{pf}

Fix  a central skew field $\wt{D}_w$ over $K_w$ with invariant $1/d$.
Set $E=:\wt{D}\m_w\times E'$, then Construction \ref{C:con2} 
gives us an \esc{E}{K_w} $X$ corresponding to  $\Gm$.
\end{Emp}

\begin{Emp} \label{E:ex2}
{\bf Second basic example.} By Brauer-Hasse-Noether theorem (see \cite[Ch. XIII, $\S$6]{We1})
there exists a unique central skew field $D^{int}$ over $K$ which is locally isomorphic to 
$D$ at all places of $K$ except $w$ and 
$\bar{w}$ and has Brauer invariant $1/d$ at $w$. 
By Landherr theorem (see \cite[Ch. 10, Thm. 2.4]{Sc}), $D^{int}$ admits 
an involution of the second kind over $F$. Fix an embedding $\be_1:K\hra\B{C}$. It induces an
archimedean completion $F_{\be_1}$ of $F$, and we have the following 

\begin{Prop} \label{P:invol}
a) There exists an involution of the second kind $\al^{int}$ of $D^{int}$
over $F$ such that:\par
\quad i) the pairs $(D,\al)\otimes_F F_u$ and $(D\ii,\al\ii)\otimes_F F_u$ are isomorphic at all places $u$ of $F$,
except $v$ and  ${\be_1}$;

\quad ii) the signature of $(D\ii,\al\ii)$ at $\be_1$ is $(d-1,1)$.

b) The group $\bold{G\ii}:=\GU{}(D^{int},\al^{int})$ is determined uniquely (up to an 
isomorphism) by conditions i),ii) of a). 
\end{Prop}

\begin{pf} 
a) follows from \cite[(2.2) and the discussion around it]{Cl} as in \cite[Prop. 2.3]{Cl}.
b) follows immediately from \cite[Ch. 10, Thm. 6.1]{Sc}.
\end{pf}

Let $\bold{G\ii}$ be as in the proposition. Then embedding $\be_1$ defines an isomorphism 
$D\ii\otimes_K K_{\be_1}\isom\Mat_d(\B{C})$, and we identify $\bold{PG\ii}(F_{\infty_1})$ with 
$\PGU{d-1,1}(\B{R})$ by the induced isomorphism. Set $\bold{G\ii}(F)_+:=\bold{G\ii}(F)\cap\bold{G\ii}(F_{\be_1})^0$.
Then $\bold{G\ii}(F)_+=\bold{G\ii}(F)$ if $d>2$, and $[\bold{G\ii}(F):\bold{G\ii}(F)_+]=2$ if $d=2$.
Set $E\ii:=\bold{G\ii}(\af)$, and  
let $E_0\ii\subset E\ii$ be the closure of $Z(\bold{G\ii}(F))\subset E\ii$.
Embed diagonally $\bold{G\ii}(F)$ into $\bold{G\ii}(F_{\infty_1})\times E\ii$ and define
$\Gm\ii$  to be the image of $\bold{G\ii}(F)_+$ under the natural projection to
$$\bold{PG\ii}(F_{\infty_1})\times(E\ii/E_0\ii)=\PGU{d-1,1}(\B{R})\times (E\ii/E_0\ii).$$

\begin{Prop} \label{P:test2}
 The subgroup $\Gm\ii$ is a cocompact lattice in $ \pgr^0\times (E\ii/E_0\ii)$,
and it has an injective projection to the first factor.
\end{Prop}

\begin{pf}
 Notice that the natural projection $E\ii/E_0\ii\to E\ii/Z(E\ii)=PE\ii$ induces
an isomorphism $\Gm\ii\isom\bold{PG\ii}(F)_+\subset\pgr^0\times PE\ii$
and that the group $Z(E\ii)/E_0\ii\cong(\B{A}_K^f)\m/\overline{K}\m$ is compact.
Therefore it will suffice to prove that  $ \bold{PG\ii} (F)$ is a cocompact 
lattice with an injective projection to the first factor of 
$\bold{PG\ii}(F_{\infty_1})\times PE\ii$.
This can be proved by exactly the same considerations as Proposition \ref{P:test1}, c).
\end{pf}

By the proposition, $\Gm\ii$ satisfies the assumptions of Construction \ref{C:con1}, 
so it determines an \esc{E\ii/E_0\ii}{\B{C}} $\wt{X}\ii$, 
which can be regarded as an \esc{E\ii}{\B{C}} with a trivial action of $E_0\ii$.
\end{Emp}

\begin{Rem} \label{R:con2}
For each $S\in\ff{E\ii}$ we have the following isomorphisms
\begin{align}
(\wt{X}\ii_{S})^{an} & \cong S\bs[B^{d-1}\times (E\ii/E_0\ii)]/\Gm\ii\notag\\ 
                     & \cong(S\cdot\overline{Z(\bold{G\ii}(F)})\bs[B^{d-1}\times 
                       \bold{G\ii}(\af)]/\bold{G\ii}(F)_{+}\notag\\
                     & \cong(S\cdot Z(\bold{G\ii}(F))\bs[B^{d-1}\times \bold{G\ii}
                       (\af)]/\bold{G\ii}(F)_{+}\notag\\
            & \cong S\bs[B^{d-1}\times \bold{G\ii}(\af)]/\bold{G\ii}(F)_{+}.\notag
\end{align}
\end{Rem}

\subsection{First Main Theorem} \label{SS:fmt}

\begin{Def} \label{D:adm}
 An isomorphism $\Phi :E\isom E\ii$ is called {\em admissible} if 
it is a product of $\bold{G}(\afv)\isom\bold{G\ii}(\afv)$, induced by some 
$\afv$-linear algebra isomorphism $D\otimes_F \afv\isom D\ii\otimes_F\afv$ 
(compare \rp{invol}), and the composition map $\wt{D}_w\m\times F_{v}\m\isom(D_{w}\ii)\m\times F_{v}\m
\isom\bold{G\ii}(F_v)$, constructed from some algebra isomorphism $\wt{D}_w\isom D\ii\otimes_K K_w$ as in \ref{E:ex1}.
\end{Def}

\begin{Emp} \label{E:fix}
Fix a field isomorphism $\B{C}\isom\B{C}_p$, whose composition with embedding $\be_1:K\hra\B{C}$ (chosen in \re{ex2})
is the natural embedding $K\hra K_w\hra\B{C}_p$.
Identifying $\B{C}$ with $\B{C}_p$ by means of this isomorphism we can view,  
in particular, $K_w$ as a subfield of $\B{C}$.
\end{Emp}

\begin{Mainf} \label{M:f}
 For some admissible isomorphism $\Phi :E\isom E\ii$ there exists a
$\Phi$-equivariant isomorphism $f_{\Phi}$ from the \esc{E}{\B{C}} $X_{\B{C}}$
to the \esc{E\ii}{\B{C}} $\wt{X}\ii$.
\end{Mainf}

\begin{Emp} \label{E:def}
Let $E_0$ be the kernel of the action of $E$ on $X$, and put $\wt{E}:=E/E_0$.
By Corollary \ref{C:rel} there exists a subgroup 
$\Dt\subset\pgr^0\times \wt{E}$ such that the \esc{\wt{E}}{\B{C}} 
$X_{\B{C}}$ corresponds to $\Dt$ by the real case of Construction \ref{C:con1}. 
By \rp{prop2}, each admissible isomorphism $\Phi:E\isom E\ii$ satisfies $\Phi(E_0)=E_0\ii$.
Hence $\Phi$ induces an isomorphism $\bar{\Phi}:\wt{E}\isom E\ii/E_0\ii$.
\end{Emp}

\begin{Thm} \label{T:complex}
 There exists an admissible isomorphism $\Phi:E\isom E\ii$ and an
inner automorphism $\varphi$ of $\PGU{d-1,1}$ such that 
$(\varphi\times\bar{\Phi})(\Dt)=\Gm\ii$.
\end{Thm}

\begin{Lem} \label{L:comp}
 \rt{complex} implies the First Main Theorem.
\end{Lem}

\begin{pf}
\rt{complex} implies that there exists a $\Phi$-equivariant analytic 
isomorphism $\wt{f}_{\Phi}:(X_{\B{C}})\an\isom (\wt{X}\ii)\an$. From the 
$\Phi$-equivariance we obtain analytic isomorphism
$\wt{f}_{\Phi,S}:(X_{S,\B{C}})\an\isom (\wt{X}\ii_{\Phi(S)})\an$ for each $\SIE$.
\rco{GAGA1} provides us with an algebraic isomorphism
$f_{\Phi,S}:X_{S,\B{C}}\isom\wt{X}\ii_{\Phi(S)}$ satisfying $(f_{\Phi,S})\an\cong
\wt{f}_{\Phi,S}$. Taking their inverse limit we obtain a 
$\Phi$-equivariant isomorphism $f_{\Phi}:=\underset{\underset{S}{\longleftarrow}}
{\lim}\,f_{\Phi,S}:X_{\B{C}}\isom\wt{X}\ii$.
\end{pf}

Thus we have reduced our  First Main Theorem to a purely group-theoretic
statement. For proving it we need to know more information about $\Dt$. First we
introduce some auxiliary notation.

\begin{Emp} \label{E:supp}
Let $\Dt'\subset\pgr^0\times PE$ and $\Dt''\subset\pgr^0\times PE'$ be the images of $\Dt$ 
under the natural projections.
Since the groups $E_0\bs Z(E)$ and $E_0\bs \wt{D}_w\m\cdot Z(E)$ are compact,
Lemma \ref{L:conn} shows that subgroups $\Dt'$ and $\Dt''$ correspond by the real case of
Construction \ref{C:con1} to the \esc{PE}{\B{C}} $X'_{\B{C}}:=Z(E)\bs X_{\B{C}}$
and to the \esc{PE'}{\B{C}} $X''_{\B{C}}:=(\wt{D}_w\m\times Z(E))\bs X_{\B{C}}$ respectively.
 The same lemma implies also that the natural projections $\Dt\to\Dt'$ and
$\Dt\to\Dt''$ are isomorphisms.

 Let $E'_0$ be the image of $E_0$ under the canonical projection to $E'$.
Let $\Gm'$ be the image of $\Gm$ under the projection $\gp\times E'\to\gp\times 
( E'_0\bs E')$. Then by \rp{prop2} c), the group $\Gm'$ corresponds by
the $p$-adic case of \rc{con1} to the \esc{E'_0\bs E'}{K_w} $X''':=\wt{D}_w\m\bs X$.
Recall also that by \rp{prop2} d) the \esc{PE'}{K_w} $X''=(\wt{D}_w\m\times Z(E))\bs X$ is 
obtained from the subgroup $P\Gm\subset\pgp\times PE'$ by the $p$-adic case of 
Construction \ref{C:con1}.

For each subset $\Theta$ of $\Dt$, $\Dt'$ or $\Dt''$ (resp. of $\Gm$, $\Gm'$ or $P\Gm$)
we denote by $\Theta_{\be}$ (resp. $\Theta_{G}$) and $\Theta_{E}$ its
projections to the first and to the second factors respectively (compare \ref{N:index}).   
\end{Emp} 

Our next task is to establish the connection between $\Dt$ and $\Gm$.
The next key proposition is the modifications of \cite[Prop. 2.6]{Ch2}. In 
it we apply Ihara's technique of elliptic elements to relate elements in $\Dt$ 
and in $\Gm$.

\begin{Prop} \label{P:corresp}
For each $\dt\in\Dt$ with elliptic projection $\dt_{\be}\in\pgr^0$, there exist 
$\gm\in\Gm$ and $\gm_D\in\wt{D}_w\m$ with $(\gmG,\gm_D)\in\gp\times\wt{D}_w\m$ 
elliptic (with respect to its action on $\Si$) and a
 representative $\wt{\dt}=(\wt{\dt}_{\be},\wt{\dt}_E)\in\gr^0\times E$ of $\dt$ satisfying the following 
conditions:\par
a) the elements $(\gm_D,\gm_E)$ and $\wt{\dt}_E$ are conjugate in $E$;\par

b) the characteristic polynomials of $\wt{\dt}_{\be}$ and $\gm_{G}$
are equal.\par

 Conversely, for each $\gm\in\Gm$ and $\gm_D\in\wt{D}_w\m$ with
$(\gmG,\gm_D)\in\gp\times\wt{D}_w\m$ elliptic, there exist $\dt\in\Dt$ with 
elliptic projection 
$\dt_{\be}\in\pgr^0$ and a representative $\wt{\dt}\in\gr^0\times E$ of $\dt$
 satisfying conditions a) and b).
\end{Prop}

\begin{pf}
 If an element $\dt_{\be}\in\Dt_{\be}$ is elliptic, then $\dt_{\be}$ has a fixed 
elliptic point $P$ on $B^{d-1}$. The action of $\dt_{\be}$ on $B^{d-1}$ coincides 
with the action
of $\dtE$ on $B^{d-1}\cong B^{d-1}\times\{1\}\subset(B^{d-1}\times\wt{E})/\Dt\cong
(X_{\B{C}})\an$, therefore $P$, viewed as a point of $(X_{\B{C}})\an$ (or of
$X(\B{C})$), is an elliptic point of $\dtE$.
Using the isomorphism $\B{C}\isom\B{C}_p$, chosen above, $P$ can be considered 
as a point of $X(\B{C}_p)$, hence
 as a point of the $p$-adic pro-analytic space $X\an$. There exists 
an element $g\in E$ such that the point $P':=g(P)$ lies in 
$\C{Y}:=\rho_1(\Si)$ in the notation of \rco{etcom}.

Let $\pi$ be the natural projection $X\to X'''$. Choose a representative 
$\wt{g}\in E$ of $g\dtE g^{-1}\in \wt{E}$. Since $\wt{g}$ fixes $P'$, 
it fixes the projection $P'':=\pi(P')\in(X'''_{\B{C}_p})\an$. Hence $\wt{g}$ stabilizes the connected 
component $\Om\Hat{\otimes}_{K_w}\B{C}_p\times\{1\}\subset(X'''_{\B{C}_p})\an$ 
containing $P''$.
By \rp{prop2} c), the image of $\wt{g}$ under the canonical projection 
$E\to E'_0\bs E'$ belongs to the projection of $\Gm'$ to 
$E'_0\bs E'$. We can therefore choose $\gm\in\Gm$ whose projection to $E'_0\bs E'$
coincides with that of $\wt{g}\in E$. Therefore $\wt{g}\cdot\gmE^{-1}$ belongs to $\wt{D}_w\m\times E'_0=
\wt{D}_w\m\cdot E_0$. Hence there exists a $\gm_D\in \wt{D}_w\m$ such that
$\wt{g}(\gm_D^{-1},\gmE^{-1})\in E_0$. It follows that $(\gm_D,\gmE)\in E$ is also a 
representative of $g\dtE g^{-1}$.

The action of $(\gm_D,\gmE)$ on the tangent space of $P'\in\C{Y}$ is conjugate to 
the action of $\dtE$ on the tangent space of $P$, therefore $P'$ is an elliptic 
point of $(\gm_D,\gmE)$. Since $\rho_1$ is \'{e}tale, one-to-one (use \rco{etcom})
and $\wt{D}_w\times\Gm$-equivariant, the action of  $(\gm_D,\gmE)$ on the tangent
space of $P'\in\C{Y}$ coincides with the action of $(\gmG,\gm_D)$ on the tangent 
space of $\rho_1^{-1}(P')\in\Si$. Therefore $\rho_1^{-1}(P')$ is an elliptic
point of  $(\gmG,\gm_D)$.
It follows that the action of $(\gmG,\gm_D)$ on the tangent space of 
$\rho^{-1}(P')\in\Si$ is 
conjugate to the action of $\dt_{\be}$ on the tangent space of $P\in B^{d-1}$. 
Using the \'etalness of the projection $\Si\to\Om$ we conclude from \rl{lintr}
that there exists a representative 
$\dt'_{\be}\in\gr^0$ of $\dt_{\be}$ such that the characteristic polynomials of
$\dt'_{\be}$ and $\gmG$ are equal. 
Hence $\wt{\dt}:=(\dt'_{\be}, g^{-1}(\gm_D,\gmE)g)$ is the required representative of $\dt$.

The proof of the opposite direction is very similar, but much easier technically.
If an element $(\gmG,\gm_D)\in\GmG\times\wt{D}_w\m$ is elliptic, then it has an 
elliptic point $Q\in\Si$. Hence $Q':=\rho_1(Q)\in X\an$ is an 
elliptic point of $(\gm_D,\gmE)\in E$. Hence $Q'$ can be considered as a point of
the complex analytic space 
$(X_{\B{C}})\an\cong (B^{d-1}\times \wt{E})/\Dt$. Choose a representative 
$(x,g)\in B^{d-1}\times E$ of $Q'$. Then the element $g(\gm_D,\gmE)g^{-1}\in E$
fixes $Q'':=g(Q')\in B^{d-1}\times\{1\}$, hence it stabilizes the connected component
$B^{d-1}\times\{1\}\subset(X_{\B{C}})\an$. It follows that the image of  
$g(\gm_D,\gmE)g^{-1}$ under the projection of $E$ to $\wt{E}$ belongs to 
$\DtE$. The rest of the proof is exactly the same as in the other direction.
\end{pf}

\begin{Cor} \label{C:ellel}
 For each $\dt\in\Dt$ with elliptic projection $\dt_{\be}\in\pgr^0$, there exists a 
representative 
$$\wt{\dt}=(\wt{\dt}_{\be},\wt{\dt}_v,\wt{f}_v,\wt{\dt}^{f;v})\in
\gr^0\times\wt{D}_w\m\times F_v\m\times \bold{G}(\afv)$$
such that\par
a) if we view $K$ as a subset of $\B{C}$, of $K_w$ and of $K\otimes_F\afv$ respectively, 
then the characteristic polynomials of $\wt{\dt}_{\be},\wt{\dt}_v$ and 
$\wt{\dt}^{f;v}$ have their coefficients in $K$ and coincide;\par

b) $\wt{f}_v$ and the similitude factor of $\wt{\dt}^{f;v}$ belong to $F$, viewed 
as a subset of $F_v\m$ and of $(\afv)\m$ respectively, and coincide.
\end{Cor}

\begin{pf}
 Take $\gm$ and $\wt{\dt}$ as in the proposition. Then the statement follows from
\rp{padell}. 
\end{pf}


\begin{Prop} \label{P:closure}
 The closure $\overline{\DtE'}$ contains $(S\wt{D}_w\m\cap T_1)\times P(\bold{G^{der}}(\afv))$.
\end{Prop}

\begin{pf}
Let $X'_0$ be the connected component of $X'_{\B{C}}$ such that $(X'_0)^{an}
\supset B^{d-1}\times \{1\}$. Then by \rp{prop1} c), $\overline{\DtE'}=
Stab_{PE'}(X'_0)$. \rp{padact} implies that the group 
$S\wt{D}_w\m\cap T_1$ acts trivially on the set of connected components of 
$X'_{\B{C}}$, therefore it remains to show only that 
$\overline{\DtE''}\supset P(\bold{G^{der}}(\afv))$. To prove it
we first observe that by the strong approximation theorem (see, for example, 
\cite[Ch. II, Thm. 6.8]{Ma}), the closure $\overline{\PG(F)}$ of $\PG(F)$ in $PE'=
\PG(\afv)$ contains $P(\bold{G^{der}}
(\afv))$. So the proposition follows from

\begin{Lem} \label{L:closure}
We have $\overline{\DtE''}=\overline{\PG(F)}$.
\end{Lem}

\begin{pf}
By \rp{prop1} c), we see that $\overline{\DtE''}$ is the stabilizer of the connected component 
$Y_{\be}$ of
$X''_{\B{C}}$ such that $(Y_{\be})^{an}\supset B^{d-1}\times \{1\}$ and that 
$\overline{\PG(F)}$  is the stabilizer of the connected component $Y_p$ of $X''_{\B{C}_p}$ such that
$(Y_p)^{an}\supset\Om\times \{1\}$.  Since the group $PE'$ acts transitively on
the set of geometrically connected components of $X''$, the subgroups  $\overline{\DtE''}$
and $\overline{\PG(F)}$ are conjugate in $PE'$. Since $\overline{\DtE''}$ contains
$P(\bold{G^{der}}(\afv))$, it is normal. So we are done.
\end{pf}\end{pf}

\subsection{Computation of Q(TrAd)} \label{SS:TrAd}

In the next subsection a field $\B{Q}(\Tr\Ad)$ (generated by the traces of 
the adjoint representation) will be a field of definition of a
certain algebraic group.

\begin{Rem} \label{R:trad}
 If $g\in \GL{d}$, then by  direct computation we obtain that 
$\Tr\Ad\,g=\Tr\,g\cdot\Tr(g^{-1})$. Hence for $g\in\PGL{d}$ we have
$\Tr\Ad\,g=\Tr\,\wt{g}\cdot\Tr(\wt{g}^{-1})-1$ for each representative
$\wt{g}\in\GL{d}$ of $g$. 
%
\end{Rem}

\begin{Prop} \label{P:trad}
We have $\B{Q}(\Tr\Ad\,\Dt_{\infty})=F\overset{\be_1}{\hra}\B{R}$.
\end{Prop}

\begin{pf}
 It follows from  \rp{corresp}, \rp{padell} and \rr{trad} that
$$\B{Q}(\Tr\Ad\,\dt_{\infty}|\,\dt_{\infty}\in\Dt_{\infty}\text { is elliptic})=
\B{Q}(\Tr\Ad\,\gmG|\,\gmG\in P\GmG\subset\pgp\text { is elliptic }).$$ Let $F'$
be the last-named field. Then $F'\subset F$, since $P\Gm=\PG(F)$ and since $\PG$ is an
algebraic group defined over $F$. It follows from the weak approximation theorem 
that for each non-archimedean prime $u\neq v$ of $F$, the closure of 
the projection to $\PG(F_u)$ of the set $\{\gm\in P\Gm|\,\gmG$ is elliptic$\}$
contains an open non-empty subset of $\PG(F_u)$. (Recall that the closure of 
$P\GmG$ in $\pgp$ contains $\psp$ by \rp{basic}, and that
the set of elliptic elements of $\psp$ is open and non-empty by \rp{ell}.) 
Therefore $F'$ is dense in each non-archimedean completion
$F_u$ of $F$ for $u\neq v$. Thus $F'$ splits completely in $F$ at almost all
places. Hence $F'=F$ (see \cite[Ch. VII, $\S$4, Thm. 9]{La}).
This part of the proof is completely identical with Cherednik's proof of
\cite[Prop. 2.7]{Ch2}. 

Now we want to prove that $\B{Q}(\Tr\Ad\,\Dt_{\infty})=\B{Q}(\Tr\Ad\,\dt_{\infty}
|\,\dt_{\infty}\in\Dt_{\infty}\text { is elliptic })$. Since the group 
$\PGU{d-1,1}$ is absolutely simple, the representation 
\qquad\qquad\qquad\linebreak
$\Ad:\pgr\to \GL{}(\Lie(\pgr))\cong \GL{d^2-1}(\B{R})$ is absolutely irreducible. 
Therefore our statement is a consequence of the following general

\begin{Lem} \label{L:eqtrad}
Let $\rho$ be an absolutely irreducible algebraic representation of $\PGU{d-1,1}$,
and let $\wt{\Dt}$ be a dense subgroup of $\pgr^0$.
Then $\B{Q}(\Tr(\rho(\wt{\Dt})))=\B{Q}(\Tr\rho(\dt)|\,\dt\in\wt{\Dt}\text { is 
elliptic })$.
\end{Lem}

\begin{pf}
Let $L$ be the last-named field. If $g\in\pgr^0$ is elliptic and $g^r$ is not 
elliptic for some $r\in\B{Z}\sm\{0\}$, then by \rl{lintr}, $g$ belongs to some 
Zariski closed proper subset of $\PGU{d-1,1}$. Therefore for each $N\in\B{N}$, there 
exists an open subset $W\subset\pgr^0$ such that for each $g\in W$ and each $r\in\B{Z}$, 
satisfying $1\leq |r|\leq N$, the element $g^r$ is elliptic. Choose  $g\in W$.
By the continuity of multiplication, there exists an open 
neighbourhood $U\subset W$ of $g$ such that for each $g_1,...,g_k\in U$ and each 
$n_1,...,n_k\in\B{Z}$, 
satisfying $n_1+...+n_k\neq 0 ,\: |n_1|+...+|n_k|\leq N$, the element
$g_1^{n_1}\cdot ...\cdot g_k^{n_k}$ is elliptic. Take $N=6m^2$, where $m$ is the dimension of $\rho$.

Since $\pgr^0$ is a connected real Lie group, it is generated by $U$. The subgroup 
$\wt{\Dt}$ is dense in $\pgr^0$ by 
\rp{basic}, therefore  $\wt{\Dt}\cap U$ generates the group $\wt{\Dt}$ 
(see \cite[Ch. IX, Lem. 3.3]{Ma}).
Since the restriction of $\rho$ to the Zariski dense subgroup $\wt{\Dt}$ is 
absolutely irreducible, Burnside's theorem (see \cite[vol. II, Ch. XVII, 130]{Wa})
implies that $\C{D}:=\dim_{\B{R}}(\Span_{\B{R}}(\rho(\wt{\Dt})))=m^2$.

Set $\wt{\Dt}^0:=\{1\}\subset \wt{\Dt}$, and for each positive integer $n$ set
$\wt{\Dt}^n:=\{g_1^{n_1}\cdot ...\cdot g_k^{n_k}|g_i\in\Dt_{\be}\cap U,
\:|n_1|+...+|n_k|\leq n\}\subset \wt{\Dt}$. Denote 
$\dim_{\B{R}}(\Span_{\B{R}}(\rho(\wt{\Dt}^n)))$ by $\C{D}_n$. Since $\wt{\Dt}=
\underset{n}{\bigcup}\wt{\Dt}^n$, we have $1=\C{D}_0\leq \C{D}_1\leq ...\leq \C{D}_n\leq ...\leq\C{D}
=\underset{n}\sup\,\C{D}_n$. Moreover, if $\C{D}_n=\C{D}_{n+1}$ for some $n$, 
then $\C{D}_n=\C{D}_{n+1}=...=\C{D}$. Therefore $\C{D}_{m^2-1}=m^2$.
 Hence there exist elements $\bar{\dt}_i\in\wt{\Dt}^{m^2-1},\; i=1,...,
m^2$ such that $\{\rho(\bar{\dt}_i)\}_i$ constitute a basis for
$\Mat_{m^2}(\B{R})$. Choose any $g\in\wt{\Dt}\cap U$ and take $\dt_i:=
g^{m^2+1}
\bar{\dt}_i$. Then $\{\rho(\dt_i)\}_i$ still constitute a basis for
$\Mat_{m^2}(\B{R})$. Each $\dt_i$ is of the form $g_1^{n_1}\cdot ...\cdot g_k^{n_k}$,
where the $g_i$'s belong to $\wt{\Dt}\cap U$ and the $n_i$'s satisfy $n_1+...+n_k\geq 2$ and 
$|n_1|+|n_2|+...+|n_k|\leq 2m^2$. In particular, each $\dt_i$ is 
elliptic, therefore $\Tr\rho(\dt_i)\in L$.

\begin{Lem} \label{L:claim}
 If for some $\dt\in\pgr^0$ the elements $\dt\dt_i$ are elliptic for all 
$i=1,...,m^2$, then $\rho(\dt)$ can be written as a linear combination 
of the $\rho(\dt_i)$'s with coefficients in $L$.
\end{Lem}

\begin{pf}
 Let $e_1,...,e_{m^2}$ be the dual basis of $\{\rho(\dt_i)\}_i$ relative to 
the bilinear
form $(x,y)\mapsto \Tr(xy)$. If $\dt$ is as in the lemma, then 
$\Tr\rho(\dt\dt_i)=\Tr(\rho(\dt)\rho(\dt_i))\in L$ for all $i=1,...m^2$.
Hence $\rho(\dt)$ can be written as a linear combination  of the $e_i$'s with 
coefficients in $L$. Therefore it is enough
 to prove that each $e_i$ can be written as a linear combination of the $\rho(\dt_i)$'s
with coefficients in $L$. The last condition is equivalent to the condition
that each $ \rho(\dt_i)$ can be written as a linear combination 
of the $e_i$'s with coefficients in $L$. Thus, as we mentioned above, to complete the
proof it is enough to show that each $\dt_i$ satisfies the conditions of the
lemma. This follows directly from the definition of the $\dt_i$'s and of $U$.  
\end{pf}

The choice of the $\dt_i$'s assures that for every $\dt\in\wt{\Dt}\cap(U\cup U^{-1})$
the elements $\dt\dt_i$ are elliptic for all 
$i=1,...,m^2$. Therefore the above lemma implies that $\rho(\dt)$ can
 be written as a linear combination of the $\rho(\dt_i)$'s with coefficients in $L$.
The set $U\cap\wt{\Dt}$ generates the group $\wt{\Dt}$, hence for every  $\dt\in\wt{\Dt}$,
the linear transformation $\rho(\dt)$ can be written as a polynomial
in  the $\rho(\dt_i)$'s with coefficients in $L$. For each $i,j,k\in\{1,...,
m^2\}$ the elements $\dt_i\dt_j\dt_k$ are elliptic, therefore by the lemma each
$\rho(\dt_i\dt_j)=\rho(\dt_i)\rho(\dt_j)$ can be written as a linear combination
 of the $\rho(\dt_k)$'s with coefficients in $L$. Hence every polynomial in the $\rho(\dt_i)$'s with
 coefficients in $L$,  can be written as a linear combination
 of the $\rho(\dt_i)$'s with coefficients in $L$. In particular, this is true for each $\rho(\dt)$ with 
$\dt\in\wt{\Dt}$. Hence $\B{Q}(\Tr\rho(\wt{\Dt}))\subset L$.
\end{pf}\end{pf}

\begin{Cor} \label{C:trad}
Suppose that a subgroup $\wt{\Dt}\subset\Dt_{\infty}$ is Zariski dense in
$\PGU{d-1,1}$ and that $\Dt_{\infty}\subset Comm_{\pgr}(\wt{\Dt})$.
Then $\B{Q}(\Tr\Ad\,\wt{\Dt})=\B{Q}(\Tr\Ad\,\Dt_{\infty})(=F)$.
\end{Cor}

\begin{pf}
Set $L:=\B{Q}(\Tr\Ad\,\wt{\Dt})$, then there exists an $L$-form $V$ of 
$\Lie(\pgr)$ preserved by $\Ad\,\wt{\Dt}$ (see \cite[Ch. VIII,
Prop. 3.22]{Ma}). Take any $\dt\in\Dt_{\be}$. Then some subgroup of
finite index $\wt{\Dt}'$ of $\wt{\Dt}$ satisfies $\dt\wt{\Dt}'\dt^{-1}\subset
\wt{\Dt}$, hence $(\Ad\,\dt)(\Ad\,\wt{\Dt}')(\Ad\,\dt)^{-1}(V)=V$. 

Since the subgroup $\wt{\Dt}'$ is also Zariski dense in $\PGU{d-1,1}$, 
Burnside's theorem implies that $\Ad\,\wt{\Dt}'$ generates $\End\,V$ as an $L$-vector 
space. Therefore \qquad\qquad\qquad\qquad\qquad\qquad\linebreak
$(\Ad\,\dt)(\End\,V)(\Ad\,\dt)^{-1}\subset \End\,V$. In other 
words, $\Ad(\Ad\,\dt)(\End\,V)=\End\,V$. Let $\bold{H}$ be the Zariski closure of $\Ad\,\wt{\Dt}
\subset \GL{}(V)$. Then $\bold{H}$ is an $L$-form of $\Ad\,\PGU{d-1,1}$, hence
$\Lie\,\bold{H}\subset \End\,V$ is an $L$-form of $\Lie(\Ad\,\PGU{d-1,1})$. 
In particular, $ \Lie\,\bold{H}=\End\,V\cap \Lie(\Ad\,\PGU{d-1,1})$, therefore 
$\Ad(\Ad\,\Dt_{\be})(\Lie\,\bold{H})=\Lie\,\bold{H}$. 
Since $\PGU{d-1,1}$ is adjoint, the homomorphism
$\ad:=\Ad_{*}:\Lie\,\PGU{d-1,1}\to \Lie(\Ad\,\PGU{d-1,1})$ is an isomorphism. 
Therefore $\wt{V}:=\ad^{-1}(\Lie\,\bold{H})$ is an $L$-form of 
$\Lie(\PGU{d-1,1})$ and $\Ad\,\Dt_{\be}\subset \GL{}(\wt{V})$. It follows that 
$\B{Q}(\Tr\Ad\,\Dt_{\be})\subset L$.
\end{pf}

\subsection{Proof of arithmeticity} \label{SS:arithm}
\begin{Emp} \label{E:Ma1}
Consider the subgroup $\Dt'\subset\pgr^0\times PE\subset\pgr\times PE$, defined in \re{supp}. For a 
finite place $u$ of $F$ let $\bold{G_u}$ be $\bold{PG_{F_u}}$ for $u\neq v$ and $\PGL{1}(\wt{D}_w)$,
viewed as an algebraic group over $F_v\cong K_w$, for $u=v$.
In what follows it will be also convenient to introduce a formal symbol $\be$ and to
write $F_{\be}$ instead of $\B{R}$ and  $\bold{G_{\be}}$ instead of $\PGU{d-1,1}$ 
(the  algebraic group over $F_{\infty}\cong\B{R}$).

Let $M$ be a finite set of non-archimedean primes of $F$, containing $v$ for 
simplicity of notation. Set $\bar{M}:=M\cup\infty
$ and choose $S\in\ff{\PG(\B{A}_F^{f;M})}$. For each subset $M'$ of $\bar{M}$, denote 
$\underset{u\in M'}{\prod}\bold{G_u}(F_u)$ by $G_{M'}$. Denote also the 
projection of $\Dt'\cap (G_{\bar{M}}\times S)$ to $G_{\bar{M}}$ by $\Dt^S$.
Let $\Dt^{S}_{\infty}$ (resp. $\DtE^S$) be the projection of $\Dt^S$ to 
$\bold{G_{\be}}(F_{\infty})$ (resp. to $G_M$).
For each $u\in\bar{M}$ and each $\dt\in\Dt^S$ denote the projection of $\dt$ to 
$\bold{G_u}(F_u)$ by $\dt_u$.
\end{Emp}

\begin{Def} \label{D:irrlat}
 A lattice $\Gm\subset G_{\bar{M}}$ is called {\em irreducible} if for every proper
 non-empty
subset $M'\subset \bar{M}$ the subgroup $(\Gm\cap G_{M'})(\Gm\cap G_{\bar{M}\sm
M'})$ is of infinite index in $\Gm$ (compare \cite[p.133]{Ma}).
\end{Def}

\begin{Def} \label{D:prQD}
 We say that a lattice $\Gm$ of $G_{\bar{M}}$ has {\em property (QD')} if 
the  closure of $\Gm\cdot \bold{G_{\be}}(F_{\infty})$ in $G_{\bar{M}}$ has a finite index.
\end{Def}

\begin{Rem} \label{R:QD}
Since the group $\PGU{d-1,1}$ is isotropic over $\B{R}$, it follows from 
\cite[p.290, Rem. (v)]{Ma} that if $\Gm$ has property (QD'), 
then it has property (QD) in the sense of Margulis (see \cite[p.289]{Ma}). 
\end{Rem}

\begin{Prop} \label{P:condit}
The subgroup $\Dt^S\subset G_{\bar{M}}$ is a finitely generated cocompact irreducible 
lattice, which is of infinite index in $Comm_{G_{\bar{M}}}(\Dt^S)$
and has property (QD').
\end{Prop}

\begin{pf}
 Observe that $\pgr\times PE=G_{\bar{M}}\times \PG(\B{A}^{f;M}_{F})$ and that  
$\Dt'$ is a cocompact lattice in $G_{\bar{M}}\times \PG(\B{A}^{f;M}_{F})$ having 
injective projection to $\pgr$, hence to $G_{\bar{M}}$. It follows from  
\rl{basic} that $\Dt^S\subset G_{\bar{M}}$ is a cocompact  lattice, which is of 
infinite index in $Comm_{G_{\bar{M}}}(\Dt^S)$. By \rp{closure} the closure of 
$\bold{G_{\be}}(F_{\infty})\cdot\Dt'$ in $G_{\bar{M}}\times 
\PG(\B{A}^{f;M}_{F})$ contains $\bold{G_{\be}}(F_{\infty})\times (S\wt{D}_w\m\cap T_1)\times P(\bold{G^{der}}(\afv))$.
 Hence the 
closure of $\bold{G_{\be}}(F_{\infty})\cdot\Dt'$ in $G_{\bar{M}}$ contains
$\bold{G_{\be}}(F_{\infty})\times(S\wt{D}_w\m\cap T_1)\times\underset{u\in M\sm\{v\}}{\prod}P(\bold{G^{der}}(F_u))$.

 In particular, $\Dt^S$ has property (QD'). Let $M'$ be a non-empty subset of
 $M$. Then $\Dt^S\cap G_{M'}=\{1\}$, because the projection of $\Dt'$ to $\pgr$ is injective. 
Suppose that $\Dt^S$ is not irreducible, then $[\Dt^S:(\Dt^S\cap
G_{\bar{M}\sm M'})]<\infty$. Hence
$$[\overline{\bold{G_{\be}}(F_{\infty})\cdot\Dt^S}:\overline{(\bold{G_{\be}}(F_{\infty}
)\cdot\Dt^S )\cap G_{\bar{M}\sm M'}}]<\infty.$$

Since $ \overline{\bold{G_{\be}}(F_{\infty})\cdot
\Dt^S}\supset \bold{G_{\be}}(F_{\infty})\times(S\wt{D}_w\m\cap T_1)\times\underset{u\in M\sm\{v\}}{\prod}P(\bold{G^{der}}(F_u))$ and\qquad\qquad\qquad\linebreak
$\overline{(\bold{G_{\be}}(F_{\infty})\cdot\Dt^S )\cap G_{\bar{M}\sm M'}}\subset 
G_{\bar{M}\sm M'}$, we get a contradiction. Since $\Dt^S$ is a cocompact lattice in $G_{\bar{M}}$, 
it is finitely generated (see \cite[Ch. IX, 3.1(v)]{Ma}).
\end{pf}

\begin{Emp} \label{E:Step2}
Now we are going to use the results of Margulis (see \cite{Ma}).
By \cite[Ch. VIII, Prop. 3.22]{Ma}, there exists a basis in 
$\Lie(\pgr)$ such that all transformations in $\Ad\,\Dt_{\infty}$ are written
in this basis  as matrices with entries in $\B{Q}(\Tr\Ad\,\Dt_{\be})=F\subset F_{\be_1}\cong\B{R}$. 
Define a rational over 
$\B{R}$ homomorphism $\varphi:\bold{G_{\be}}\to \GL{d^2-1}$ by assigning to
$g\in\bold{G_{\be}}$ the matrix of $\Ad\,g$ in the above basis.
It follows that $\varphi(\Dt_{\infty})\subset \GL{d^2-1}(F)$. 
Let $\bold{H}$ be the Zariski closure of $\varphi(\Dt_{\infty})$;
then $\bold{H}$ is an algebraic group, defined over $F$, and $\varphi(\Dt_{\infty}) \subset \bold{H}(F)$. 
Since $\Dt_{\infty}$ is Zariski 
dense in $\bold{G_{\be}}$ and since the group $\bold{G_{\be}}=\PGU{d-1,1}$ is adjoint,
$\varphi$ induces an isomorphism $\PGU{d-1,1}\isom \bold{H_{F_{\be_1}}}$. In
particular, $\bold{H}$ is an $F$-form of $\PGU{d-1,1}$.

 By \rp{condit}, $\Dt^S$ satisfies the conditions of Theorem (B) of 
\cite[p.298]{Ma}, therefore it is arithmetic in the sense of \cite[p.292]{Ma}.
The group $\Dt_{\be}^S$ is Zariski dense in $\bold{G_{\be}}$ 
(see \cite[Ch. IX, Lem. 2.1]{Ma}), and $\Dt_{\be}\subset 
Comm_{\bold{G_{\be}}(F_{\be})}(\Dt^S_{\be})$. Hence $\varphi(\Dt^S_{\be})$ is
Zariski dense in $\bold{H}$, and $\B{Q}(\Tr\Ad\,\Dt^S_{\be})=F$ (by \rcr{trad}).

It follows (see Margulis' proof \cite[pp.307--311]{Ma})
that the following conditions are satisfied: \par

a) The group $\bold{H}(F_{\infty_i})$ is compact for each $i=2,...,g$.\par

b) There exists a unique bijection $I$ from $M$ to a (finite) set of
non-archimedean primes of $F$ satisfying the following property: for each
$u\in M$ there exists a continuous isomorphism $\omega_u:F_u\isom F_{I(u)}$ and 
an $\omega_u$-algebraic isomorphism $\tau_u:\bold{G_u}\isom \bold{H}$ (that is $\tau_u$ becomes
an isomorphism of algebraic groups over $F_u$ after the identification of $F_u$ 
with $F_{I(u)}$ by means of $\omega_u$) such that $\tau_u(\dt_u)
=\varphi(\dt_{\infty})\in\bold{H}(F)\subset  \bold{H}(F_{I(u)})$ for all $\dt\in\Dt^{S}$. 
Since the subgroup $\Dt^{S}_{u}$
is Zariski dense in $\bold{G_u}$ (see \cite[Ch. IX, Lem. 2.1]{Ma}), $\tau_u$ is unique.\par

c) Let $\tau_M:\prod_{u\in M}\bold{G_u}(F_u)\isom
\prod_{u\in M}\bold{H}(F_{I(u)})$ be the product of the $\tau_u$'s. 
Put $\C{O}_{F,I(M)}:=\{f\in F| f\in \C{O}_{F_u}$ for each finite prime
$u\notin I(M)$ of $F\}$. Then the 
subgroup $\tau(\DtE^{S})\subset \bold{H}(F)$ is commensurable with 
$\bold{H}(\C{O}_{F,I(M)})$.

Taking $M$ larger and larger we conclude from b) that there
exists a unique one-to-one surjective map $I$ of the set of all non-archimedean primes of $F$
into itself such that for each prime $u$ of $F$ there exists a continuous isomorphism $\omega_u:F_u\isom F_{I(u)}$
and a unique $\omega_u$-algebraic isomorphism $\tau_u:\bold{G_u}\isom \bold{H}$
such that $\tau_u(\dt_u)=\varphi(\dt_{\infty})\in \bold{H}(F)\subset 
\bold{H}(F_{I(u)})$ for all $\dt\in\Dt'$.
The maps $\tau_u$ combined together for all non-archimedean 
primes $u$ of $F$ give us a continuous isomorphism $\tau:\prod_u
\bold{G_u}(F_u)\isom \prod_u\bold{H}(F_u)$ such that $\tau(\DtE')\subset \bold{H}(F)\subset\bold{H}(\af)\subset\prod_u\bold{H}(F_u)$.
By c), the subgroup $\tau(\DtE'\cap S)$ is commensurable 
with $\bold{H}(\C{O}_F)$ for each $S\in\ff{PE}$.
\end{Emp}

\subsection{Determination of $\bold{H}$} 

\begin{Emp} \label{E:class}
 Recall that $\bold{H}$ is an $F$-form of $\PGU{d-1,1}$. In particular, 
it is a form of $\PGL{d}$. By the classification of simple algebraic groups (see \cite{Ti}),
there exists a quadratic extension $F'$ of $F$ and a central simple algebra 
$D'$ over $F'$ of dimension $d^2$ (defined up to a replacement $D'\mapsto (D')\op$) with an involution of
the second kind $\al'$ over $F$ such that $\bold{H}\cong\PGU{}(D',\al')$. Moreover,
$F'$ is uniquely determined if $d>2$ and can be chosen arbitrary 
if $d=2$. We denote the group $\GU{}(D',\al')$ by $\bold{G'}$ and 
will not distinguish between $\bold{H}$ and $\bold{PG'}$.
\end{Emp}

\begin{Cl} 
For each non-archimedean prime $u$ of $F$ we have $I(u)=u$ and $\omega_u$ 
is the identity.
\end{Cl}

\begin{pf}
Since the map $\tau_u:\bold{G_u}\isom \bold{H}$ is $\omega_u$-algebraic, 
we have $\Tr\Ad(\tau_u(g))=\omega_u(\Tr\Ad(g))$ for each $g\in \bold{G_u}(F_u)$.
Hence for each $\dt\in\Dt'$ we have
$$\Tr\Ad(\varphi\times\tau)(\dt)=(\Tr\Ad(\dt_{\be});...,\omega_u(\Tr\Ad(\dt_u)),
...)\in (F_{\be_1};...,F_{I(u)},...).$$

Recall that  $(\varphi\times\tau)(\Dt')\subset \bold{H}(F)$, hence
$\Tr\Ad((\varphi\times\tau)\Dt')\subset F$. On the other hand,
\rco{ellel} implies that $\Tr\Ad(\dt)\in F\subset F_{\infty_1}\times\B{A}
_{F}^{f}$  for each $\dt\in\Dt'$ with elliptic $\dt_{\infty}$. In particular,  for such $\dt$'s 
we have $\Tr\Ad(\dt_{\be})=\Tr\Ad(\dt_u)\in F$ for each $u$.  
Since we showed in the proof of \rp{trad} that $\B{Q}(\Tr\Ad(\dt_{\be})| \dt_{\infty}$ is elliptic$)=F$,
we conclude from the above that the restriction of each $\omega_u:F_u\isom F_{I(u)}$
to $F$ is the identity. Since each $\omega_u$ is continuous, the claim follows.
\end{pf}

\begin{Emp}
Next we will show that in the case $d>2$ we have $F'=K$. Indeed, if a prime $u$ of $F$ splits in $K$,
then $\bold{PG'}(F_u)\cong\bold{G_u}(F_u)\cong P\bar{D}_{u}\m$ for some central simple algebra
$\bar{D}_u$ over $F_u$. It follows that $u$ splits in $F'$. By \cite[Ch. VII, $\S$4, Thm. 9]{La},
$F'=K$. As we mentioned before, we may take $F'=K$ also in the case $d=2$.
\end{Emp}
\begin{Prop} \label{P:adeles}
The map $\tau$ induces a continuous isomorphism $PE\isom \bold{H}(\af)$.
\end{Prop}

\begin{pf}
Since $PE\cong P\wt{D}_w\m\times PE'$ and $\bold{H}(\af)\cong\bold{H}(F_v)\times\bold{H}(\afv)$, 
we need only to show that $\tau^v:\prod_{u\neq v}
\bold{G_u}(F_u)\isom \prod_{u\neq v}\bold{H}(F_u)$
 induces a continuous isomorphism $PE'\isom \bold{H}(\afv)$.

First we claim that $\tau^v$ induces a continuous map from $\overline{\Dt'}\subset PE'$ to $\bold{H}(\afv)$.
In fact, let a sequence $\{\dt_n\}_n\subset\Dt'_E$ converges to $g\in PE'$. Then the sequence 
$\{\dt_n\dt_{n+1}^{-1}\}_n$ converges to $1$. Therefore for each $S\in\ff{PE'}$ there 
exists $N_S\in \B{N}$ such that $\dt_n\dt_{n+1}^{-1}\in\DtE''\cap S$ (hence 
$\tau^v(\dt_n\dt_{n+1}^{-1})\in\tau^v(\DtE''\cap S)$) for all $n\geq N_S$.
Since $\tau^v(\DtE''\cap S)$ is commensurable with $\bold{H}(\C{O}_F)$, it is 
contained in a compact subset of $\bold{H}(\afv)$. Therefore 
the sequence $\{\tau^v(\dt_n\dt_{n+1}^{-1})\}_n\subset\bold{H}(\afv)$ has a limit 
point. Let $h$ be some limit point of $\{\tau^v(\dt_n\dt_{n+1}^{-1})\}_n$, and let
$\{\tau^v(\dt_{n_i}\dt_{n_i +1}^{-1})\}_i$ be a subsequence, converging to $h$. 
Then for each prime $u\neq v$ of $F$ we have 
$$h_u=\lim_{i\to\infty}\tau_u((\dt_{n_i}\dt_{n_i +1}^{-1})_u)
=\tau_u(\lim_{i\to\infty}(\dt_{n_i}\dt_{n_i +1}^{-1})_u)=1,$$
because $\tau_u$ is continuous. It follows that $1$ is the only limit point of 
$\{\tau^v(\dt_n\dt_{n+1}^{-1})\}_n$, therefore the sequence $\{\tau^v(\dt_n)
\tau^v(\dt_{n+1}^{-1})\}_n$ converges to $1$. Now by similar arguments we see that
the sequence $\{\tau^v(\dt_n)\}_n$ converges to  $\tau^v(g)\in\bold{H}(\afv)$. 

Moreover, the same arguments also imply that if we show that $\tau^v(PE')=\bold{H}(\afv)$,
then the continuity of $\tau^v$ and of $(\tau^v)^{-1}$ will follow automatically.

Observe that for each non-archimedean place $u$ we have $\bold{G}(F_u)
^{der}=\bold{G^{der}}(F_u)$ (resp. $\bold{G'}(F_u)^{der}=\bold{(G')^{der}}
(F_u)$) (see \cite[1.3.4 and Thm. 6.5]{PR} in the anisotropic and 
\cite[Thm. 7.1 and 7.5]{PR} in the isotropic cases respectively). Therefore $\tau^v$ 
induces an isomorphism of derived groups 
$\prod_{u\neq v} P(\bold{G^{der}}(F_u))\isom\prod_{u\neq v}P(\bold{(G')^{der}}(F_u))$. 

By \rp{closure}, $\overline{\DtE''}\supset P(\bold{G^{der}}(\afv))=
\PG(\afv)\cap\prod_{u\neq v} P(\bold{G^{der}}(F_u))$.
Hence by the facts shown above,
$$\tau^v(P(\bold{G^{der}}(\afv)))\subset\bold{PG'}(\afv)
\cap\prod_{u\neq v}P(\bold{(G')^{der}}(F_u))=
P(\bold{(G')^{der}}(\afv)).$$
In particular, $\tau^v(\prod_{u\neq v}\,P(\bold{G^{der}}(\C{O}_{F_u})))=
\prod_{u\neq v}\tau_u(P(\bold{G^{der}}(\C{O}_{F_u})))\subset 
P(\bold{(G')^{der}}(\afv))$. It follows that
$\tau_u(P(\bold{G^{der}}(\C{O}_{F_u}))\subset P(\bold{(G')^{der}}(\C{O}_{F_u}))$ for almost all $u\neq v$. 
Since each $\tau_u$ is algebraic, the subgroups $\tau_u(P(\bold{G^{der}}(\C{O}_{F_u}))$  and 
$P(\bold{(G')^{der}}(\C{O}_{F_u}))$ are conjugate (hence equal) for 
almost all $u\neq v$. It follows that $\tau^v(P(\bold{G^{der}}(\afv)))=P(\bold{(G')^{der}}(\afv))$.

Therefore to complete the proof it will suffice to show that $\PG(\afv)$ (resp. $\bold{PG'}(\afv)$)
is the normalizer of $P(\bold{G^{der}}(\afv))$ in the product ${\prod}_{u\neq v}\bold{PG}(F_u)$, and 
similarly for $\bold{PG'}$. Since $\bold{PG}(\afv)$ is the restricted
topological product of the $\bold{PG}(F_u)$'s with respect to the 
$\bold{PG}(\C{O}_{F_u})$'s, 
it remains to show that the normalizer of $P(\bold{G^{der}}(\C{O}_{F_u}))$ in 
$\bold{PG}(F_u)$ is $\bold{PG}(\C{O}_{F_u})$ for almost all $u$. This can be done by direct calculation.
\end{pf}

We will use the same letter $\tau$ to denote the isomorphism between $PE$ and 
$\bold{PG'}(\afv)$.

\begin{Emp} \label{E:fst3}
Notice that  regular function $t:=\Tr^d/\Det$ on $\GL{d}$ defines a function on
$\PGL{d}$. Moreover, an algebraic automorphism $\psi$ of $\PGL{d}$  is inner if and only if it satisfies
$t\circ\psi=t$. Therefore there is a unique choice of an algebra $D'$ defining $\bold{G'}$ (see \re{class})
such that the function $t':=\Tr^d/\Det$ on $\bold{PG'}$, defined by the natural embedding $\bold{G'}(F)\hra D'$, 
satisfies $t'\circ\varphi=t$.
\end{Emp}

\begin{Prop} \label{P:homom}
We have $D'\cong D\ii,\;\bold{G'}\cong \bold{G\ii}$, and $\tau$ is induced by some 
admissible isomorphism.
\end{Prop}

\begin{pf}
 By \rcr{ellel}, for each $\dt\in\Dt'$ with elliptic 
$\dt_{\infty}$ we have $t(\dt_{\infty})=t(\dt_{v})=t(\dt^{f;v})\in K$. Since
$(\varphi\times\tau)(\dt)\in \bold{PG'}(F)\subset\bold{PG'}(F_{\be_1}\times\af)$, 
we have $t((\varphi\times\tau)\dt)\in K$. By our assumption, 
$t(\varphi(\dt_{\infty}))=t(\dt_{\infty})$ for all $\dt\in\Dt'$. 
Hence for each $\dt\in\Dt'$ with elliptic $\dt_{\infty}$ we have 
$t(\tau_u(\dt_u))=t(\dt_{u})$ for each non-archimedean prime $u$ of $F$.

Recall that the algebraic isomorphism $\tau_u:\PG(F_u)\to \bold{PG'}(F_u)$ for $u\neq v$ 
is induced either by an $F_u$-linear isomorphism
$D\otimes_{F}F_u\isom D'\otimes_{F}F_u$ or by an $F_u$-linear isomorphism 
$D\otimes_{F}F_u\isom(D')\op\otimes_{F}F_u$, composed 
with an inverse map ($g\mapsto g^{-1}$). In the first case we have
$t(\tau_u(g_u))=t(g_{u})$ for all $g_u\in \bold{G_u}(F_u)$, and in the second one
$t(\tau_u(g_u))=t(g_{u}^{-1})$ for all $g_u\in \bold{G_u}(F_u)$. 

To exclude the second possibility we need to show the existence of a
$\dt\in\Dt'$ with elliptic $\dt_{\infty}$ such that $t(\dt_{\infty})\neq t(\dt_
{\infty}
^{-1})$. Since the condition $t(g)=t(g^{-1})$ is Zariski 
closed and non-trivial and since the closure of all elliptic elements of
$\Dt_{\infty}\in\pgr^0$ contains an open non-empty set, we are done. 

It follows 
that $D'$ is locally isomorphic to $D\ii$ at every non-archimedean place of $K$, 
except possibly at $w$ and $\bar{w}$, and that the map $\tau^v:\PG(\afv)\isom
\bold{PG'}(\afv)$ is induced by some admissible isomorphism. 
To prove the statement for the $v$-component we copy the above  proof replacing 
$\PG(F_u)$ by $\PGL{1}(\wt{D}_w)$ and $D\otimes_{F}F_u$ by $\wt{D}_w\oplus\wt{D}_w\op$.

Since $D'$ and $D\ii$ are locally isomorphic at all places, they are isomorphic.
We showed before that $\bold{PG'}(F_{\infty_1})\cong\pgr$, and that for each 
$i=2,...,g$ the group $\bold{PG'}(F_{\infty_i})$ is compact and, therefore, is
isomorphic to $\PGU{d}(\B{R})$. \rp{invol}, b) then implies that 
$\bold{G'}\cong\bold{G\ii}$. 
\end{pf}

From now on we identify $\bold{G'}$ with $\bold{G\ii}$.

\subsection{Completion of the proof} \label{compl}
 Our next task is to prove the following

\begin{Prop} \label{P:equality}
We have $(\varphi\times\tau)(\Dt')=\bold{PG\ii}(F)_{+}$.
\end{Prop}

\begin{pf}
 First observe that $(\varphi\times\tau)(\Dt')_{\infty}=\varphi(\Dt_{\infty})\
\subset\varphi(\pgr^0)=\bold{PG\ii}(F_{\infty_1})^0$, therefore $(\varphi\times\tau)
(\Dt')\subset \bold{PG\ii}(F)_{+}$ and $(\varphi\times\tau^v)(\Dt'')\subset
\bold{PG\ii}(F)_{+}$. Since the projection of $\bold{PG\ii}(F)$ to $\bold{PG\ii}(F_{\infty_1})\times
\bold{PG\ii}(\afv)$ is injective, it remains to show that
\begin{equation} \label{EE:index}
[\bold{PG\ii}(F):\bold{PG\ii}(F)_{+}]=[\bold{PG\ii}(F):(\varphi\times\tau^v)(\Dt'')].
\end{equation}

We are going to use of Kottwitz' results described in \ref{SS:inner}.
Recall that $\bold{PG\ii}$ is an inner form of $\PG$. Let $\omega_{\bold{PG}}$ and 
$\omega_{\bold{PG\ii}}$ be non-zero invariant differential forms
of top degree on  $\PG$ and $\bold{PG\ii}$ respectively, connected with one another by some inner
twist as in \re{Tam}. They define invariant measures $|\omega_{\bold{PG}}|$
and $|\omega_{\bold{PG\ii}}|$ on  $\PG(F_u)$ and
$\bold{PG\ii}(F_u)$ for every completion 
$F_u$ of $F$ and product measures on $\PG(\B{A}_F)$ and $\bold{PG\ii}(\B{A}_F)$  respectively
(see \cite[Ch. 2]{We2}). It follows from Weil's conjecture on Tamagawa
numbers and from Ono's result (see Ono's appendix to \cite{We2}) that
\begin{equation} \label{EE:first}
|\omega_{\bold{PG\ii}}|(\bold{PG\ii}(\B{A}_F)/\bold{PG\ii}(F))=|\omega_{\PG}|(\PG(\B{A}_F)/\PG(F)).
\end{equation}

\begin{Lem} \label{L:prod}
Let $A$ and $B$ be locally compact groups, let $S$ be a compact and open subgroup
of $A$, and let $\Gm$ be a lattice in $A\times B$ with injective projection to
$B$. Then for every right invariant measures $\mu_A$ on $A$ and $\mu_B$ on $B$
we have 
$$(\mu_A\times\mu_B)([A\times B]/\Gm)=\mu_A(S)\cdot\mu_B([(S\bs A)\times B]/\Gm).$$
\end{Lem}

\begin{pf}
Let $\Gm_A$ be the projection of $\Gm$ to $A$. Choose representatives 
$\{a_i\}_{i\in I}$ of the double classes $S\bs A/\Gm_{A}$. For each $i\in I$ let
$\Gm_i$ be the projection of the subgroup $(a_i^{-1}Sa_i\times B)\cap\Gm$ to $B$.
Then $\Gm_i$ is a lattice in $B$, therefore there exists a measurable subset  
$U_i$ of $B$ such that $B$ is the disjoint union 
$\underset{\gm\in\Gm_i}{\coprod}U_i\gm$. Since $\Gm$ has an injective projection to 
$B$, we have $A\times B=\underset{\gm\in\Gm}{\coprod}\underset{i\in I}{\coprod}
(Sa_i\times U_i)\gm$. Then
\begin{align}
(\mu_A\times\mu_B)([A\times B]/\Gm) & =\underset{i}\sum\,\mu_A(Sa_i)\cdot
\mu_B(U_i)=\mu_A(S)\cdot\underset{i}\sum\,\mu_B(U_i)\notag\\
   & =\mu_A(S)\cdot\underset{i}\sum\,\mu_B(a_i\times U_i)=\mu_A(S)\cdot
   \mu_B([(S\bs A)\times B]/\Gm)\notag.
\end{align}\end{pf}

 By the lemma, for each $S\in\ff{\PG(\afv)}$ the left hand side of (\ref{EE:first}) is 
equal to
\begin{equation} \label{EE:second}
\prod_{i=2}^g |\omega_{\bold{PG\ii}}|(\bold{PG\ii}(F_{\be_i}))\cdot|\omega_{\bold{PG\ii}}|(\bold{PG\ii}(F_v)) 
\cdot|\omega_{\bold{PG\ii}}|(\tau^v(S))\cdot\cr
|\omega_{\bold{PG\ii}}|(\tau^v(S)\bs[\bold{PG\ii}(F_{\infty_1})\times 
\bold{PG\ii}(\afv)]/\bold{PG\ii}(F))
\end{equation}
and the right  hand side of (\ref{EE:first}) is equal to
\begin{equation} \label{EE:third}
\prod_{i=1}^g |\omega_{\PG}|(\PG(F_{\be_i}))\cdot|\omega_{\PG}|(S) 
\cdot|\omega_{\PG}|(S\bs\PG(\af)/\PG(F)).
\end{equation}

By definition, $|\omega_{\bold{PG\ii}}|(\bold{PG\ii}(F_{\be_i}))=
|\omega_{\PG}|(\PG(F_{\be_i}))$ for each $i=2,...,d$ and 
$|\omega_{\bold{PG\ii}}|(\tau^v(S))=|\omega_{\PG}|(S)$ for each $S\in\ff{\PG(\afv)}$. 

Since the expressions of (\ref{EE:second}) and (\ref{EE:third}) are equal,
\rp{compat} and \rr{compat} imply that
\begin{equation} \label{EE:fourth}
\mu_{\bold{PG\ii_{F_{\infty_1}}}}(\tau^v(S)\bs[\bold{PG\ii}(F_{\infty_1})\times 
\bold{PG\ii}(\afv)]/\bold{PG\ii}(F)_+)=\cr
=d\cdot\mu_{\bold{PG_{F_v}}}(S\bs\PG(\af)/\PG(F))
\end{equation}
($'+'$ was added to multiply the left hand side by $2$ when  $d=2$).

If $S$ is sufficiently small, then  for each
$a\in\PG(\afv)$ the group $a^{-1}Sa\cap \PG(F)$ is torsion-free by \rp{basic}.
Let $Y_{a^{-1}Sa}$ be the projective variety over $K_w$ such that $Y_{a^{-1}Sa}\an
\cong({a^{-1}Sa}\cap\PG(F))_v\bs\Om$.
By Kurihara's result (see \cite[Thm. 2.2.8]{Ku}) $c_{d-1}(T_{Y_{a^{-1}Sa}})=
\chi_E({a^{-1}Sa}\cap \PG(F))\cdot c_{d-1}(T_{\B{P}^{d-1}})$,
where $c_{d-1}(T_{Y_{a^{-1}Sa}})$ (resp. $c_{d-1}(T_{\B{P}^{d-1}}))$ is the
$(d-1)$-st Chern 
class of the tangent bundle of $Y_{a^{-1}Sa}$ (resp. $\B{P}^{d-1}$). Notice that 
$c_{d-1}(T_{\B{P}^{d-1}})=d$, hence  $c_{d-1}(T_{Y_{a^{-1}Sa}})=
d\cdot\mu_{\bold{PG_{F_v}}}(({a^{-1}Sa}\cap\PG(F))_v\bs\PG(F_v))$.

Since $(Y_{{a^{-1}Sa},\B{C}})^{an}\cong\Dt''_{a^{-1}Sa}\bs B^{d-1}$, we have
$c_{d-1}(T_{Y_{a^{-1}Sa}})=c_{d-1}(T_{(Y_{a^{-1}Sa})_{\B{C}}})=\chi_E(\Dt''_
{a^{-1}Sa}\bs B^{d-1})$ (see
for example \cite[Prop. 11.24 and (20.10.6)]{BT}). The last expression is equal to
$\chi_E(\Dt''_{a^{-1}Sa})=\mu_{\PGU{d-1,1}}(\Dt_{a^{-1}Sa}''\bs\pgr)$. This shows 
that for each  $a\in\PG(\afv)$ we have
$$d\cdot\mu_{\bold{PG_{F_v}}}(({a^{-1}Sa}\cap\PG(F))_v\bs\PG(F_v))=\mu_{\PGU{d-1,1}}
(\Dt_{a^{-1}Sa}''\bs\pgr).$$

Summing this equality for $a$ running over a set of representatives of double
classes \qquad\qquad\linebreak $S\bs\PG(\afv)/\PG(F)$ we obtain that
$$d\cdot\mu_{\bold{PG_{F_v}}}(S\bs\PG(\af)/\PG(F))=
\mu_{\PGU{d-1,1}}(S\bs[\pgr\times \PG(\afv)]/\Dt'').$$

Since the right hand side of the last expression is equal to 
$$\mu_{\bold{PG\ii_{F_{\infty_1}}}}(\tau^v(S)\bs[\bold{PG\ii}(F_{\infty_1})\times 
\bold{PG\ii}(\afv)]/(\varphi\times\tau^v)(\Dt'')),$$
we conclude (\ref{EE:index}) from (\ref{EE:fourth}).
\end{pf}
 
\begin{Emp} \label{E:complet}
 By \rp{homom} there exists an admissible 
isomorphism $\Phi:E\isom E\ii$, inducing the isomorphism 
$\tau:PE\isom PE\ii$.
Choose $\dt\in\Dt$ with elliptic  $\dt_{\be}\in\Dt_{\be}$ and $\Tr\Ad(\dt_{\infty})\neq 
-1$. Choose its 
representative $\wt{\dt}\in\gr^0\times E$ as in \rcr{ellel}. Then 
$(\Tr\,\wt{\dt})(\Tr\,\wt{\dt}^{-1})\in K\m$. Let $\dt'$ be the projection of 
$\wt{\dt}$ to $\pgr^0\times E$. Set $\wt{\gm}:=(\varphi\times\Phi)(\dt')$, and let 
$\wt{\gm}_E$ be its projection to $E$.

By the definition of admissible map, $\Tr(\wt{\gm}_E)\in 
K\m$. Let $\bar{\dt}$ be the image of $\dt$ in $\Dt'$, then $\bar{\gm}:=
(\tau\times \varphi)(\bar{\dt})$ belongs to $\bold{PG\ii}(F)_{+}$. Let $\gm'\in 
\bold{G\ii}(F)_{+}$ be 
some representative of $\bar{\gm}$, then $\wt{\gm}_{E}^{-1}\gmE'\in Z(E\ii)$.
Therefore $\wt{\gm}_{E}^{-1}\gmE'=(\Tr\,\wt{\gm}_{E})^{-1}(\Tr\,\gmE')\in K\m=
Z(\bold{G\ii}(F))$. Thus $\wt{\gm}_{E}$ and $\gmE'$ have equal projections to 
$\pgr^0\times(E\ii/E_0\ii)$, hence $(\varphi\times\bar{\Phi})(\dt)\in\Gm\ii$. 

The condition $\{\dt_{\infty}$ is elliptic and $\Tr\Ad(\dt_{\infty})\neq -1\}$
is open and non-empty, therefore the above $\dt$'s generate the 
whole group $\Dt\cong\Dt_{\infty}$ (see \cite[Ch. IX, Lem. 3.3]{Ma}). 
It follows that $(\varphi\times\bar{\Phi})(\Dt)\subset\Gm\ii$. 
Since the projection $\pi:\Gm\ii\to \bold{PG\ii}(F)_{+}$ is an isomorphism,
\rp{equality} implies that $(\varphi\times\bar{\Phi})(\Dt)=
\Gm\ii$. This completes the proof of \rt{complex} and of the First Main Theorem. 
\end{Emp}

\section{Theorem on the $p$-adic uniformization}

 The First Main Theorem implies that for some admissible isomorphism $\Phi:
E\isom E\ii$ there exists a $\Phi$-equivariant $\B{C}$\/-rational isomorphism 
$f_{\Phi}:X_{\B{C}}\isom \wt{X}\ii$. Therefore for some $\B{C}/K_w$-descent
$X\ii$ of the \esc{E\ii}{\B{C}} $\wt{X}\ii$, $f_{\Phi}$ induces a 
$K_w$-rational isomorphism
$X\isom X\ii$. To describe $X\ii$ we need some preparations, following 
\cite{De} (see also \cite{Mi})

\subsection{Technical preliminaries} \label{SS:tecprel}

In this subsection we recall basic notions related to Shimura varieties
and give their explicit description in the cases we are interested in.

\begin{Emp} \label{E:Shim}
First we realize $\wt{X}\ii$ as a Shimura variety. Set 
$\bold{H\ii}:=\bold{R}_{F/\B{Q}}\bold{G\ii}$. Then $\bold{H\ii}$ is a reductive group 
over $\B{Q}$ such that $\bold{H\ii}(\B{A}^f)=\bold{G\ii}(\af)$ and 
$\bold{H\ii}(\B{R})=\prod_{i=1}^r\bold{G\ii}(F_{\be_i})$.
Put $\bold{S}:=\bold{R}_{\B{C}/\B{R}}\B{G}_m$, and let $h$ be a homomorphism $\bold{S}\to \bold{H\ii}_{\B{R}}$
such that for each $z\in\B{C}\m\cong\bold{S}(\B{R})$ we have
$$ h(z)=(diag(1,...,1,z/\bar{z})^{-1};I_d;...;I_d)\in\prod_{i=1}^r\bold{G\ii}(F_{\be_i}),$$ 
using the identification of $\bold{G\ii}(F_{\be_1})$ with $\gr$, chosen in \re{ex2}.
Then the conjugacy class $M\ii$ of $h$ in $\bold{H\ii}(\B{R})$ is isomorphic to $B^{d-1}$ if $d>2$
and to $\B{C}\sm\B{R}$ if $d=2$. 
The pair $(\bold{H\ii},M\ii)$ satisfies Deligne's axioms (see \cite[1.5 and 2.1]{De} or 
\cite[II, 2.1]{Mi}), and 
the Shimura variety  $M_{\B{C}}(\bold{H\ii},M\ii)$, corresponding to $(\bold{H\ii},M\ii)$,
is isomorphic to $\wt{X}\ii$. 
\end{Emp}

\begin{Emp} \label{E:split}
For each pair $(\bold{H\ii},M\ii)$ as above there is a number field 
$E(\bold{H\ii},M\ii)\subset\B{C}$, called the reflex field of $(\bold{H\ii},M\ii)$, 
which is defined as follows (compare \cite[1.2, 1.3 and 3.7]{De}).
The group $\Hom_{\B{C}}(\bold{S}_{\B{C}},(\B{G}_m)_{\B{C}})$ is a 
free abelian
group of rank 2 with generators $z$ and $\bar{z}$ such that if $i:\bold{S}(\B{R})
\hra \bold{S}(\B{C})$ is the natural inclusion, then for each $w\in\B{C}\m\cong\bold{S}
(\B{C})$ we have $z\circ i(w)=z$ and $\bar{z}\circ i(w)=\bar{w}$. Let 
$r:(\B{G}_m)_{\B{C}}\to\bold{S}_{\B{C}}$ be the algebraic homomorphism such that 
$(z^p \bar{z}^q)\circ r(x)=x^p$. Then $E(\bold{H\ii},M\ii)$ is the field of definition 
of the conjugacy class of the composition map $r'':(\B{G}_m)_{\B{C}}\overset{r}{\lra}\bold{S}_{\B{C}}
\overset{h}{\lra}\bold{H\ii}_{\B{C}}$.
\end{Emp}

\begin{Prop} \label{P:refl1}
We have $E(\bold{H\ii},M\ii)=K$, where the latter is viewed as a subfield of $\B{C}$ through the embedding
$\be_1$, chosen in  \re{ex2}.
\end{Prop}

\begin{pf}
Note that $\bold{H\ii}(\B{C})$ is naturally embedded into $\GL{d}(\B{C})^{2g}$ so that 
each factor corresponds to an embedding of $K$ into $\B{C}$. Supposing that the first and the second
factors corresponds to our fixed embedding and to its complex conjugate respectively we have
$$ r''(z)=(diag(1,...,1,z^{-1});diag(1,...,1,z);I_d;...;I_d)\text{ for each }z\in\B{C}\m.$$
Therefore the reflex field $E(\bold{H\ii},M\ii)$ contains $K\subset\B{C}$. On the other hand, 
the Skolem-Noether theorem implies that for each $\sigma\in\Aut_K(\B{C})$ the homomorphism
${\sigma}(r'')$ is conjugate to $r''$. This implies the assertion.
\end{pf}

\begin{Emp} \label{E:maxtor}
Let $\bold{T}\subset\bold{H\ii}$ be a maximal torus of $\bold{H\ii}$, 
defined over $\B{Q}$,
such that some conjugate $h'\in M\ii$ of $h$ in $\bold{H\ii}(\B{R})$ factors through 
$\bold{T}_{\B{R}}$. Then we have a natural 
embedding $i_{\bold{T}}:M_{\B{C}}(\bold{T},h')\hra M_{\B{C}}(\bold{H\ii},M\ii)$, 
where $M_{\B{C}}(\bold{T},h')$ is the Shimura variety corresponding to 
$(\bold{T},h')$. Since $\bold{T}$ is commutative, the reflex field 
$E_{\bold{T}}:=E(\bold{T},h')$ of $(\bold{T},h')$ is the field of 
definition of the morphism $r'':(\B{G}_m)_{\B{C}}\overset{r}{\lra}\bold{S}_{\B{C}}
\overset{h'}{\lra}T_{\B{C}}$. Hence $r''$ defines a morphism of algebraic groups
over $\B{Q}$ 
$$r':E_{\bold{T}}^{*}:=\bold{R}_{E_{\bold{T}}/\B{Q}}(\B{G}_m)
\overset{\bold{R}_{E_{\bold{T}}/\B{Q}}(r'')}{\lra}\bold{R}_{E_{\bold{T}}/\B{Q}}
(\bold{T})\overset{\text{N}_{E_{\bold{T}}/\B{Q}}}{\lra}\bold{T}.$$ 
Notice that $E_{\bold{T}}\supset E(\bold{H\ii},M\ii)$. 
Let $\theta_{E_{\bold{T}}}$ be the Artin isomorphism of global class field 
theory sending the uniformizer to the arithmetic Frobenius automorphism. 
Let $\la_{\bold{T}}:
\Gal(E_{\bold{T}}^{ab}/E_{\bold{T}})\to \bold{T}(\B{A}^f)/\overline
{\bold{T}(\B{Q})}$ be the composition map $$\Gal(E_{\bold{T}}^{ab}/
E_{\bold{T}})\overset{\theta_{E_{\bold{T}}}^{-1}}
{\overset{\sim}{\lra}}E_{\bold{T}}^{*}(\B{R})^{\circ}\bs E_{\bold{T}}^{*}(\B{A})
/\overline{E_{\bold{T}}^{*}(\B{Q})}\overset{r'}{\lra}\bold{T}(\B{R})^{\circ}\bs \bold{T}(\B{A})/\overline
{\bold{T}(\B{Q})}\overset{\text{proj}}{\lra}\bold{T}(\B{A}^f)/\overline{\bold{T}
(\B{Q})}.$$
For each $E'\supset E(\bold{H\ii},M\ii)$ we denote the composition map 
$\Gal(E_{\bold{T}}^{ab}\cdot E'/E')\overset{\text{res}}{\lra}\Gal(E_{\bold{T}}
^{ab}/E_{\bold{T}})\overset{\la_{\bold{T}}}{\lra}\bold{T}(\B{A}^f)/
\overline{\bold{T}(\B{Q})}$ by $\la_{\bold{T},E'}$.
\end{Emp}

\begin{Lem} \label{L:maxtor}
 Each maximal torus $\bold{T}$ of $\bold{H\ii}$, defined over $\B{Q}$, is equal to
the  intersection of $\bold{H\ii}$ with $\bold{R}_{L/\B{Q}}\B{G}_m$ for a unique maximal 
commutative subfield $L$ of $D\ii$. (In such
a situation we will call $\bold{T}$ an $L$-torus). In this case, $\al\ii$ induces 
a nontrivial automorphism of $L$, and the subgroup 
$\bold{T}(\B{Q})\subset L\m\cong\bold{R}_{L/K}\B{G}_m(K)$ is Zariski dense in 
$\bold{R}_{L/K}\B{G}_m$.
\end{Lem}

\begin{pf}
Let $L$ be the subalgebra of $D\ii$ spanned over $K$ by $\bold{T}(\B{Q})\subset\bold{H\ii}(\B{Q})\subset D\ii$, then $L$ is a commutative subfield and 
$\bold{T}(\B{Q})\subset \bold{H\ii}(\B{Q})\cap \bold{R}_{L/\B{Q}}\B{G}_m(\B{Q})$.
Since $\bold{T}$ is connected and $\B{Q}$ is infinite and 
perfect, the subgroup $\bold{T}(\B{Q})$ is Zariski dense in
 $\bold{T}$ (see \cite[Ch. V, Cor. 13.3]{Bo}). It follows that 
$\bold{T}\subset\bold{H\ii}\cap\bold{R}_{L/\B{Q}}\B{G}_m$. 
Since $\bold{T}$ is maximal, $L$ have to be maximal and $\bold{T}=\bold{H\ii}\cap\bold{R}_{L/\B{Q}}\B{G}_m$.

For each $g\in\bold{T}(\B{Q})$ we have $\al\ii(g)\in g^{-1}F\m\subset\bold{T}(\B{Q})$, 
so that $\al\ii(L)=L$. To prove the 
last assertion we observe that there exists a maximal $F$-rational subtorus
$\bold{T'}$ of $\bold{G\ii}$ such that $\bold{T}=\bold{R}_{F/\B{Q}}(\bold{T'})$. 
Then the subgroup $\bold{T}(\B{Q})=\bold{T'}(F)$ is Zariski dense in
$\bold{T'_K}\cong\bold{R}_{L/K}\B{G}_m\times(\B{G}_m)_K$. Hence its projection to $\bold{R}_{L/K}\B{G}_m$
is also Zariski dense.
\end{pf}

\begin{Emp} \label{E:refl}
Now we want to calculate the reflex field $E_{\bold{T}}$. Observe that 
$L\otimes_{\B{Q}}\B{C}\subset D\ii\otimes_{\B{Q}}\B{C}
\cong \Mat_d(\B{C})^{2g}$. Possibly after a conjugation we may assume that
$L\otimes_{\B{Q}}\B{C}$ is the subalgebra of diagonal matrices of 
$\Mat_d(\B{C})^{2g}$. Then each diagonal entry of each of the $2g$ copies of
$\Mat_d(\B{C})$ corresponds to an embedding of 
$L$ into $\B{C}$, and the map $r'':(\B{G}_m)_{\B{C}}\to\bold{T}_{\B{C}}$ is as 
follows:
$$r''(z)=(diag(1,...,1,z^{-1});diag(1,...,1,z);I_d;...;I_d).$$
Let $\iota_1$ be the embedding $L\hra\B{C}$, corresponding to the right low entry of the 
first matrix, then the right low entry of the second matrix corresponds to the 
embedding $\bar{\iota}_1:=\iota_1\circ\al\ii$. Now we embed $L$ into $\B{C}$ via $\iota_1$.
\end{Emp}

\begin{Prop} \label{P:refl2}
We have $E_{\bold{T}}=L\subset\B{C}$, and $r':E_{\bold{T}}^{*}\to\bold{T}$
is characterized by $r'(l)=l^{-1}\cdot\al\ii(l)$ for each
$l\in E_{\bold{T}}^{*}(\B{Q})\subset L\m$.
\end{Prop}

\begin{pf}
As was noted before, $E_{\bold{T}}\supset E(\bold{H\ii},M\ii)$. Hence by \rp{refl1}, 
$E_{\bold{T}}\supset K$. By the definition,
$\sigma(r''(z)) =r''(\sigma(z))$ for each $\sigma\in \Aut_{E_{\bold{T}}}(\B{C})$, 
hence the group $\Aut_{E_{\bold{T}}}(\B{C})$ must stabilize $\iota_1$, so that
$E_{\bold{T}}\supset L$. Finally, it is clear that $r''$ is defined over 
$L\subset\B{C}$. For each $l\in E_{\bold{T}}=L$ we have 
$$r'(l)=\text{N}_{L/\B{Q}}(diag(1,...,1,l^{-1});
diag(1,...,1,l);I_d;...;I_d)=l^{-1}\cdot\al\ii(l).$$
\end{pf}

Set $L_w:=L\otimes_{K}K_w\subset D_w\ii:=D\ii\otimes_{K}K_w$. Since
$D_w\ii$ is a division algebra, $L_w$ is a field extension of $K_w$ of 
degree $d$, and $L_w=L\cdot K_w$.

\begin{Lem} \label{L:field}
The following relations hold:

a) $E_{\bold{T}}\cdot K_w=L_w$;\quad b) $E_{\bold{T}}^{ab}\cdot K_w=(L_w)^{ab}$.
\end{Lem}

\begin{pf}
a) was proved above.\par

b) The group $\Gal(E_{\bold{T}}^{ab}\cdot K_w/E_{\bold{T}}\cdot K_w)$ is abelian, 
hence $L_w\subset E_{\bold{T}}^{ab}\cdot K_w\subset L_w^{ab}$. 
By the class field theory, the composition of the canonical projections
$\Gal((L_w)^{ab}/L_w)\to\Gal(E_{\bold{T}}^{ab}\cdot K_w/E_{\bold{T}}\cdot K_w)
\to\Gal(E_{\bold{T}}^{ab}/E_{\bold{T}})\cong\Gal(L^{ab}/L)$ is injective 
(use, for example, \cite[Ch. VII, Prop. 6.2]{CF}). Therefore we have the required 
equality.
\end{pf}

\begin{Prop} \label{P:reclaw}
For each $l\in L_w\m\subset (D\ii_w)\m$ the element
$$(l^{-1},1,1)\in(D_w\ii)\m\times F_v\m\times\bold{G\ii}(\afv)\cong\bold{H\ii}(\B{A}^f)$$
belongs to $\bold{T}(\B{A}^f)$, and its equivalence class in $\bold{T}(\B{A}^f)/\overline{\bold{T}(\B{Q})}$
is $\la_{\bold{T},K_w}(\theta_{L_w}(l))$.
\end{Prop}

\begin{pf}
The statement follows immediately from the explicit formulas of \rp{refl2} using the connection between local
and global Artin maps.
\end{pf}

\begin{Def} \label{D:sppoint}
A point $x\in M_{\B{C}}(\bold{H\ii},M\ii)(\B{C})$ is called {\em ($\bold{T}$-)special} if 
$x\in i_{\bold{T}}(M_{\B{C}}(\bold{T},h)(\B{C}))$.
\end{Def}

\begin{Rem} \label{R:spec}
The group $\bold{T}(\B{A}^f)$ acts naturally on the set of 
$\bold{T}$-special points and
the group $\bold{T}(\B{Q})$ acts on it trivially. Hence by 
continuity the closure 
$\overline{\bold{T}(\B{Q})}\subset\bold{T}(\B{A}^f)$ acts trivially on the set of 
$\bold{T}$-special points, therefore the action of
$\bold{T}(\B{A}^f)/\overline{\bold{T}(\B{Q})}$ on it is well-defined.
\end{Rem}

\begin{Def} \label{D:canmod}
Let  $K'\supset E(\bold{H\ii},M\ii)$ be a subfield of $\B{C}$. A $\B{C}/K'$-descent of 
the \esc{\bold{H\ii}(\B{A}^f)}{\B{C}} $M_{\B{C}}(\bold{H\ii},M\ii)$
is called {\em weakly-canonical} if for each maximal torus $\bold{T}\hra\bold{H\ii}$ as above,
each $\bold{T}$-special point $x$ is defined over
$E_{\bold{T}}^{ab}\cdot K'$, and for each 
$\sigma\in \Gal(E_{\bold{T}}^{ab}\cdot K'/E_{\bold{T}}\cdot K')$ we 
have $\sigma(x)=\la_{\bold{T},K'}(\sigma)(x)$.
\end{Def}

\begin{Rem} \label{R:canmod}
Our definition of the canonical model coincides with that of \cite{Mi3}, which differs from those of \cite{De2} and 
\cite{Mi} (see the discussion in \cite[1.10]{Mi3}). The seeming difference (by sign) between our reciprocity map and 
that of \cite{Mi3} is due to the fact that we consider left action of the adelic group whereas Milne considers right
action.  
\end{Rem}

\begin{Prop} \label{P:canmod}
 For each field $K'$ satisfying $E(\bold{H\ii},M\ii)\subset K'\subset\B{C}$ there exists a 
unique (up to an isomorphism) weakly-canonical $\B{C}/K'$-descent of the \hfill\linebreak
\esc{\bold{H\ii}(\B{A}^f)}{\B{C}} $M_{\B{C}}(\bold{H\ii},M\ii)$.
\end{Prop}

\begin{pf}
 Uniqueness is proved in \cite[5.4]{De}, for the existence see \cite[6.4]{De}
or \cite[II, Thm. 5.5]{Mi}.
\end{pf}

By \rp{refl1}, we have $E(\bold{H\ii},M\ii)=K\subset K_w\subset\B{C}$ 
(in our convention \re{fix}). Hence by \rp{canmod}, the \esc{E\ii}{\B{C}} 
$\wt{X}\ii$ has a unique weakly-canonical $\B{C}/K_w$-descent $X\ii$.

\subsection{Theorem on the $p$-adic uniformization} \label{SS:padun}

Now we are ready to formulate our 

\begin{Mains} \label{M:second}
For each admissible isomorphism $\Phi:E\isom E\ii$ there exists a 
$\Phi$-equivariant 
isomorphism $f_{\Phi}$ from the \esc{E}{K_w} $X$ to the \esc{E\ii}{K_w} $X\ii$.
\end{Mains}

\begin{Cor} \label{C:Main}
After the identification of $E$ with $E\ii$ by means of $\Phi$, we have for each 
$\SIE$ of the form $T_n\times S'$ with $S'\in\ff{E'}$ an isomorphism of 
$K_w$-analytic spaces $\varphi_S:(X_S\ii)\an\isom\gp\bs(\SiN\times(S'\bs 
\bold{G}(\af)/\bold{G}(F)))$. These isomorphisms commute with the natural projections for 
$T\supset 
S$ and with the action of $E\overset {\Phi}{\overset{\sim}{\lra}}E\ii$. 
\end{Cor}

\begin{pf}
(of the Second Main Theorem)\par
{\bf Step 1.} We want to prove that for $\Phi$ and $f_{\Phi}$ as in the First 
Main Theorem, the $\B{C}/K_w$-descent of $\wt{X}\ii$ corresponding to 
$X$ is weakly-canonical. For this we have to show that for each maximal torus 
$\bold{T}\hra\bold{H\ii}$ as in
\re{maxtor} and each $\bold{T}$-special point $x=f_{\Phi}(y)\in 
M_{\B{C}}(\bold{H\ii},M\ii)(\B{C})=\wt{X}\ii$  we have:\par

a) $y\in X(\B{C_p})$ is defined over $E_{\bold{T}}^{ab}\cdot K_w$;\par
b) $\sigma(y)=\Phi^{-1}(\la_{\bold{T},M}(\sigma))(y)$ for each 
$\sigma\in\Gal(E_{\bold{T}}^{ab}\cdot K_w/E_{\bold{T}}\cdot K_w)$.

By \rp{reclaw}, \rl{field} and the definition of admissible map, it will suffice 
to show that when $L_w$ is embedded into $\wt{D}_w$ by means of 
the isomorphism $D_w\ii\isom\wt{D}_w$ from Definition \ref{D:adm} we have
\begin{equation} \label{E:cond}
\text{i) every point } y\in X(\B{C_p}),\text{ fixed by } \Phi^{-1}(\bold{T}(\B{Q})),
\text{ is rational over }(L_w)^{ab};\cr
\text{ii) }\theta_{L_w}(l)(y)=l^{-1}(y)\text{ for each }l\in L_w\m\subset\wt{D}_w\m\cong\wt{D}_w\m\times\{1\}\subset\wt{D}_w\m\times E'.
\end{equation}

Let $(x,a)\in\Si\times E'$ be a representative of 
$y\in X(\B{C}_p)$. Then $(\sigma(x),a)$ is a representative of 
$\sigma(y)$ for each (not necessarily continuous) $\sigma\in\Aut_{K_w}(\B{C}_p)$. Recall that for each
embedding $L_w\hra\Mat_d(K_w)$ there exists an $(L_w\m\times L_w\m)$-equivariant $L_w$-rational
embedding $\wt{\imath}:\si{1}{L_w}\hra \Si$. 

\begin{Prop} \label{P:Drin}
There exists an embedding $L_w\hra \Mat_d(K_w)$ such that the image of the 
corresponding $\wt{\imath}:\si{1}{L_w}\hra \Si$ contains $x$.
\end{Prop}

\begin{pf}
Let $x'\in\Om$ be the projection of $x$. Then $y':=[(x',a)]\in(\wt{D}_w\m\bs X)
\an\cong(\Om\times(E')^{disc})/\bold{G}(F)\overline{Z(\bold{G}(F))}$ (use
\rp{prop2}) is the projection of $y$. Since 
$\Phi^{-1}(\bold{T}(\B{Q}))$ stabilizes $y$, it also stabilizes $y'$, therefore the projection of
$a^{-1}\Phi^{-1}(\bold{T}(\B{Q}))a$ to $E'$ is contained in $\bold{G}(F)\overline{Z(\bold{G}(F))}\subset E'$.
In other words for each $t\in \bold{T}(\B{Q})\subset D\ii$  we have $\pr_{E'}(a^{-1}\Phi^{-1}(t)a)=g\cdot z$ 
for some $g\in \bold{G}(F)$ and some $z\in\overline{Z(\bold{G}(F))}$.

 Since $\Phi$ is induced by some algebra isomorphism $D(\afv)\isom D\ii(\afv)$, we have
 $\Tr\,t=\Tr(a^{-1}\Phi^{-1}(t)a)=(\Tr\,g)\cdot z$. Therefore for $t$'s with non-zero trace
we get $z=(\Tr\,t)\cdot (\Tr\,g)^{-1}\in K\m\subset(\B{A}_K^{f;w,\bar{w}})\m$. This means that 
$\pr_{E'}(a^{-1}\Phi^{-1}(t)a)\in \bold{G}(F)\subset E'$.
As the set of all $t$'s in $\bold{T}(\B{Q})\subset D\ii$ with $\Tr\,t\neq 0$ 
generates $L\subset D\ii$ as an algebra, the map $l\mapsto\pr_{E'}(a^{-1}\Phi^{-1}(l)a)$
defines embeddings $L\hra D$ and $L_w\hra D\otimes_K K_w\cong \Mat_d(K_w)$. 

This shows that $\pr_{E'}(a^{-1}\Phi^{-1}( \bold{T}(\B{Q}))a)\subset\bold{G}(F)\subset E'$, so that
$a^{-1}\Phi^{-1}(\bold{T}(\B{Q}))a\subset \wt{D}_w\m\times\GmE\subset \wt{D}_w\m\times E'=E$.
Hence $\Phi^{-1}(\bold{T}(\B{Q}))$ preserves $\rho_a(\Si)\subset 
(X_{\B{C}_p})\an$ (in the notation of \rco{etcom}). Moreover, it follows from 
the definition of the embeddings $L_w\hra\wt{D}_w$ and $L_w\hra\Mat_d(K_w)$ that
for each $t\in \bold{T}(\B{Q})\subset L\m\subset L_w\m$ the image
of $a^{-1}\Phi^{-1}(t)a\subset\wt{D}_w\m\times\GmE$ under the
canonical map $\wt{D}_w\m\times\GmE\isom\wt{D}_w\m\times\GmG\subset\wt{D}_w\m
\times\gp$ is equal to $(t,t)$. Since
$y$ is a fixed point $\Phi^{-1}(\bold{T}(\B{Q}))$, we conclude from the above that 
$(t,t)(x)=x$ for every $t\in \bold{T}(\B{Q})$.
Noticing that $\bold{T}(\B{Q})$ is Zariski dense in $\bold{R}_{L/K}\B{G}_m$ by \rl{maxtor}
and that $\bold{R}_{L/K}\B{G}_m\otimes_K K_w\cong\bold{R}_{L_w/K_w}\B{G}_m$,
\rl{image} completes the proof.
\end{pf}

 Since $\wt{\imath}$ is $(L_w\m\times L_w\m)$-equivariant and $L_w$-rational,
the proposition together with \rl{Galact} imply (\re{cond}).
In other words, we have proved that for some admissible isomorphism $\Phi:E\isom E\ii$ 
there exists a $\Phi$-equivariant $K_w$-linear isomorphism $f_{\Phi}:X\isom X\ii$.

{\bf Step 2.} Let $\Psi$ be another admissible isomorphism $E\isom E\ii$.
The definition of admissibility together with the theorem of Skolem-Noether imply 
that $\Psi\circ\Phi^{-1}:E\ii\isom E\ii$ is an inner automorphism, so that there 
exists $g_{\Psi}\in E\ii$ such that $\Psi\circ\Phi^{-1}(g)=g_{\Psi}g g_{\Psi}
^{-1}$ for all $g\in E\ii$. Take $f_{\Psi}:X\overset{f_{\Phi}}{\lra}X\ii\overset
{g_{\Psi}}{\lra}X\ii$. Then for each $g\in E$ we have $f_{\Psi}\circ g=g_{\Psi}\circ
f_{\Phi}\circ g=g_{\Psi}\circ\Phi(g)\circ f_{\Phi}=(g_{\Psi}\circ\Phi(g)\circ 
g_{\Psi}^{-1})\circ (g_{\Psi}\circ f_{\Phi})=(\Psi\circ\Phi^{-1})(\Phi(g))\circ
f_{\Psi}=\Psi(g)\circ f_{\Psi}$, that is $f_{\Psi}$ is a $\Psi$-equivariant
isomorphism. This completes the proof of the Second Main Theorem.
\end{pf}

\section{ $P$-adic uniformization of automorphic vector bundles}

In the previous section we proved that the Shimura varieties corresponding to
the pairs $(\bold{H\ii},M\ii)$ have $p$-adic uniformization. Our next task is to show the analogous result
for automorphic vector bundles.

\subsection{Equivariant vector bundles} \label{SS:eqvb}


\begin{Emp} \label{E:int}
 Set $\bold{H}:=\bold{R}_{F/\B{Q}}\bold{G}$. Then  for some algebraic group $\wt{\bold{H}}$ over $K_w$
we have natural isomorphisms
$\bold{H_{K_w}}\cong\GL{d}\times\bold{\wt{H}}$, 
$\bold{PH_{K_w}}\cong\PGL{d}\times\bold{P\wt{H}}$ and $\bold{PH\ii_{K_w}}\cong 
\PGL{1}(D_w\ii)\times\bold{P\wt{H}}$, where the first factors 
correspond to the natural embedding $F\hra K_w$. Using these decompositions
let $\bold{PH_{K_w}}$ acts on $\B{P}^{d-1}_{K_w}$ through the natural 
action of the first factor and the trivial action of the second one, and let $\bold{H}(K_w)$, $\PH(K_w)$ and 
$\bold{PH\ii}(\B{R})^0\cong\pgr^0\times \PGU{d}(\B{R})^{g-1}$ act similarly on $\Si$, on $\Om$ and 
on $B^{d-1}$ respectively. Let $\beta_{\B{R}}$ be the natural 
embedding $B^{d-1}\hra (\B{P}_{\B{C}}^{d-1})\an$ and let $\beta_w$ (resp. $\beta_{w,n}$) 
be the composition 
of  the natural projection $\Si\to\Om$ (resp. $\SiN\to\Om$) and the natural
embedding $\Om\hra(\B{P}_{K_w}^{d-1})\an$.

Let $\pi\in K_w$ be a uniformizer, let $\wt{\Pi}$ be an element of $\GL{d}(K_w)$ satisfying 
$\wt{\Pi}^d=\pi$, and  let $\Pi'\in\pgp$
be the projection of $\wt{\Pi}$. Set $\Pi:=(\Pi',1)\in\PGL{d}(K_w)\times\bold{P\wt{H}}(K_w)\cong\PH(K_w)$.
Let $K_w^{(d)}$ be the unramified field extension of $K_w$ of degree $d$.
Since the Brauer invariant of $D_w\ii$ is $1/d$,  the group $\bold{PH\ii_{K_w}}$ is 
isomorphic to the quotient of $\bold{PH_{K_w}}\otimes_{K_w}K_w^{(d)}$ by the equivalence
relation $\Fr(x)\sim\Pi^{-1}x\Pi$, where $\Fr\in\Gal(K_w^{(d)}/K_w)$ is the 
Frobenius automorphism. For each scheme $Y$ over $K_w$ on which $\bold{PH_{K_w}}$ acts 
$K_w$-rationally define a twist $Y\tw :=(\Fr(x)\sim\Pi^{-1}x)\bs Y\otimes_{K_w}K_w^{(d)}$. Then 
$Y\otimes_{K_w}K_w^{(d)}\cong Y\tw\otimes_{K_w}K_w^{(d)}$, and the natural action of 
$\bold{PH\ii_{K_w}}$ on it is $K_w$-rational.

Let $W$ be a $\bold{PH_{K_w}}$-equivariant vector bundle on 
$\B{P}_{K_w}^{d-1}$, that is a vector bundle on $\B{P}_{K_w}^{d-1}$ equipped with an action of 
the group $\bold{PH_{K_w}}$, lifting its action on $\B{P}_{K_w}^{d-1}$.
Then $(W\tw,p\tw)$ is a $\bold{PH\ii_{K_w}}$-equivariant vector bundle on $(\B{P}_{K_w}^{d-1})\tw$,
 and $\beta_{\B{R}}^{*}((W\tw_{\B{C}})\an)$ (resp. $\beta_w^{*}(W\an)$, $\beta_{w,n}^*(W\an)$)
is a $\bold{PH\ii}(\B{R})^0$- (resp. $\bold{H}(K_w)$-)equivariant analytic vector bundle on $B^{d-1}$
(resp. $\Si$, $\SiN$).

For each $\SIE$ (resp. $S\in\ff{E\ii}$) consider a double quotient
\qquad\qquad\qquad\linebreak
$\wt{V}_S:=S\bs [\beta_w^{*}(W\an) \times E']/\Gm$ (resp. $\wt{\wt{V}}\ii_S:=
S\bs [\beta_{\B{R}}^{*}(W\tw_{\B{C}})\an \times E\ii]/\Gm\ii$).
\end{Emp}

\begin{Prop} \label{P:anv}
For each $S\in \ff{E}$ (resp. $S\in\ff{E\ii}$) $\:\wt{V}_{S}$ (resp $\wt{\wt{V}}\ii_S$)
has a natural structure of an affine scheme $V_S$ over $X_S$ (resp. $\wt{V}\ii_S$ over $\wt{X}\ii_S$). 
Moreover, $V_S$ (resp. $\wt{V}\ii_S$) is a vector bundle on $X_S$  (resp. $\wt{X}\ii_S$) if
$S$ is sufficiently small.
\end{Prop}

\begin{pf}
We give the proof in the $p$-adic case. The complex case is similar, but easier.

I) First we take $S$ of the form $T_n\times S'$ 
with sufficiently small $S'\in\ff{E'}$. Then 
$\wt{V}_{S}$ is a finite disjoint union of quotients of the form 
$\Gm_{aS'a^{-1}}\bs\beta_{w,n}^{*}(W\an)$ with some $a\in E'$.
Since  the projection $\SiN\to \Om$ factors through each 
$\Gm_{aS'a^{-1},0}\bs\SiN$ (in the notation of the proof of \rp{sav2}), the quotient $\Gm_{aS'a^{-1},0}\bs\beta_{w,n}^{*}(W\an)$ 
is naturally an analytic vector bundle on $\Gm_{aS'a^{-1},0}\bs\SiN$.
Now (as in the proof of \rp{sav2}) the quotient vector bundle  
$P\Gm_{a\bar{S}a^{-1}}\bs(\Gm_{aS'a^{-1},0}\bs\beta_{w,n}^{*}(W\an))
\cong\Gm_{aS'a^{-1}}\bs\beta_{w,n}^{*}(W\an)$ on $\Gm_{aS'a^{-1}}\bs\SiN$ is obtained by gluing.
For the algebraization we use \rcr{GAGA2} a).

II) For each $T\in\ff{E}$ there exists a normal subgroup of the form
$S=T_n\times S'$, where $S'\in\ff{E'}$ is sufficiently small. Then by the same
considerations as in \rp{sav1}, $V_T$ can be defined as $(T/S)\bs V_S$
(using \rcr{GAGA2}, a)).

III) Suppose that $V_{S_1}$ and $V_{S_2}$, constructed in I) and II), are vector bundles on $
X_{S_1}$ and $X_{S_2}$ respectively for $S_1\subset S_2$ in $\ff{E}$. Then
the natural morphism $f:V_{S_1}\to V_{S_2}\times_{X_{S_2}}X_{S_1}$
of vector bundles on $X_{S_1}$ induces an isomorphism on each fiber. 
Hence it is an isomorphism.

IV) Suppose that $T\subset S$ in $\ff{E}$ and that $V_S$
is a vector bundle on $X_S$. Choose  a normal subgroup $S_0\in\ff{E}$ of $T$ 
such that $V_{S_0}$ is a vector bundle on $X_{S_0}$. Then
$V_T=(T/S_0)\bs V_{S_0}\cong (T/S_0)\bs V_S\times_{X_S}X_{S_0}\cong
V_S\times_{X_S}((T/S_0)\bs X_{S_0})\cong V_S\times_{X_S}X_T$, so $V_T$ is a vector
bundle on $X_T$.    
\end{pf}

\begin{Emp} \label{E:alv}
Choose $\SIE$ (resp. $S\in\ff{E\ii}$) sufficiently small. Then $V_S$
(resp. $\wt{V}_S\ii$) is a vector bundle on $X_S$ (resp. $\wt{X}_S\ii$).
Thus $V:=V_S\times_{X_S}X$ (resp. $\wt{V}\ii:=\wt{V}_S\ii\times_{\wt{X}\ii_S}\wt{X}\ii$)
is a vector bundle on $X$ (resp. $\wt{X}\ii$). By Step III)
of the proof, $V$ (resp. $\wt{V}\ii$) does not depend on $S$.
Each $g\in E$ (resp. $g\in E\ii$) defines an isomorphism 
$V\an_S\isom V\an_{gSg^{-1}}$ (resp. $(\wt{V}\ii_S)\an\isom(\wt{V}\ii_{gSg^{-1}})\an$).
Therefore by \rcr{GAGA2} a), $g$ defines an isomorphism $V_S\isom V_{gSg^{-1}}$
(resp. $\wt{V}\ii_S\isom \wt{V}\ii_{gSg^{-1}}$). The product of this isomorphism and 
the action of $g$ on $X$ (resp. $\wt{X}\ii$) gives us an isomorphism
$g:V=V_S\times_{X_S}X\isom V_{gSg^{-1}}\times_{X_{gSg^{-1}}}X=V$
(resp. $g:\wt{V}\ii\isom \wt{V}\ii$). Thus we have constructed an algebraic action of $E$
(resp. of $E\ii$) on $V$ (resp. $\wt{V}\ii$), satisfying $S\bs V\cong V_S$ for all
$S\in\ff{E}$  (resp. $S\bs \wt{V}\ii\cong \wt{V}\ii_S$ for all
$S\in\ff{E\ii}$). Moreover, $V=\underset{\underset{S}{\longleftarrow}}{\lim}\,V_S$
and $\wt{V}\ii=\underset{\underset{S}{\longleftarrow}}{\lim}\,\wt{V}\ii_S$.

By \cite{Mi}, there exists a unique canonical model $V\ii$
of $\wt{V}\ii$ over $K_w$ (the definition of the canonical model will
be explained in the last paragraph of the proof of \rp{eq}) such that $V\ii$ is an 
$E\ii$-equivariant vector bundle on $X\ii$.
\end{Emp}

Our main task is to prove the following

\begin{Maint} \label{M:th}
 For each admissible isomorphism $\Phi:E\isom E\ii$, each isomorphism $f_{\Phi}$
from the First or the Second Main Theorem can be lifted to a 
$\Phi$-equivariant isomorphism $f_{\Phi,V}:V\isom V\ii$.
\end{Maint}
We will prove this theorem, using standard principal bundles(=torsors) 
(see \cite[Ch. III, $\S$3]{Mi}).

\subsection{Equivariant torsors} \label{SS:contor}

\begin{Emp} \label{E:not}
For each $\SIE$ (resp. $S\in\ff{E\ii}$) consider the double quotient
$\wt{P}_S:=S\bs [\Si\times (\bold{PH_{K_w}})\an \times E']/\Gm$ (resp.
$\wt{P}\ii_S:=
S\bs [B^{d-1}\times (\bold{PH\ii}_{\B{C}})\an\times E\ii]/\Gm\ii$).
\end{Emp}

\begin{Prop} \label{P:altor}
For each $\SIE$ (resp. $S\in\ff{E\ii}$) $\wt{P}_S$ (resp. $\wt{\wt{P}}_S\ii$)
has a natural structure of an affine scheme $P_S$ over $X_S$ (resp. $\wt{P}_S\ii$ over $\wt{X}_S\ii$). 
Moreover, $P_S$ is a $\bold{PH_{K_w}}$-torsor over $X_S$ (resp. $\wt{P}_S\ii$ is a 
$\bold{PH\ii}_{\B{C}}$-torsor over $\wt{X}_S\ii$) if $S$ is sufficiently small.
\end{Prop}

\begin{pf}
is almost identical to that of \rp{anv} (using \rp{torsor} and \rl{tor} instead of \rco{GAGA2} a) and 
arguments of step III) respectively).
\end{pf}

\begin{Emp}
Arguing as in \re{alv} and using \rco{torsor} we obtain an $E$-equivariant 
$\bold{PH_{K_w}}$-torsor $P=\underset{\underset{S}{\longleftarrow}}{\lim}\,P_S$ 
over $X$  (resp. an $E\ii$-equivariant $\bold{PH\ii}_{\B{C}}$-torsor
$\wt{P}\ii=\underset{\underset{S}{\longleftarrow}}{\lim}\,\wt{P}\ii_S$ over 
$\wt{X}\ii$). By \cite[III, Thm. 4.3]{Mi}, there exists a unique 
canonical model $P\ii$ of $\wt{P}\ii$ over $K_w$ (the definition
will be explained in \rco{canm}) such that $P\ii$ is an $E\ii$-equivariant
$\bold{PH_{K_w}\ii}$-torsor over $X\ii$. Let $\pi:P\to X$ and $\pi\ii:P\ii\to X\ii$ be
the natural projections. Denote also the natural projection from the
$\bold{PH_{K_w}\ii}$-torsor $P\tw$ to $X$ by $\pi\tw$.
\end{Emp}

\begin{Mainfo} \label{M:fo}
 For each admissible isomorphism $\Phi :E\isom E\ii$, each isomorphism $f_{\Phi}$
from the First or the Second Main Theorems can be lifted to a 
$\Phi$-equivariant isomorphism $f_{\Phi,P}:P\tw\isom P\ii$ of 
$\bold{PH_{K_w}\ii}$-torsors.
\end{Mainfo}

\subsection{Connection between the Main Theorems} \label{SS:conmain}

\begin{Prop} \label{P:eq}
The Fourth Main Theorem implies the Third one.
\end{Prop}

\begin{pf}
Consider the pro-analytic maps $\wt{\rho}':[\Si\times (\bold{PH_{K_w}})\an\times
(E')^{disc}]/\Gm\to(\B{P}_{K_w}^{d-1})\an$ and $(\wt{\wt{\rho}}\,')\ii:[B^{d-1}
\times(\bold{PH\ii_{\B{C}}})\an\times(E\ii)^{disc}]/\Gm\ii\to(\B{P}_{\B{C}}^{d-1})\an$ given by
$\wt{\rho}'(x,g,e)=g\beta_w(x)$ and $(\wt{\wt{\rho}}\,')\ii(x,g,e)=g\beta_{\B{C}}(x)$.
Then $\wt{\rho}'$ (resp. $(\wt{\wt{\rho}}\,')\ii$) is $(\bold{PH_{K_w}})\an$- 
(resp. $(\bold{PH\ii}_{\B{C}})\an$-) equivariant and commutes with the action of $E$ 
(resp. $E\ii$). Hence it defines an equivariant analytic map $\wt{\rho}:P\an\to
(\B{P}_{K_w}^{d-1})\an$ 
(resp. $\wt{\wt{\rho}}\ii:(\wt{P}\ii)\an\to(\B{P}_{\B{C}}^{d-1})\an$).

\begin{Prop} \label{P:mor}
There exists a unique algebraic morphism $\rho:P\to \B{P}_{K_w}^{d-1}$ (resp. 
$\wt{\rho}\ii:\wt{P}\ii\to \B{P}_{\B{C}}^{d-1}$) such that 
$\rho\an\cong\wt{\rho}$ (resp. $(\wt{\rho}\ii)\an\cong\wt{\wt{\rho}}\ii$).
\end{Prop}

\begin{pf}
We prove the statement for $\rho$ (in the second case the proof is 
exactly the same). We have to show that the graph $Gr(\wt{\rho})\subset
P\an\times (\B{P}_{K_w}^{d-1})\an$ corresponds to an algebraic subscheme. 
For each $\SIE$ let $\wt{\rho}_S:P_S\an\to (\B{P}_{K_w}^{d-1})\an$ be the morphism
induced by $\wt{\rho}$. Since $Gr(\wt{\rho})=\underset{\underset{S}
{\longleftarrow}}{\lim}\,Gr(\wt{\rho}_S)\subset (\underset{\underset{S}
{\longleftarrow}}{\lim}\,P_S\an)\times (\B{P}_{K_w}^{d-1})\an$, 
it remains to show that the graph $Gr(\wt{\rho}_S)\subset
P_S\an\times (\B{P}_{K_w}^{d-1})\an$ corresponds to a unique algebraic subvariety
for each $S$ sufficiently small. 

Take $S$ so small that $X_S$ is smooth, then by \rl{descent} b) there exists
a quotient  $Q_S:=\bold{PH_{K_w}}\bs(P_S\times\B{P}_{K_w}^{d-1})$ by the
diagonal action of $\bold{PH_{K_w}}$. Moreover, $Q_S$ is a $\B{P}^{d-1}$-bundle
on $X_S$, hence it is projective over $K_w$. Let $\al:P_S\times\B{P}_{K_w}^{d-1}\to Q_S$ be 
the natural projection. Since $\wt{\rho}_S$ is 
$(\bold{PH_{K_w}})\an$-equivariant, $Gr(\wt{\rho}_S)$ is invariant under the 
diagonal action of $(\bold{PH_{K_w}})\an$. Therefore the quotient
 $\wt{Q}:=(\bold{PH_{K_w}})\an\bs Gr(\wt{\rho}_S)$ is a closed
analytic subspace of $Q_S\an$, so that it is algebraic (see \rco{GAGA1}).
It follows that its inverse image $\al^{-1}(\wt{Q})=Gr(\wt{\rho}_S)$
is also algebraic. The uniqueness is clear.
\end{pf}

\begin{Cl} \label{C:un}
The map $\wt{\wt{\rho}}\ii$ is the only $(\bold{PH\ii}_{\B{C}})\an\times E\ii$-equivariant
analytic map from $(\wt{P}\ii)\an$ to $(\B{P}_{\B{C}}^{d-1})\an$.
\end{Cl}

\begin{pf}
Let $\rho':(\wt{P}\ii)\an\to (\B{P}_{\B{C}}^{d-1})\an$ be any such map. 
Composing it with the natural $(\bold{PH\ii}_{\B{C}})\an\times E\ii$-equivariant 
projection
$[B^{d-1}\times (\bold{PH\ii}_{\B{C}})\an\times(E\ii)^{disc}]/\Gm\ii\to(P\ii)\an$
and identifying a complex analytic space with the set of its $\B{C}$-rational
points, we obtain a $\bold{PH\ii}(\B{C})\times E\ii$-equivariant analytic map 
$\rho'':[B^{d-1}\times\bold{PH\ii}({\B{C}})\times(E\ii)^{disc}]/\Gm\ii\to
(\B{P}_{\B{C}}^{d-1})\an$.
Let $\rho_0$ be the restriction of $\rho$ to $B^{d-1}\isom B^{d-1}\times\{1\}\times\{1\}$. 
Then $\rho''(x,g,e)=g\rho_0(x)$ for all $x\in B^{d-1},\;g\in(\bold{PH\ii}_{\B{C}})\an$ and $e\in E\ii$. 
Therefore $\gm\rho_0(x)=\rho_0(\gm x)$ for all $\gm\in\Gm\ii$ and $x\in B^{d-1}$.
Since the subgroup $\Gm\ii$ is dense in $\pgr^0$, we obtain by continuity that
$\gm\rho_0(x)=\rho_0(\gm x)$ for all $\gm\in\pgr^0$ and $x\in B^{d-1}$. 
In particular, for the origin $0\in B^{d-1}$ we get $Stab_{\pgr^0}(0)\subset 
Stab_{\pgr^0}(\rho_0(0))$. The subgroup $Stab_{\pgr^0}(0)$ stabilizes precisely one 
point $(0:... :0:1)\in\B{P}^{d-1}(\B{C})$ if $d>2$ and two points $(0:1)$ and 
$(1:0)$ in $\B{P}^1(\B{C})$ if $d=2$. The case $\rho_0(0)=(1:0)$ is 
impossible,
because  identifying $\B{P}^1(\B{C})$ with $\overline{\B{C}}=\B{C}\cup\infty$ by
$(x:y)\mapsto x/y$ we would get in this case $\rho_0(z)=1/\bar{z}$ for all $z\in B^1$, contradicting
the analyticity of $\rho_0$. We conclude that $\rho_0(0)=(0:...:0:1)$. 
Hence $\rho_0=\beta_{\B{R}}$ and $\rho'=\wt{\wt{\rho}}\ii$.
\end{pf}

\begin{Emp} \label{E:rho}
Next we show that the map $\wt{\rho}\ii:P\ii_{\B{C}}\to 
(\B{P}_{K_w}^{d-1})\tw_{\B{C}}$ is $K_w$-rational. Recall that the map $\wt{\rho}\ii$
is $\bold{PH\ii}_{\B{C}}$-equivariant and that the actions of the group 
$\bold{PH\ii_{K_w}}$ on both $P\ii$ and $(\B{P}_{K_w}^{d-1})\tw$ are 
$K_w$-rational.
Therefore for each $\sigma\in\Aut(\B{C}/K_w)$ the analytic map 
${\sigma}(\wt{\rho}\ii)\an$ is $(\bold{PH\ii}_{\B{C}})\an\times E\ii$-equivariant, 
hence it coincides with $(\wt{\rho}\ii)\an=\wt{\wt{\rho}}\ii$. By the uniqueness of the algebraic structure,
${\sigma}(\wt{\rho}\ii)=\wt{\rho}\ii$.

It follows that $\wt{\rho}\ii$ defines a $\bold{PH_{K_w}\ii}
\times E\ii$-equivariant map $\rho\ii:P\ii\to (\B{P}_{K_w}^{d-1})\tw$. 
Notice also that $\rho$ defines a $\bold{PH_{K_w}\ii}\times E$-equivariant map 
$\rho\tw:P\tw\to (\B{P}_{K_w}^{d-1})\tw$.
\end{Emp}

Suppose that the Fourth Main Theorem holds, then

\begin{Lem} \label{L:eqv}
We have $\rho\ii\circ f_{\Phi,P}=\rho\tw$.
\end{Lem}

\begin{pf} 
By the claim, $\wt{\wt{\rho}}\ii\cong(\rho_{\B{C}}\ii)\an$ is
equal to $(\rho\tw_{\B{C}}\circ(f_{\Phi,P})_{\B{C}}^{-1})\an$.   
From the uniqueness of algebraic structures we conclude that 
$\rho_{\B{C}}\ii=\rho\tw_{\B{C}}\circ (f_{\Phi,P})_{\B{C}}^{-1}$.
Now we descend to $K_w$ as in \re{rho}.
\end{pf}

It follows from the definitions that $\rho^*(W)\cong\pi^*(V)$ 
(hence $\rho^*(W)\tw\cong(\pi\tw)^*(V)$) and $(\wt{\rho}\ii)^*(W\tw
_{\B{C}})\cong (\pi\ii_{\B{C}})^*(\wt{V}\ii)$. \rl{descent} allows us to
define $V\ii$ by the requirement that $(\pi\ii)^*(V\ii)\cong (\rho\ii)^*(W\tw)$.
(By the definition, this is the canonical model of $\wt{V}\ii$ on $X\ii$).
\rl{eqv} implies that $f_{\Phi,P}$ can be lifted to the $\Phi$-equivariant
isomorphism $\rho^*(W)\tw\cong(\rho\tw)^*(W\tw)\isom(\rho\ii)^*(W\tw)$, 
commuting with the $\bold{PH\ii_{K_w}}$-action. This gives us the
$\bold{PH_{K_w}\ii}$-equivariant isomorphism $(\pi\tw)^*(V)\isom(\pi\ii)^*(V\ii)$.
Hence the Third Main Theorem follows from \rl{descent}.
\end{pf}

\begin{Rem} \label{R:impl}
Tannakian arguments can be used to show (see \rt{torsor} and the discussion 
around it) that 
the Third Main Theorem implies the Fourth one. We will not use this implication.
\end{Rem}

\subsection{Reduction of the problem} \label{SS:redpr}

\begin{Emp} \label{E:pf}
Now we start the proof of the Fourth Main Theorem. For simplicity of notation we
identify $E$ with $E\ii$ by means of $\Phi$ and $X$ with $X\ii$ by means of 
$f_{\Phi}$.
%
%
Recall that $P_{S,\B{C}}$ is a 
$\PH_{\B{C}}$-torsor over $X_{S,\B{C}}$ for all  sufficiently small $\SIE$,
hence $(P_{S,\B{C}})\an$ is a $(\PH_{\B{C}})\an$-torsor over 
$(X_{S,\B{C}})\an$ and $(P_{\B{C}})\an=(P_{S,\B{C}})\an\times _{(X_{S,\B{C}})\an}(X_{\B{C}})\an$ is a 
$(\PH_{\B{C}})\an$-torsor over 
$(X_{\B{C}})\an\cong[B^{d-1}\times (E\ii/E\ii_0)^{disc}]/P\Gm\ii$. Set 
$Y:=(\pi\an)^{-1}(B^{d-1}\times\{1\})\subset (P_{\B{C}})\an$. Then $Y$ is a 
$(\PH_{\B{C}})\an$-torsor over $B^{d-1}$. Recall that $E_0=E\ii_0$ acts 
trivially on $P$, hence $(P_{\B{C}})\an\cong(Y\times (E\ii/E\ii_0)^{disc})
/P\Gm\ii$.
\end{Emp}

\begin{Prop} \label{P:decomp}
There exists a homomorphism $j:P\Gm\ii\to\PH(\B{C})$ and an isomorphism
$(P_{\B{C}})\an\isom (B^{d-1}\times (\PH_{\B{C}})\an\times (E\ii/E\ii_0)^{disc})/
P\Gm\ii$ such that $(x,h,g)\gm=(\gm^{-1}_{\be}x,hj(\gm),g\gmE)$ for all 
$x\in B^{d-1},\;h\in (\PH_{\B{C}})\an,\; g\in E\ii/E\ii_0$ and $\gm\in P\Gm\ii$.
\end{Prop}

\begin{pf}
The proposition asserts that there exists a decomposition 
$Y\cong B^{d-1}\times (\PH_{\B{C}})\an$ such that the group $P\Gm\ii$ acts on 
$B^{d-1}\times (\PH_{\B{C}})\an$ by the product of actions on factors. 

The trivial connection on $\Si\times (\bold{PH_{K_w}})\an\to\Si$ is 
$\Gm\times\wt{D}_w\m$-invariant, therefore it defines a natural $E$-invariant flat
connection $\wt{\C{H}}$ on the $(\bold{PH_{K_w}})\an$-torsor 
$[\Si\times(\bold{PH_{K_w}})\an\times E^{disc}]/\Gm$ over 
$[\Si\times E^{disc}]/\Gm$. Since for all sufficiently small $\SIE$
the projection $(\Si\times(\bold{PH_{K_w}})\an\times E^{disc})/\Gm\to P\an_S$ 
is \'etale, it induces an isomorphism of tangent spaces up to an extension 
of scalars. Hence $\wt{\C{H}}$ induces a flat connection 
$\wt{\C{H}}_S$ on $P\an_S$. By the definition, $\wt{\C{H}}_S$ is a 
$(\bold{PH_{K_w}})\an$-invariant analytic vector subbundle of $(T_{P_S})\an$, 
therefore \rl{descent} and \rco{GAGA2} imply the existence of a unique flat connection
$\C{H}_S$ on $P_S$ such that $\wt{\C{H}}_S\cong \C{H}_S\an$. 
Since the projection $\pi_S:P\to P_S$ is \'etale, $\C{H}_S$ defines a unique flat
connection $\C{H}$ on $P$ satisfying $(\pi_S)_*(\C{H})=\C{H}_S$. Moreover, 
$\C{H}$ is $E$-equivariant and does not depend on $S$. 

The connection $\C{H}$ determines flat connections $\C{H}_{\B{C}}$ on 
$P_{\B{C}}$ and $(\C{H}_{\B{C}})\an$ on $(P_{\B{C}})\an$. Let $\C{H}'$ be the 
restriction of $(\C{H}_{\B{C}})\an$ to $Y$. Then $\C{H}'$ is a
$P\Gm\ii$-invariant flat connection on the 
$(\PH_{\B{C}})\an$-principal bundle $Y$ over the simply connected complex manifold
$B^{d-1}$. By \rl{trivcon}, there exists a decomposition $Y\cong B^{d-1}\times (\PH_{\B{C}})\an$
such that the corresponding action of $ P\Gm\ii$  on $B^{d-1}\times (\PH_{\B{C}})\an$
preserves the trivial connection. 

For each $\gm\in P\Gm\ii$ let 
$\wt{\gm}:B^{d-1}\times(\PH_{\B{C}})\an\to(\PH_{\B{C}})\an$ be the analytic map such that
$\gm(x,h)=(\gm(x),\wt{\gm}(x,h))$ for all $x\in B^{d-1}$ and $h\in(\PH_{\B{C}})\an$.
Since the action of $P\Gm\ii$ preserves the trivial connection, we have 
$\partial\wt{\gm}/\partial x\equiv 0$ for each $\gm\in P\Gm\ii$. Hence analytic   
$\wt{\gm}$'s depend only on $h$. Since the action of $P\Gm\ii$ commutes 
with the action of $(\PH_{\B{C}})\an$, we have $\wt{\gm}(h)=h\wt{\gm}(1)$ for all $h\in(\PH_{\B{C}})\an$
and $\gm\in P\Gm\ii$. Therefore the map $\gm\mapsto\wt{\gm}(1)^{-1}$ is the required homomorphism.
\end{pf}

\begin{Thm} \label{T:diag}
There exists an inner isomorphism (=inner twisting) $\Phi_{\B{C}}:\PH_{\B{C}}\isom \bold{PH\ii}_{\B{C}}$
such that $j\circ\Phi_{\B{C}}:P\Gm\ii=\bold{PG\ii}(F)_+\to \bold{PH\ii}(\B{C})\cong
\bold{PG\ii}(\B{C}\otimes_{\B{Q}}F)$ is induced by the natural (diagonal) embedding 
$F\hra \B{C}\otimes_{\B{Q}}F\cong \B{C}^g$.
\end{Thm}

\begin{Rem}
Algebraization considerations as in \rl{comp} (using \rp{torsor} instead of 
\rco{GAGA1}) show that \rt{diag} implies the existence of a $\Phi$-equivariant 
isomorphism $P\tw_{\B{C}}\isom\wt{P}\ii$, lifting $\wt{f}_{\Phi}$.
\end{Rem}

\subsection{Proof of density} \label{SS:density}

To prove \rt{diag} we will use Margulis' results. For this we first show that 
the subgroup $j(P\Gm\ii)$ is sufficiently large.
We start with the following technical

\begin{Lem} \label{L:comm}
let $n$ and $d$ be positive integers. For each $i=1,...,n$ we denote by $\pr_i$
the projection to the $i$-th factor.

a) Let $\C{G}_1,...,\C{G}_n$ be Lie algebras, and let $\C{H}$ be an ideal in 
the Lie algebra $\C{G}=\prod_{i=1}^n \C{G}_i$.
Then $\C{H}\supset\prod_{i=1}^n [\pr_i\,\C{H},\C{G}_i]$.\par

b) Let $\Dt$ be a subgroup of $\PGL{d}(\B{C})^n$. Suppose that $\pr_i(\Dt)$
is infinite for every $i=1,...,n$. If $\wt{\Dt}:=Comm_{\PGL{d}(\B{C})^n}(\Dt)$ is 
Zariski dense in $(\PGL{d})^n$, then the same is true for $\Dt$.\par

c) If a subgroup $\Dt\subset \PGU{d}(\B{R})^n$ is Zariski dense (in 
$(\PGU{d})^n$), then it is dense.
\end{Lem}

\begin{pf}
a) If $x=(x_1,...,x_n)\in\prod_{i=1}^n\C{G}_i$ belongs to $\C{H}$, then $[x,y_i]=(0,...,
[x_i,y_i],...,0)=[\pr_i\,x,y_i]\in\C{H}$ for all $y_i\in \C{G}_i$.\par

b) Let $\bold{J}$ be the Zariski closure of $\Dt$ in $(\PGL{d})^n$, then 
$\dt\bold{J}\dt^{-1}\cap\bold{J}$ is an 
algebraic subgroup of finite index in $\bold{J}$ for each $\dt\in\wt{\Dt}$. 
Hence $\dt \bold{J}^0\dt^{-1}=\bold{J}^0$. In particular, the subgroup
$\Ad\,\wt{\Dt}$ 
stabilizes $\Lie\,\bold{J}^0\subset\Lie(\PGL{d})^n$. Since $\wt{\Dt}$ is Zariski dense 
in $(\PGL{d})^n$, the Lie algebra $\Lie\,\bold{J}^0=\Lie\,\bold{J}$ is an ideal in 
$\Lie(\PGL{d})^n$. By our assumption, $\pr_i(\bold{J})$ is an infinite algebraic group 
for each $i=1,...,n$, therefore $\pr_i(\Lie\,\bold{J})\neq 0$ is an ideal in a simple Lie 
algebra $\Lie(\PGL{d})$. Therefore a) implies that $\Lie\,\bold{J}=\Lie(\PGL{d})^n$. Since 
the group $(\PGL{d})^n$ is connected, $\bold{J}=(\PGL{d})^n$.
\par

c) Let $M$ be the closure of $\Dt$ in $\PGU{d}(\B{R})^n$. Then $M$ is a Lie 
subgroup of the Lie group $\PGU{d}(\B{R})^n$. Hence  $\Lie\,M$ is an 
$\Ad\,M$-invariant subspace of $\Lie(\PGU{d}(\B{R}))^n$. Since the adjoint 
representation is algebraic, $\Lie\,M$ is an ideal in $\Lie(\PGU{d}(\B{R}))^n$. 
Since $M$ is compact, it has a finite number of connected components.
Hence $M^0$ is also Zariski dense, therefore it is not contained in 
$\PGU{d}(\B{R})^{i-1}\times\{1\}\times\PGU{d}(\B{R})^{n-i}$ for
any $i=1,...,n$. It follows that $\Lie\,M=\Lie\,M^0$ is not contained in 
$\Lie(\PGU{d}(\B{R}))^{i-1}\times\{0\}\times\Lie(\PGU{d}(\B{R}))^{n-i}$, so that
$\pr_i(\Lie\,M)\neq 0$. Now the assertion follows exactly in the same way as in b).
\end{pf}
 
\begin{Prop} \label{P:Zar}
The subgroup $j(P\Gm\ii)$ is Zariski dense in $\PH_{\B{C}}$.
\end{Prop}

\begin{pf}
Let $\bold{G'}\subset \PH_{\B{C}}$ be the Zariski closure of $j(P\Gm\ii)$. Then 
$\wt{R}:=(B^{d-1}\times (\bold{G'})\an\times(\bold{PG\ii}(\afv))^{disc})/P\Gm\ii$ 
is a $\bold{PG\ii}(\afv)$-invariant $(\bold{G'})\an$-subtorsor of the 
$(\PH_{\B{C}})\an$-torsor $(P''_{\B{C}})\an=((\bold{G\ii}(F_v)\times Z(E\ii))\bs 
P_{\B{C}})\an\cong[B^{d-1}\times(\PH_{\B{C}})\an\times(\bold{PG\ii}(\afv))^{disc}]
/P\Gm\ii$ over $(X''_{\B{C}})\an\cong [B^{d-1}
\times(\bold{PG\ii}(\afv))^{disc}]/P\Gm\ii$. Hence by \rp{torsor} there exists an 
algebraic $\bold{G'}$-subtorsor $R$ of $P''_{\B{C}}$ such that $R\an\cong\wt{R}$. 
Using our identification of $\B{C}_p$ with $\B{C}$, we obtain a closed analytic 
subspace $(\bold{R}_{\B{C}_p})\an\subset(P''_{\B{C}_p})\an\cong
(\Om\Hat{\otimes}_{K_w}\B{C}_p\times(\PH_{\B{C}_p})\an\times(PE')^{disc})/P\Gm$. 
Recall that $P\Gm=\PH(\B{Q})$ is naturally embedded into $\PH(\B{C}_p)$.

\begin{Lem} \label{L:Zar}
The subgroup generated by the elements of $P\Gm$ with elliptic projections to
$\pgp$ is Zariski dense in $\PH_{\B{C}_p}$.
\end{Lem}

\begin{pf}
The subgroup of $\pgp$ generated by the set of all elliptic elements is open
and normal, because a conjugate of an elliptic element is elliptic. Hence it 
contains $\psp$. The subgroup $P\GmG\cap\psp$ is dense in
$\psp$. Therefore by \cite[Ch. IX, Lem. 3.3]{Ma}, the subgroup of $P\GmG$ 
generated by all elliptic elements of $P\GmG$ contains $P\GmG\cap\psp$. 
In particular, it has finite index in $P\GmG=\PH(\B{Q})$. Since $\PH$ is connected,
the statement follows from \cite[Ch. V, Cor. 18.3]{Bo}.
\end{pf}

If $\bold{G'}\neq \PH_{\B{C}}$, then by the lemma there exists $\gm\in P\Gm$ with elliptic projection to $\pgr$ whose image $\gm_p\in\PH(\B{C}_p)$ does not
belong to $\bold{G'}(\B{C}_p)$. Let $x\in\Om\Hat{\otimes}_{K_w}\B{C}_p\times\{1\}
\subset (X''_{\B{C}_p})\an$ be an elliptic point of $\gm_E\in PE'$, and let $\wt{x}$
be an arbitrary point of $(R_{\B{C}_p})\an$, lying over $x$. Then 
$\gm_E(\wt{x})=\gm_p(\wt{x})$ is another point of $(R_{\B{C}_p})\an$, 
lying over $x$. Hence $\gm_p$ must belong to $\bold{G'}(\B{C}_p)$, contradicting to our choice
of $\gm$.
\end{pf}

Recall that we defined  in \rp{mor} the algebraic $\bold{PH_{K_w}}\times E$-equivariant map 
$\rho:P\to\B{P}_{K_w}^{d-1}$. Identify $\bold{PH}_{\B{C}}$ with $(\PGL{d})^g$ in
such a way  that the first factor corresponds to the embedding $\be_1:K\hra\B{C}$.
Denote by $j_k:P\Gm\ii\to \PGL{d}(\B{C})$ the composition of $j$ with the projection $\pr_k$ of 
$\PGL{d}(\B{C})^g$ to its 
$k$-th factor. Denote also by $\pr^i$ the projection of $\PGL{d}(\B{C})^g$ 
to the product of all factors except the $i$-th. We will sometimes identify $P\Gm\ii$ with its projection
$P\Gm\ii_{\be}\subset\pgr^0$.

\begin{Prop} \label{P:compact}
The subgroup $j_1(P\Gm\ii)$ is not relatively compact in $\PGL{d}(\B{C})$.
\end{Prop}

\begin{pf}
If not, then $j_1(P\Gm\ii)$ is contained in some maximal compact subgroup of 
$\PGL{d}(\B{C})$ (see for example \cite[Prop. 3.11]{PR}). After a suitable 
conjugation we may assume that $j_1(P\Gm\ii)\subset \PGU{d}(\B{R})$. By \rp{Zar},
 $j_1(P\Gm\ii)$ is Zariski dense in $\PGL{d}$, hence it is infinite. Therefore,
by \rl{comm}, $j_1(P\Gm\ii)$ is dense in $\PGU{d}(\B{R})$.

Consider the map $\rho:P(\B{C})\to\B{P}^{d-1}(\B{C})$ and its restriction $\rho_0$ to
$B^{d-1}\times\{1\}\cong B^{d-1}$. Then, as in the proof of \rcl{un}, 
$\rho_0:B^{d-1}\to\B{P}^{d-1}(\B{C})$ satisfies 
$\rho_0(\gm x)=j_1(\gm)\rho_0(x)$ for all $x\in B^{d-1}$ and $\gm\in P\Gm\ii$.
The group $\PGU{d}(\B{R})$ acts transitively on $\B{P}^{d-1}(\B{C})$, hence
$\rho_0(B^{d-1})$ is dense in $\B{P}^{d-1}(\B{C})$.

Now we want to prove that $j_1:P\Gm\ii\to\PGU{d}(\B{R})$ can be extended to a
continuous homomorphism $\wt{j}_1:\pgr^0\to\PGU{d}(\B{R})$.
For each $g\in\pgr^0$ choose a sequence $\{\gm_n\}_n\subset P\Gm\ii\subset\pgr^0$ 
converging to $g$. Since $\PGU{d}(\B{R})$ is compact, there exists a subsequence $\{\gm_{n_k}\}_k
\subset \{\gm_n\}_n$ such that $\{j_1(\gm_{n_k})\}_k$ converges to some $a\in\PGU{d}(\B{R})$.
Then $\rho_0(gx)=\underset{k}{\lim}\,\rho_0(\gm_{n_k}(x))=(\underset{k}{\lim}\,
j_1(\gm_{n_k}))\rho_0(x)=a\rho_0(x)$ for all $x\in B^{d-1}$.
In follows that $a=a(g)$ depends only on $g$,
since $\rho_0(B^{d-1})$ is dense in $\B{P}^{d-1}(\B{C})$
and since the group $\PGL{d}(\B{C})$ acts faithfully on $\B{P}^{d-1}(\B{C})$.
In particular, $a(g)$ does not depend on the choice of  $\{\gm_n\}_n$,
and $a(g)=\underset{n}{\lim}\,j_1(\gm_n)$. It follows that $\wt{j}_1:=a$ is the required
extension.

Since $\pgr^0$ is simple and $j_1(P\Gm\ii$) is dense, $\wt{j}_1$ must be injective
and surjective. Hence is it an isomorphism, a contradiction.
\end{pf}

\begin{Prop} \label{P:inj}
For each $i=1,...,g$ the homomorphism $j_i:P\Gm\ii\to \PGL{d}(\B{C})$ is 
injective.
\end{Prop}

\begin{pf}
Suppose that for some $i$ the subgroup $\Dt_i:=\Ker(j_i)$ is non-trivial. 
Then $\Dt_i$ is a normal subgroup of $P\Gm\ii\cong P\Gm\ii_{\be}\subset\pgr^0$. Hence 
the closure of $(\Dt_i)_{\be}$ is a non-trivial normal
subgroup of a simple group $\pgr^0$. Therefore the projection $(\Dt_i)_{\be}$ is dense in $\pgr^0$.
Hence there exists an element $\dt\in\Dt_i$ with elliptic projection 
$\dt_{\be}\in\pgr^0$. Therefore the element $(j(\dt),\dt_E)\in\PGL{d}(\B{C})^g
\times PE'$ has a fixed point $[y,g,e]\in P''(\B{C}_p)=(\Om(\B{C}_p)\times \PGL{d}(\B{C}_p)^g
\times PE')/P\Gm$. Hence $(g^{-1}j(\dt)g,e^{-1}\dt_E e)$ 
stabilizes $[y,1,1]\in P''(\B{C}_p)$. It follows that $ e^{-1}\dt_E e=\gmE$ for some 
$\gm\in P\Gm=\bold{PH}(\B{Q})$ and that $g^{-1}j(\dt)g\in\bold{PH}(\B{C}_p)$ is the image of $\gm$.
Hence $j_k(\dt)\neq 1$ for all $k$, contradicting to our assumption.
\end{pf}

\subsection{Use of rigidity} \label{SS:rig}

Now we are going to use the following theorem of Margulis \cite[Ch. VII, Thm. 5.6]
{Ma}.

\begin{Thm} \label{T:super}
Let $L$ be a local field, let $\bold{J}$ be a connected absolutely simple adjoint
$L$-group, and let $A$ be a finite set. For each $\al\in A$ let $k_{\al}$ be a local
field and let $\bold{G_{\al}}$ be an adjoint absolutely simple $k_{\al}$-isotropic
group. Set $G:=\underset{\al\in A}{\prod}\bold{G_{\al}}(k_{\al})$.
Let $\Gm$ be an irreducible lattice in $G$, and let $\Lambda$ be a subgroup of $Comm_G(\Gm)$.
Suppose that $rank\,G:=\underset{\al\in A}{\sum}rank_{k_{\al}}\bold{G_{\al}}\geq 2$.

If the image of a homomorphism $\tau:\Lambda\to\bold{J}(L)$ is  Zariski dense in $\bold{J}$
and not
relatively compact in $\bold{J}(L)$, then there exists a unique $\al\in A$, 
a continuous homomorphism $\theta:k_{\al}\to L$ and a unique $\theta$-algebraic 
isomorphism $\eta:\bold{G_{\al}}\isom\bold{J}$ such that 
$\tau(\la)=\eta(\theta(\pr_{\al}(\la))$ for all $\la\in\Lambda$.
\end{Thm}

\begin{Emp} \label{E:super}
We use the notation of \re{Ma1} with $\Dt'=P\Gm\ii$. Take any $M$ and $S$ such 
that $rank\, G_{\bar{M}}\geq 2$. Then by \rp{condit}, $\Gm:=\Dt^S$ is an 
irreducible lattice in $G_{\bar{M}}$. We will try to apply \rt{super} in the 
following situation. Take $G=G_{\bar{M}}$, $\Lambda$ be the projection of 
$\Dt'$ to $G_{\bar{M}}$, $L=\B{C}$, $\bold{J}=(\PGL{d})_{\B{C}}$ and
 $\tau$ be the homomorphism $j_i:P\Gm\ii\to\PGL{d}(\B{C})$ for some 
$i\in\{1,...,g\}$. Consider first $i=1$. By \rp{Zar} and \rp{compact},
$\tau=j_1$ satisfies the conditions of \rt{super}, hence there exists an algebraic
 isomorphism $\eta_1:(\PGU{d-1,1})_{\B{C}}\isom (\PGL{d})_{\B{C}}$ such that 
$j_1(\gm)=\eta_1(\gm_{\be})$ for all $\gm\in P\Gm\ii$.

Now take $i\geq 2$. Suppose that $j_i(P\Gm\ii)$ is not relatively compact. Then 
using again \rp{Zar} we conclude from \rt{super} that there exists an algebraic isomorphism
$\eta_i:(\PGU{d-1,1})_{\B{C}}\isom (\PGL{d})_{\B{C}}$ such that 
$j_i(\gm)=\eta_i(\gm_{\be})$ for all $\gm\in P\Gm\ii$.
In particular, $j(P\Gm\ii)$ is not Zariski dense in $(\PGL{d})^g$.
This contradicts to \rp{Zar}. Therefore after a suitable conjugation we may 
assume that $j_i(P\Gm\ii)\subset \PGU{d}(\B{R})$ for all $i=2,...,g$.
\end{Emp}

It follows that up to an algebraic automorphism of $(\PGL{d})^g$, 
$j(P\Gm\ii)\subset\pgr\times\PGU{d}(\B{R})^{g-1}\cong
\bold{PH\ii}(\B{R})$ and that $j_1$ is the natural embedding 
$\bold{PG\ii}(F)_+\hra \bold{PG\ii}(F_{\be_1})$. Therefore $j$ together with the
natural embedding $\bold{PG\ii}(F)_+\hra\bold{PG\ii}(\afv)$ embed $P\Gm\ii$ into
$\pgr^0\times\PGU{d}(\B{R})^{g-1}\times\bold{PG\ii}(\afv)$.

\begin{Lem} \label{L:clpr}
The closure of the projection of $P\Gm\ii$ to $\PGU{d}(\B{R})^{g-1}\times 
\bold{PG\ii}(\afv)$ contains $\PGU{d}(\B{R})^{g-1}\times P(\bold{(G\ii)^{der}}
(\afv))$.
\end{Lem}

\begin{pf}
Let $(g_{\be},g_f)$ be an element of $\PGU{d}(\B{R})^{g-1}\times 
P(\bold{(G\ii)^{der}}(\afv))$, let $U\subset \PGU{d}(\B{R})^{g-1}$ be an open 
neighbourhood of $g_{\be}$, and let $S\in\ff{\bold{PG\ii}(\afv)}$. We have to show 
that $P\Gm\ii\cap(\pgr\times U\times g_f S)\neq\emp$. By the strong approximation 
theorem there exists a $\gm\in P\Gm\ii$, whose projection to $\bold{PG\ii}(\afv)$ 
belongs to $g_f S$. Let $\gm'$ be the projection of $\gm^{-1}$ to 
$\PGU{d}(\B{R})^{g-1}$.

Since $j(P\Gm\ii)$ belongs to the commensurator of $\DtS:=j(P\Gm\ii_S)$ in 
$\PGL{d}(\B{C})^g$, \rp{Zar}, \rp{inj} and \rl{comm}, b), c) imply that $\pr^1(\DtS)$ 
is dense in $\PGU{d}(\B{R})^{g-1}$. It follows that there exists $\dt\in P\Gm\ii$ 
whose projection to $\PGU{d}(\B{R})^{g-1}\times\bold{PG\ii}(\afv)$ belongs to 
$\gm'U\times S$. Then $\gm\dt$ belongs to $P\Gm\ii\cap(\pgr^0\times U\times g_fS)$.
\end{pf}

\begin{Emp} \label{E:comp}
Now we proceed as in the proof of \rt{complex}. Let $M$ and $S$ be as in \re{Ma1}, and
let $P\GmS\ii$ be the projection of $P\Gm\ii_{(S)}:=P\Gm\ii\cap 
\PGU{d}(\B{R})^{g-1}\times G_{\bar{M}}\times S$ to $\PGU{d}(\B{R})^{g-1}\times 
G_{\bar{M}}$. The proof of \rp{condit} holds in our case, hence $P\GmS\ii$ is 
arithmetic. 
It follows that there exists a permutation $\sigma$ of the set
$\{2,...,g\}$ such that for every $i=2,...,g$ there exists a unique algebraic 
isomorphism $r_i:\bold{PG_{F_{\be_i}}}\isom\PGU{d}$ satisfying
$r_i(\gm)=j_{\sigma(i)}(\gm)$ for each $\gm\in P\Gm\ii_{(S)}$. In particular,  
$\sigma$ and the $r_i$'s do not depend on $M$ and $S$.
Since $P\Gm\ii=\underset{M,S}{\cup} P\Gm\ii_{(S)}$, we then have 
$r_i(\gm)=j_{\sigma(i)}(\gm)$ for all $i\in\{2,...,g\}$ and $\gm\in P\Gm\ii$. This shows 
the existence of an algebraic isomorphism $\Phi_{\B{C}}$ wich will satisfy \rt{diag}, 
if we show that it is inner. But this can be immediately shown by the standard argument using elliptic elements 
and function $t$ defined in \re{fst3} (compare for example the proofs of 
\rp{homom} and \rp{inj}).
\end{Emp}

\subsection{Rationality question} \label{SS:Rat}
Consider the $(\bold{PH\ii}_{\B{C}})\an$-torsor $(\wt{P}\ii)\an\cong[B^{d-1}
\times(\bold{PH\ii}_{\B{C}})\an\times (E\ii/E\ii_0)^{disc}]/P\Gm\ii$
over $(\wt{X}\ii)\an$. As in the $p$-adic case, it has a canonical flat connection 
$\wt{\wt{\C{H}}}\ii$.
The same considerations as in the $p$-adic case (see the proof of \rp{decomp})
show that there exists a unique connection $\wt{\C{H}}\ii$ on $\wt{P}\ii$
such that  $(\wt{\C{H}}\ii)\an\cong\wt{\wt{\C{H}}}\ii$. It follows from the proofs of 
\rp{decomp} and \rt{diag} that $(\wt{f}_{\Phi,P})_*(\C{H}_{\B{C}})=\wt{\C{H}}\ii$.

\begin{Lem} \label{L:uncom}
If an analytic automorphism $\varphi:(\wt{P}\ii)\an\isom (\wt{P}\ii)\an$
commutes with the action of $(\bold{PH\ii}_{\B{C}})\an\times E\ii$, preserves 
$\wt{\C{H}}\ii$ and induces the identity map on $(\wt{X}\ii)\an=
(\PH_{\B{C}})\an\bs (\wt{P}\ii)\an$, then $\varphi$ is the identity.
\end{Lem} 

\begin{pf}
Recall that $(\wt{P}\ii)\an\cong[B^{d-1}\times(\bold{PH\ii}_{\B{C}})\an\times 
(E\ii/E\ii_0)^{disc}]/P\Gm\ii$. Since $\varphi$ induces the identity map on 
$(\wt{X}\ii)\an$, there exists a holomorphic map $\psi:B^{d-1}\to 
(\bold{PH\ii}_{\B{C}})\an$ such that $\varphi [x,1,1]=[x,\psi(x),1]$ for all 
$x\in B^{d-1}$. Since $\varphi$ preserves $\wt{\C{H}}\ii$, we have $\partial\psi/\partial x\equiv 0$. Hence $\psi$ is a constant, say $a$. Then $\varphi [x,h,e]=
[x,ha,e]$ for all $x\in B^{d-1},\;h\in(\bold{PH\ii}_{\B{C}})\an$ and $e\in E\ii$. 
In particular, $\varphi [x,1,1]=\varphi [\gm_{\be}^{-1}x,j(\gm),\gm_E]= 
[\gm_{\be}^{-1}x,j(\gm)a,\gm_E]=[x,j(\gm)aj(\gm)^{-1},1]$
for all $\gm\in P\Gm\ii$. Therefore $j(\gm)aj(\gm)^{-1}=a$ for all 
$\gm\in\Gm\ii$. Since $P\Gm\ii$ is Zariski dense in $\bold{PH\ii}_{\B{C}}$,
 $a=1$.
\end{pf} 

\begin {Cor} \label {C:canm}
The torsor $\wt{P}\ii$ has a unique $E\ii$-equivariant structure $P\ii$ of a 
$\bold{PH_{K_w}\ii}$-torsor over $X\ii$ such that there exists a connection 
$\C{H}\ii$ on $P\ii$ satisfying $\C{H}_{\B{C}}\ii\cong\wt{\C{H}}\ii$.
($P\ii$ is called the canonical model of $\wt{P}\ii$ over $X\ii$.)
\end{Cor} 

\begin {pf}
The existence is proved in \cite[III, $\S$3]{Mi}. Suppose that $P'$ and $P''$
are two structures satisfying the above conditions. Let
$f:P'_{\B{C}}\isom\wt{P}\ii\isom P''_{\B{C}}$ be the natural isomorphism. 
For each $\sigma\in\Aut_{K_w}(\B{C})$ set $\varphi_{\sigma}:={\sigma}(f)^{-1}
\circ f$. Then the automorphism $(\varphi_{\sigma})\an$ of $(P'_{\B{C}})\an\cong 
(\wt{P}\ii)\an$ satisfies the assumptions of the lemma. 
Hence $(\varphi_{\sigma})\an$ is the identity, so that ${\sigma}(f)=f$ for
all $\sigma\in\Aut_{K_w}(\B{C})$. It follows that $P'=P''$. 
\end {pf} 

To finish the proof of the Fourth (and the Third) Main Theorem it remains to show 
that the homomorphism $\wt{f}_{\Phi,P}:P\tw_{\B{C}}\isom P\ii_{\B{C}}$ is 
$K_w$-linear. 
Since $(\wt{f}_{\Phi,P})_*(\C{H}\tw _{\B{C}})=(\wt{f}_{\Phi,P})_*(\C{H}_{\B{C}})
=\C{H}\ii _{\B{C}}$, this follows from the lemma by the same considerations 
as the corollary.


\end{document}